%% file: Homology.tex
\documentclass[a4paper,12pt]{article}

\usepackage{amscd}
\usepackage{amsfonts}
\usepackage{amsmath}
\usepackage{amssymb}
\usepackage{amsthm}
\usepackage{ascmac}
\usepackage[T1]{fontenc}
\usepackage{here}
\usepackage{mathrsfs}
\usepackage{slashbox}
\usepackage{txfonts}
\usepackage[all]{xy}

\allowdisplaybreaks

\theoremstyle{plain}
\newtheorem{thm}{Theorem}[section]
\newtheorem{lmm}[thm]{Lemma}
\newtheorem{prp}[thm]{Proposition}
\newtheorem{crl}[thm]{Corollary}
\theoremstyle{definition}
\newtheorem{dfn}[thm]{Definition}
\newtheorem{rmk}[thm]{Remark}
\newtheorem{exm}[thm]{Example}

\newcommand{\vs}{\vspace{0.2in}}
\newcommand{\ba}{\begin{array*}}
\newcommand{\ea}{\end{array*}}
\newcommand{\be}{\begin{eqnarray*}}
\newcommand{\ee}{\end{eqnarray*}}
\newcommand{\bi}{\begin{itemize}}
\newcommand{\ei}{\end{itemize}}

\def\ens#1{{\mathchoice{\left\{ #1 \right\}}{\{ #1 \}}{\{ #1 \}}{\{ #1 \}}}}
\def\set#1#2{{\mathchoice{\left\{ #1 \middle| #2 \right\}}{\{ #1 \mid #2 \}}{\{ #1 \mid #2 \}}{\{ #1 \mid #2 \}}}}
\def\r#1{\text{\rm #1}}

\def\Bigv#1{\left| #1 \right|}
\def\v#1{{\mathchoice{\Bigv{#1}}{| #1 |}{| #1 |}{| #1 |}}}
\def\Bign#1{\left\| #1 \right\|}
\def\n#1{{\mathchoice{\Bign{#1}}{\| #1 \|}{\| #1 \|}{\| #1 \|}}}

\def\ol#1{\overline{#1}{}}
\def\tl#1{\tilde{#1}{}}
\def\ul#1{\underline{#1}{}}

\newcommand{\colim}{\operatorname*{colim}}

\newcommand{\bA}{\mathbb{A}}

\newcommand{\bN}{\mathbb{N}}

\newcommand{\bQ}{\mathbb{Q}}
\newcommand{\bR}{\mathbb{R}}

\newcommand{\bZ}{\mathbb{Z}}

\newcommand{\cA}{\mathscr{A}}

\newcommand{\cC}{\mathscr{C}}

\newcommand{\cM}{\mathscr{M}}
\newcommand{\cN}{\mathscr{N}}
\newcommand{\cO}{\mathscr{O}}

\newcommand{\cS}{\mathscr{S}}

\newcommand{\rH}{\r{H}}

\newcommand{\rZ}{\r{Z}}

\newcommand{\N}{\bN}
\newcommand{\Q}{\bQ}
\newcommand{\R}{\bR}
\newcommand{\Z}{\bZ}

\newcommand{\BdR}{\r{B}_{\r{dR}}}

\newcommand{\Gm}{\mathbb{G}_{\r{m}}}
\newcommand{\Qp}{\mathbb{Q}_p}

\newcommand{\Zp}{\mathbb{Z}_p}

\newcommand{\Ab}{\r{Ab}}
\newcommand{\alg}{\r{alg}}
\newcommand{\Alg}{\r{Alg}}

\newcommand{\can}{\r{can}}

\newcommand{\Ch}{\r{Ch}}

\newcommand{\Frac}{\r{Frac}}

\newcommand{\Gal}{\r{Gal}}

\newcommand{\Hom}{\r{Hom}}

\newcommand{\id}{\r{id}}
\newcommand{\im}{\r{im}}

\newcommand{\Int}{\r{Int}}
\newcommand{\inv}{\r{inv}}

\newcommand{\ob}{\r{ob}}
\newcommand{\op}{\r{op}}
\newcommand{\pr}{\r{pr}}

\newcommand{\sgn}{\r{sgn}}

\newcommand{\Spa}{\r{Spa}}

\newcommand{\sSet}{\r{sSet}}

\newcommand{\ZFC}{\r{ZFC}}

\newcommand{\Teichmuller}{Teichm\"uller }

\newcommand{\ad}{\r{ad}}
\newcommand{\Adic}{\r{Adic}}
\newcommand{\Conv}{\r{Conv}}
\newcommand{\ED}{\r{ED}}
\newcommand{\ev}{\r{ev}}
\newcommand{\Fil}{\r{Fil}}
\newcommand{\Germ}{\r{Germ}}
\newcommand{\ID}{\r{ID}}

\newcommand{\PreAdic}{\r{PreAdic}}

\title{Non-Archimedean Analytic Singular Homology Based On Cosimplicial Perfectoid Spaces and Integration along Cycles}
\author{Tomoki Mihara}
\date{}

\begin{document}

\maketitle
\input{Abstract}
\tableofcontents

\input{Introduction}
\input{Preliminaries}
\input{Group}
\input{Polytope}
\input{Cycle}
\input{Integration}

\input{References}

\end{document}

%% file: Abstract.tex
\begin{abstract}
We introduce singular homology for non-Archimedean analytic spaces using a cosimplicial perfectoid space as a Galois representation. We define an integration along a cycle, which gives a pairing with the singular homology and the space of differential forms.
\end{abstract}

%% file: Introduction.tex
\section{Introduction}
\label{Introduction}

In this paper, we invented a singular homology of an adic space over a non-Archimedean field $k$ as a Galois representation, and an integration of a differential form along a cycle as a $\BdR$-valued functional which might be helpful when one needs to check the non-triviality of a cycle. The formulation is based on a quite simple idea: Consider the simplicial homology associated to an explicit perfectoid cosimplicial object $\Delta_{K/k}^{\bullet}$ defined by the completed group algebra of the group $D^n_k$ of ``exponential maps on the standard $n$-simplex $\Delta^n$'' equipped with the supremum norm for each $n \in \N$, and define the integration of a differential form $\omega$ along a cycle $\gamma$ as the integration of the inverse image $\gamma^* \omega$ defined in a way relating an ``exponential map'' $\ul{q} \in \Hom_{\Ab}(\Z[p^{-1}],(k^{\alg})^{\times})$ on the standard $1$-simplex $\Delta^1$ given as a system of $p$-power roots of a $q \in k^{\times}$ to a certain $p$-adic period $\log \ul{q} \in \BdR$ as a non-Archimedean analogue of the Archimedean analytic formula ``$\int_{0}^{1} a^x dx = \frac{1}{\log a}$''.

\vs
In order to realise the natural idea in a rigorous way, we invent several methods to construct sheafy Banach algebras (Theorem \ref{sheafy group}, Theorem \ref{sheafy quotient}, and Theorem \ref{sheafy germ}). For example, we introduced a notion of a $p$-divisible normed group, and defined its completed group algebra as a sheafy Banach algebra (Theorem \ref{perfectoid}). For example, we introduce the notion of a $k$-exponential map on a non-Archimedean analogue of the notion of a convex subset of Euclidean spaces, and the set of $k$-exponential maps forms a $p$-divisible normed group.

\vs
Like the classical exponential map satisfies the exponential law, a $k$-exponential map satisfies a property related to a group homomorphism. In particular, an $(n+1)$-tuple of systems of $p$-power roots of non-zero elements of $k$ defines a $k$-exponential map on the non-Archimedean counterpart of the standard $n$-simplex $\Delta^n_{\Lambda}$ for any $n \in \N$. A $p$-power root $\ul{q}$ of a $q \in k^{\times}$ forms an element of $\varprojlim_{x \mapsto x^p} \cC \cong \cC^{\flat}$ (cf.\ \cite{Sch12} Lemma 3.4 (iii)), where $\cC$ denotes the completion of an algebraic closure $k^{\alg}$ of $k$ and $\cC^{\flat}$ the tilt of the perfectoid field $\cC$. In particular, it lifts to an element of the Witt ring $W(R)$ of the valuation ring $R$ of $\cC^{\flat}$ when $\v{q} \leq 1$, and hence to an element of $\BdR$. That is why $\log \ul{q}$ makes sense in $\BdR$ after an appropriate scaling into $1 + \Fil^1 \BdR$. For more details, see Example \ref{Tate curve}.

\vs
We summarise the construction of this paper. First, \S \ref{Preliminaries} consists of two subsections. In \ref{Convention and Terminology}, we introduce conventions which will help us to shorten notation and standard terminologies on norms. In \S \ref{Adic Ring}, we briefly recall adic rings and perfectoid algebras.

\vs
Secondly, \S \ref{Perfectoid Group Algebra} consists of two subsections. In \S \ref{Completed Group Algebra}, we introduce a notion of a submultiplicatively normed Abelian group, and its completed group algebra as a Banach algebra. In \S \ref{p-divisible normed Group}, we introduce a notion of a $p$-divisible normed group as a special submultiplicatively normed Abelian group, and show that its completed group algebra is sheafy.

\vs
Thirdly, \S \ref{Perfectoid Algebra Associated to Polytope} consists of two subsections. In \ref{Affine Germ}, we introduce a notion of an affine germ, and a notion of an affine polytope as a special affine germ and also as a non-Archimedean counterpart of a polytope in the usual sense. In particular, we introduce a non-Archimedean counterpart $\Delta^n_{\Lambda}$ of the standard $n$-simplex for an $n \in \N$. In \S \ref{Exponential Map}, we introduce a notion of a $k$-exponential map on an affine germ, and construct a $p$-divisible normed group using $k$-exponential maps.

\vs
Fourthly, \S \ref{Analytic Cycle} consists of two subsections. In \S \ref{Analytic Standard Simplex}, we study the completed group algebra $\cA_K(D^n_k)$ over an extension $K/k$ associated to the $p$-divisible normed group $D^n_k$ of $k$-exponential maps on $\Delta^n_{\Lambda}$ and its specific quotient $\cA_{K/k}^n$ for an $n \in \N$. In \S \ref{Analytic Singular Homology}, we introduce a notion of analytic singular homology using $\cA_{K/k}^n$.

\vs
Finally, \S \ref{p-adic Period and Integration} consists of three subsections. In \S \ref{Rigid Analytic Singular Simplex}, we introduce a notion of a rigid analytic singular $n$-simplex as a special morphism from the adic space associated to $\cA_{K/k}^n$. In \S \ref{Integration along Rigid Analytic Cycles}, we introduce an integration of a differential form along a cycle. In \S \ref{Appendix}, we collect several properties of the topological ring structure of period rings for the sake of the future study.

%% file: Preliminaries.tex
\section{Preliminaries}
\label{Preliminaries}

We introduce several notations, and briefly recall adic rings and perfectoid algebras.

\subsection{Convention and Terminology}
\label{Convention and Terminology}

We denote by $\N$ the set of natural numbers, which coincides with the least infinite ordinal $\omega$. Let $n \in \N$. Then $n$ coincides with the finite ordinal $n$ given as the set of natural numbers smaller than $n$, and the formula $i \in n$ is formally equivalent to ``$i$ is a natural number smaller than $n$''. We denote by $[n]$ the successor ordinal $n+1$ of $n$, which coincides with the set of natural numbers smaller than or equal to $n$. We apply these set-theoretic conventions in order to make convention and discussion on sums, products, and so on shorter. For convenience, we set $[-1] \coloneqq \emptyset$.

\vs
For a set $X$ and an $x \in X$, we denote by $\delta_{x/X} \colon X \to [1]$ the characteristic function of $\ens{x} \subset X$, and frequently regard its value as an element of an arbitrary ring $R$ through the natural map $[1] \hookrightarrow \Z \to R$.

\vs
For sets $X$ and $Y$, we denote by $X^Y$ the set of maps $Y \to X$. In particular, for a set $X$ and a natural number $n$, $X^n$ is a set of maps $x$ assigning an $x(i) \in X$ to each $i \in n$, i.e.\ natural number $i$ smaller than $n$, under the convention explained above. For a pair $X, Y$ of sets and a $y \in T$, we denote by $\pr_{X,y/Y}$ the $y$-th canonical projection $X^Y \to X$.

\vs
We denote by $\Ab$ the big category of Abelian groups and group homomorphisms. We deal with such big categories in a standard way to handle definable classes in $\ZFC$ set theory. Let $C$ be a category. We denote by $\ob(C)$ the class of objects of $C$. For $(X,Y) \in \ob(C)^2$, we denote by $\Hom_C(X,Y)$ the set of morphisms $X \to Y$ in $C$.

\vs
Whenever we refer to the supremum or the infimum, we intend those operations in $[0,\infty]$. In particular, we have $\sup \emptyset = 0$ and $\inf \emptyset = \infty$. When we refer to a valuation field, we always assume that the codomain of the valuation is $[0,\infty)$.

\vs
A {\it weighted set} is a set $X$ equipped with a map $\n{\bullet}_X \colon X \to [0,\infty)$ called {\it the weight of $X$}. Many notions appearing in this paper, e.g.\ a Banach space and a $p$-divisible normed group, naturally form weighted sets.

\vs
For a weighted set $X$ and a set $S$, we call the map
\be
\n{\bullet}_{\sup} \colon X^S \to [0,\infty], \ f \mapsto \sup_{s \in S} \n{f(s)}_X
\ee
{\it the supremum norm}. A map $f \in X \to S$ is said to be {\it bounded} (resp.\ {\it short}) if its supremum norm is finite (resp.\ smaller than or equal to $1$). 

\vs
For weighted sets $X$ and $Y$, we call the map
\be
\n{\bullet}_{\op} \colon X^Y \to [0,\infty], \ f \mapsto \inf \set{C \in [0,\infty]}{\forall y \in Y, \n{f(y)}_X \leq C \n{y}_Y}
\ee
{\it the operator norm}. A map $f \in X \to Y$ is said to be {\it bounded} (resp.\ short) if its operator norm is finite (resp.\ smaller than or equal to $1$), and is said to be {\it isometric} if $\n{f(x)}_Y = \n{x}$ for any $x \in X$. We note that when the weight of $Y$ is the constant map whose value is $1$, then the operator norm of a map $X \to Y$ coincides with its supremum norm. In this sense, the operator norm is a generalisation of the supremum norm.

\vs
Typical examples of a weighted set are a complete valuation field and a Banach space. For details on those notions, see \cite{BGR84} and \cite{Ber90}.

\subsection{Adic Ring}
\label{Adic Ring}

We recall a connection between seminormed rings and adic rings. Throughout this paper, a seminorm on a ring is assumed to be submultiplicative, but is not assumed to satisfy the equality $\n{1} = 1$. For a seminormed ring $R$, we denote by $R^{\circ} \subset R$ the subset of power-bounded elements, by $R_{\leq 1} \subset R^{\circ}$ the closed unit disc, by $R_{< 1} \subset R_{\leq 1}$ the open unit disc, $\Int(R_{\leq 1}) \subset R$ the integral closure of $R_{\leq 1}$ in $R$, and by $R_{\ad}$ the pair $(R,\Int(R_{\leq 1}))$.

\vs
Let $k$ be a valuation field with a non-discrete valuation. Then $k_{\ad}$ is the affinoid ring $(k,k_{\leq 1})$, and $A_{\ad}$ forms an affinoid $k$-algebra for any Banach $k$-algebra $A$ by \cite{Mih16} Proposition 1.13. We note that we assumed the completeness of $k$ in \cite{Mih16}, but the proof of \cite{Mih16} Proposition 1.13 did not use the assumption. We call $A_{\ad}$ {\it the affinoid ring over $k$ associated to $A$}.

\vs
Let $p$ denote the characteristic of the residue field $k_{\leq 1}/k_{< 1}$. We say that $k$ is a {\it perfectoid field} (cf.\ \cite{Sch12} Definition 1.2) if the Frobenius map on $k^{\circ}/p k^{\circ}$ is surjective. Suppose that $k$ is a perfectoid field. We recall three distinct notions of perfectoid.
\bi
\item[-] A Banach $k$-algebra $A$ is said to be a {\it perfectoid $k$-algebra} (cf.\ \cite{Sch12} Definition 1.6) if $A^{\circ}$ is bounded in $A$ and the Frobenius map on $A^{\circ}/p A^{\circ}$ is surjective.
\item[-] An affinoid ring $R = (R^{\triangleleft},R^{+})$ over $k_{\ad}$ is said to be a {\it perfectoid affinoid $k$-algebra} (cf.\ \cite{Sch12} p.\ 247) if $R^{\triangleleft}$ is the underlying topological ring of a perfectoid $k$-algebra and $R^{+}$ is open in $(R^{\triangleleft})^{\circ}$. 
\item[-]  A topological ring $R$ is said to be a {\it perfectoid ring} (cf.\ \cite{Fon13} \S 1.1) if $R^{\circ}$ is bounded in $R$ and there exists a topologically nilpotent unit $\varpi$ of $R$ satisfying $p \in \varpi^p R$ such that the Frobenius map $R^{\circ}/\varpi R^{\circ} \to R^{\circ}/\varpi^p R^{\circ}$ is surjective.
\ei
The following is probably well-known to experts, but we write down a full proof because of the lack of the knowledge of an appropriate reference.

\begin{prp}
\label{perfectoid Banach vs perfectoid affinoid}
Suppose that $k$ is a perfectoid field. For any Banach $k$-algebra $A$, the following are equivalent:
\begin{itemize}
\item[(i)] The spectral seminorm of $A$ is equivalent to $\n{\bullet}_A$, and $A$ forms a perfectoid $k$-algebra with respect to the spectral seminorm.
\item[(ii)] The affinoid ring $A_{\ad}$ over $k$ is a perfectoid affinoid $k$-algebra.
\item[(iii)] The underlying topological ring of $A$ is a perfectoid ring.
\end{itemize}
\end{prp}

\begin{proof}
We denote by $A'$ the seminormed $k$-algebra given as the underlying $k$-algebra of $A$ equipped with the spectral seminorm of $A$. First, assume (i). We show (ii) and (iii). Since $(A_{\ad})^{\triangleleft}$ is the underlying topological ring of $A$, and it coincides with that of the perfectoid $k$-algebra $A'$ by (i). Since $A_{\ad}^+ = \Int(A_{\leq 1})$ is a subring of $A^{\circ}$ containing the open subring $A_{\leq 1}$, it is also open in $((A_{\ad})^{\triangleleft})^{\circ} = A^{\circ}$. We conclude (ii).

\vs
By the non-discreteness of the valuation of $k$ and the surjectivity of the Frobenius map on $k^{\circ}/p k^{\circ}$, there exists a $\varpi \in k^{\times}$ with $\v{p} \leq \v{\varpi^p} < 1$. The surjectivity of the Frobenius map of $(A')^{\circ}/p (A')^{\circ}= A^{\circ}/p A^{\circ}$ implies that of the Frobenius map $A^{\circ}/\varpi A^{\circ} \to A^{\circ}/\varpi^p A^{\circ}$ by $p A^{\circ} \subset \varpi^p A^{\circ}$. We conclude (iii).

\vs
Next, assume (ii). Then there exists a norm $x$ on the underlying $k$-algebra of $A$ such that the normed $k$-algebra $A''$ given as the underlying $k$-algebra of $A$ equipped with $x$ is a perfectoid $k$-algebra and $(A_{\ad})^{\triangleleft} = A$ shares the underlying topological $k$-algebra with $A''$. Since the valuation of $k$ is non-trivial and the identity map $A \to A''$ is continuous, the identity map is bounded by \cite{BGR84} 2.1.8 Corollary 4 (ii). Since $A''$ is uniform, the boundedness of the identity map $A \to A''$ implies that of the identity map $A' \to A''$. Therefore, Banach's open mapping theorem (\cite{Bou53} Theorem I.3.3/1) implies that the identity map $A' \to A''$ is a homeomorphism. Since the valuation of $k$ is non-trivial and the identity map $A'' \to A'$ is continuous, it is bounded. Since $A'$ and $A''$ are uniform, the boundedness of the identities $A' \to A''$ and $A'' \to A'$ implies $A' = A''$, i.e.\ $\n{\bullet}_{A'} = \n{\bullet}_{A''}$. This implies that $\n{\bullet}_{A'}$ is equivalent to $\n{\bullet}_A$ and $A'$ is a perfectoid $k$-algebra. We conclude (i).

\vs
Finally, assume (iii). By \cite{Mih16} Proposition 1.17, the boundedness of $A^{\circ}$ implies that the spectral seminorm $\rho_A$ of $A$ is equivalent to $\n{\bullet}_A$. Take a topologically nilpotent unit $\varpi$ of $A$ satisfying $p \in \varpi^p A$ such that the Frobenius map $A^{\circ}/\varpi A^{\circ} \to A^{\circ}/\varpi^p A^{\circ}$ is surjective. We show the surjectivity of the Frobenius map of $(A')^{\circ}/p (A')^{\circ} = A^{\circ}/p A^{\circ}$. For this purpose, it suffices to show that $\inf_{g \in A^{\circ}} \rho_A(f - g^p) \leq \v{p}$ for any $f \in A^{\circ}$, i.e.\ $\inf_{g \in A^{\circ}} \rho_A(f - g^p) < \v{p} + \epsilon$ for any $(f,\epsilon) \in A^{\circ} + (0,\infty)$. Take an $N \in \N$ with $\v{\varpi}^{pN} < \v{p} + \epsilon$. By the inductive application of the surjectivity of the Frobenius map $A^{\circ}/\varpi A^{\circ} \to A^{\circ}/\varpi^p A^{\circ}$, there exists a $g \in (A^{\circ})^N$ such that $\rho_A(f - \sum_{i \in N} \varpi^{pi} g(i)^p) \leq \v{\varpi}^{pN}$. We obtain
\be
\rho_A \left( f - \left( \sum_{i \in N} \varpi^i g(i) \right)^p \right) \leq \max \ens{\v{p},\v{\varpi}^{pN}} < \v{p} + \epsilon,
\ee
and hence $\inf_{g \in A^{\circ}} \rho_A(f - g^p) < \v{p} + \epsilon$. We conclude (i).
\end{proof}

%% file: Group.tex
\section{Perfectoid Group Algebra}
\label{Perfectoid Group Algebra}

We construct a perfectoid algebra by introducing a completed group algebra of a $p$-divisible normed group.

\subsection{Completed Group Algebra}
\label{Completed Group Algebra}

Let $R$ be a non-Archimedean seminormed (resp.\ normed) ring with a multiplicative seminorm. We always assume a seminormed $R$-module to be non-Archimedean. Every seminormed $R$-module forms a weighted set, and the forgetful construction admits a left adjoint given by the free construction defined in the following way:

\vs
Let $X$ be a weighted set. We define $R^{\oplus X} \subset R^X$ as the $R$-submodule of maps $f \colon X \to R$ such that $f(x) = 0$ for all but finitely many $x \in X$, and equip it with the map
\be
\n{\bullet}_{R^{\oplus X}} \colon R^{\oplus X} \to [0,\infty), \ f \mapsto \sup_{x \in X} \v{f(x)}_R \n{x}_X,
\ee
for which $R^{\oplus X}$ forms a seminormed $R$-module. If $X$ is non-empty, the finiteness of non-trivial coefficients of an element of $R^{\oplus X}$ ensures that the supremum in the definition of the seminorm can be replaced by the maximum. For an $x \in X$, we abbreviate to $x$ the element of $R^{\oplus X}$ given as the ($R$-valued) characteristic function $\delta_{x/X}$ of $\ens{x} \subset X$ as long as there is no ambiguity. Then every $f \in R^{\oplus X}$ can be expressed as the essentially finite sum $\sum_{x \in X} f(x) x$.

\vs
We denote by $R^{\hat{\oplus} X}$ the completion of $R^{\oplus X}$. We recall the universality of the weighted completion $R^{\hat{\oplus} X}$ on short maps.

\begin{prp}
\label{universality of weighted sum}
Let $M$ be a Banach $R$-module, and $F \colon X \to M$ a short map. Then there exists a unique short $R$-linear homomorphism $R^{\hat{\oplus} F} \colon R^{\hat{\oplus} X} \to M$ such that for any $x \in X$, the equality $R^{\hat{\oplus} F}(x) = F(x)$ holds.
\end{prp}

\begin{proof}
By the universality of the direct sum $R^{\oplus X}$, there exists a unique $R$-linear homomorphism $R^{\oplus F} \colon R^{\oplus X} \to M$ such that for any $x \in X$, the equality $R^{\oplus F}(x) = F(x)$ holds. For any $f \in R^{\oplus X}$, we have
\be
\n{R^{\oplus F}(f)}_M = \n{R^{\oplus F} \left( \sum_{x \in X} f(x) x \right)}_M = \n{\sum_{x \in X} f(x) F(x)}_M \leq \sup_{x \in X} \v{f(x)} \ \n{F(x)}_M \leq \sup_{x \in X} \v{f(x)} \ \n{x}_X = \n{f}_{R^{\oplus X}}
\ee
by the strong triangular inequality of $\n{\bullet}_M$. This implies that $R^{\oplus F}$ is short. Therefore, it extends to a unique short $R$-linear homomorphism $R^{\hat{\oplus} F} \colon R^{\hat{\oplus} X} \to M$ by the universality of the completion.
\end{proof}

In the following in this section, suppose that $R$ is complete. For an $x \in X$, we also abbreviate to $x$ the image of $x \in R^{\oplus X}$ in $R^{\hat{\oplus} X}$. We put $X_+ \coloneqq \set{x \in X}{\n{x}_X > 0}$. For any $x \in X_+$, the evaluation map $R^{\oplus X} \to R$ at $x$ is short, and hence uniquely extends to a continuous $R$-linear homomorphism $R^{\hat{\oplus} X_+} \to R$ by the universality of the completion. Therefore, the map
\be
R^{\oplus X} \to R^{X_+}, \ f \mapsto (f(x))_{x \in X}
\ee
uniquely extends to a continuous $R$-linear homomorphism $R^{\hat{\oplus} X} \to R^{X_+}$ by the universality of the direct product $R^{X_+}$ with respect to the direct product topology. Since the resulting homomorphism is in fact injective, we regard $R^{\hat{\oplus} X}$ as an $R$-submodule of $R^{X_+}$ so that the convention $f(x)$ makes sense for any $(f,x) \in R^{\hat{\oplus} X} \times X_+$.

\begin{rmk}
\label{convergence along finite subsets}
Every $f \in R^{\hat{\oplus} X}$ can be expressed as the formal sum $\sum_{x \in X_+} f(x) x$ such that the net $(\sum_{x \in S} f(x) x)_S$ indexed by the set of finite subsets $S \subset X_+$ directed by inclusion converges to $f$.
\end{rmk}

A {\it submultiplicatively normed Abelian group} is an Abelian group $D$, to which we apply the multiplicative convention for convenience, equipped with a map $\n{\bullet}_D \colon D \to (0,\infty)$ such that for any $d \in D^2$, the inequality $\n{d(0)d(1)}_D \leq \n{d(0)}_D \n{d(1)}_D$ holds. The difference between the notion of a normed Abelian group and that of a submultiplicatively normed Abelian group is that the former one satisfies $\n{0} = 0$ while the latter one satisfies the positivity. A typical example of a submultiplicatively normed Abelian group is the multiplicative group of a non-zero normed ring.

\vs
Let $D$ be a submultiplicatively normed Abelian group. We equip the group algebra $R[D]$ with $\n{\bullet}_{R^{\oplus D}}$ under the formulation where the underlying $R$-module of the group algebra $R[D]$ is given as $R^{\oplus D}$, and denote the norm by $\n{\bullet}_{R[D]}$. Then $R[D]$ forms a normed $R$-algebra because of the submultiplicativity and the positivity of the norm of $D$.

\begin{dfn}
\label{completed group algebra}
We denote by $\cA_R(D)$ the completion of $R[D]$.
\end{dfn}

By the definition of $\n{\bullet}_{R[D]}$, the underlying Banach $R$-module of $\cA_R(D)$ coincides with $R^{\hat{\oplus} D}$. Therefore, the convention $f(d) \in R$ makes sense for each $(f,d) \in \cA_R(D) \times D$. Since the seminorm of $R[D]$ is a norm, we identify $R[D]$ with its image in $\cA_R(D)$. We show the universality of the completed group algebra $\cA_R(D)$.

\begin{prp}
\label{universality of completed group algebra}
Let $A$ be a Banach $R$-algebra, and $F \colon D \to A^{\times}$ a short group homomorphism. Then there exists a unique short $R$-algebra homomorphism $\cA_R(F) \colon \cA_R(D) \to A$ such that for any $d \in D$, the equality $\cA_R(F)(d) = F(d)$ holds.
\end{prp}

\begin{proof}
It suffices to show that the unique short $R$-linear extension $\cA_R(F) \colon \cA_R(D) \to A$ of $F$ given by Proposition \ref{universality of weighted sum} preserves the multiplication. We denote by $1_D$ the identity of $D$. Then we have $\cA_R(f)(1_D) = F(1_D) = 1 \in A^{\times}$ since $F$ is a group homomorphism.

\vs
We denote by $M \subset \cA_R(D)^2$ the $R$-submodule
\be
\set{r \in \cA_R(D)^2}{\cA_R(F)(r(0) r(1)) = \cA_R(F)(r(0)) \cA_R(F)(r(1))}.
\ee
By the continuity of the multiplications of $\cA_R(D)$ and $A$, $M$ is a closed subset with respect to the direct product topology. Since the composite of $\cA_R(F)$ and the canonical homomorphism $R[D] \to \cA_R(D)$ is the $R$-algebra homomorphism $R[D] \to A$ associated to $F$ by the universality of the group algebra, it is an $R$-algebra homomorphism. Therefore $M$ contains the image of $R[D]^2$. Since $R[D]$ is dense in $\cA_R(D)$, $M$ is dense in $\cA_R(D)^2$. Thus $M$ coincides with $\cA_R(D)^2$, i.e.\ $\cA_R(F)$ preserves the multiplication.
\end{proof}

\subsection{$p$-divisible normed Group}
\label{p-divisible normed Group}

Henceforth, let $p$ denote a prime number, and $k$ a complete valuation field with non-trivial valuation and residue characteristic $p$.

\vs
A {\it $p$-divisible normed group} is a submultiplicatively normed Abelian group $D$ satisfying the following:
\bi
\item[(1)] For any $d \in D$, there exists an $d_1 \in D$ such that $d_1^p = d$.
\item[(2)] For any $d \in D$, the equality $\n{d^p}_D = \n{d}_D^p$ holds.
\ei
The following examples are what we actually consider in this paper:

\begin{exm}
\label{p-divisible normed group}
\bi
\item[(1)] Let $K$ be a complete valuation field. If every element of $K$ has a $p$-th root in $K$, then $K^{\times}$ forms a $p$-divisible normed group with respect to the restriction of the valuation.
\item[(2)] Let $X$ be a set, $D$ be a $p$-divisible normed group, and $G \subset D^X$ a $p$-divisible subgroup. If $\sup_{x \in X} \n{g(x)}_D < \infty$ for any $g \in G$, then $G$ forms a $p$-divisible normed group with respect to the restriction of the supremum norm.
\ei
\end{exm}

Let $D$ be a $p$-divisible normed group. We observe the Banach $k$-algebra $\cA_k(D)$ constructed in Definition \ref{completed group algebra}.

\begin{thm}
\label{perfectoid}
The following hold:
\bi
\item[(1)] The norm $\n{\bullet}_{\cA_k(D)}$ is power-multiplicative.
\item[(2)] If $k$ is a perfectoid field, then $\cA_k(D)$ is a perfectoid $k$-algebra.
\ei
\end{thm}

\begin{proof}
We show the assertion (1). Let $f \in k[D]$. We show $\n{f^p}_{k[D]} = \n{f}_{k[D]}^p$. We may assume $f \neq 0$ by $\n{0^p}_{k[D]} = 0 = \n{0}_{k[D]}^p$. Then, we have $\n{f}_{k[D]} \neq 0$. Put $\tl{f} \coloneqq \sum_{d \in D} f(d)^p d^p$. Since $D$ is a $p$-divisible normed group, we have $\n{\tl{f}}_{k[D]} = \n{f}_{k[D]}^p \neq 0$. We obtain $\v{p} \n{\tl{f}}_{k[D]} < \n{\tl{f}}_{k[D]} = \n{f}_{k[D]}^p$. By the repetitive application of the inequality $\n{(g(0)+g(1))^p - g(0)^p - g(1)^p)}_{k[D]} \leq \v{p} (\max \ens{\n{g(0)}_{k[D]},\n{g(1)}_{k[D]}})^p$ for any $g \in k[D]^2$, we have $\n{\tl{f} - f^p}_{k[D]} \leq \v{p} \n{f}_{k[D]}^p$. By $\v{p} \n{f}_{k[D]}^p = \v{p} \n{\tl{f}}_{k[D]} < \n{\tl{f}}_{k[D]}$, we obtain $\n{f^p}_{k[D]} = \n{\tl{f} - (\tl{f}-f^p)}_{k[D]} = \n{\tl{f}}_{k[D]} = \n{f}_{k[D]}^p$. It implies that $\n{\bullet}_{k[D]}$ is power-multiplicative, and hence so is $\n{\bullet}_{\cA_k(D)}$ by the continuity of $\n{\bullet}_{\cA_k(D)}$ and the denseness of $k[D] \subset \cA_k(D)$.

\vs
We show the assertion (2). We obtain $\cA_k(D)^{\circ} = \cA_k(D)_{\leq 1}$ by the assertion (1) and \cite{Mih16} Corollary 1.11. In particular, $\cA_k(D)^{\circ}$ is an open bounded subset of $\cA_k(D)$. Therefore it suffices to show that the surjectivity of the Frobenius on $\cA_k(D)_{\leq 1}/p \cA_k(D)_{\leq 1}$. Let $\ol{f} \in \cA_k(D)_{\leq 1}/p \cA_k(D)_{\leq 1}$. Take an $f \in \cA(D)_{\leq 1}$ representing $\ol{f}$. For each $d \in D$ with $f(d) = 0$, put $g_d \coloneqq 0 \in k$. For each $d \in D$ with $f(d) \neq 0$, there is a $g_d \in k$ such that $\v{g_d^p - f(d)} < \v{f(d)}$ and $\v{g_d^p - f(d)} \leq \v{p} \n{d}_D^{-1}$ by the perfectoidness of $k$. For each $d \in D$, there is an $e_d \in D$ such that $e_d^p = d$ by the $p$-divisibility of $D$. For any $d \in D$, we have $\v{g_d^p} = \v{f(d)}$ by the choice of $g_d$, and hence $\v{g_d} \n{e_d}_D = \v{f}^{\frac{1}{p}} \n{d}_D^{\frac{1}{p}}$. Therefore, the formal sum $\sum_{d \in D} g_d e_d$ converges to a $g \in \cA_k(D)$ in the sense of Remark \ref{convergence along finite subsets}. We have $\n{g}_{\cA_k(D)} = \sup_{d \in D} \v{g_d} \n{e_d}_D = \sup_{d \in D} (\v{f(d)} \n{d}_D)^{\frac{1}{p}} = \n{f}_{\cA_k(D)}^{\frac{1}{p}} \leq 1$. Put $\ol{g} \coloneqq g + p \cA_k(D)_{\leq 1} \in \cA_k(D)_{\leq 1}/p \cA_k(D)_{\leq 1}$. For any $d \in D$, we have $\n{(g_d e_d)^p - f(d) d}_{\cA_k(D)} = \v{g_d^d - f(d)} \n{d}_D \leq \v{p}$. We obtain $\ol{g}^p = \sum_{d \in D} (g_d e_d)^p + p \cA_k(D)_{\leq 1} = \sum_{d \in D} f(d) d + p \cA_k(D)_{\leq 1} = \ol{f}$.
\end{proof}

Now we introduce an analytic space $D_k^{\vee}$, which will be used to give a non-trivial example of a cycle for analytic singular homology.

\begin{dfn}
\label{perfectoid dual group}
We denote by $D_k^{\vee}$ the affinoid pre-adic space $\Spa(\cA_k(D)_{\ad})$.
\end{dfn}

We recall that the affine scheme associated to the group algebra of an abstract Abelian group $G$ plays a role of the dual in the algebraic setting, because of the fact that the algebra contravariantly represents the group of characters $G \to \Gm$. Similarly, the affinoid pre-adic space $D_k^{\vee}$ associated to the completed group algebra $\cA_k(D)$ plays a role of the dual in the analytic setting, because of the fact that the Banach algebra contravariantly represents the group of short characters $D \to \Gm$ by Proposition \ref{universality of completed group algebra}. Moreover, the dual $D_k^{\vee}$ actually has a structure of a geometric object in the following sense:

\begin{thm}
\label{sheafy group}
\bi
\item[(1)] If $k$ is a perfectoid field, then $D_k^{\vee}$ is an affinoid perfectoid adic space over $\Spa(k_{\ad})$.
\item[(2)] If $k$ is spherically complete, then $D_k^{\vee}$ is an affinoid adic space over $\Spa(k_{\ad})$.
\ei
\end{thm}

In order to show Theorem \ref{sheafy group}, we prepare lemmata on a scalar extension and Tate's acyclicity.

\begin{lmm}
\label{scalar extension}
Let $F_1/F_0$ be an extension of valuation fields with non-trivial valuations, and $V$ a normed $F_0$-vector space. If $F_0$ is spherically complete, then the composite of the embedding $V \hookrightarrow F_1 \otimes_{F_0} V$ associated to the inclusion $F_0 \hookrightarrow F_1$ and the completion $\hat{\id}_{F_1 \otimes_{F_0} V} \colon F_1 \otimes_{F_0} V \to F_1 \hat{\otimes}_{F_0} V$ of $\id_{F_1 \otimes_{F_0} V}$ is injective.
\end{lmm}

We note that $F_1 \otimes_{F_0} V$ is a seminormed $F_0$-vector space. If the seminorm of $F_1 \otimes_{F_0} V$ were not a norm, then $\hat{\id}_{F_1 \otimes_{F_0} V}$ would not be injective. Therefore the injectivity of the composite is not obvious.

\begin{proof}
It suffices to show that the embedding $V \hookrightarrow F_1 \otimes_{F_0} V$ is isometric. By the definition of the tensor seminorm, it suffices to show that the inclusion $V \hookrightarrow W \otimes_{F_0} V$ is isometric for any finite dimensional $F_0$-vector subspace $W \subset F_1$ containing $F_0$. Since $F_0$ is spherically complete, $F_1$ is an $F_0$-Cartesian normed $F_0$-vector space by \cite{BGR84} 2.4.4 Proposition 2, and hence there exists an $F_0$-vector subspace $W_0 \subset W$ such that the addition $F_0 \oplus W_0 \to W$ is an isometric $F_0$-linear isomorphism. Therefore the addition $V \oplus (W_0 \otimes_{F_0} V) \to W \otimes_{F_0} V$ is an isometric $F_0$-linear isomorphism. It implies that the inclusion $V \hookrightarrow W \otimes_{F_0} V$ is isometric.
\end{proof}

An affinoid ring $R = (R^{\triangleleft},R^{+})$ is said to be {\it uniform} if $(R^{\triangleleft})^{\circ}$ is bounded in $R^{\triangleleft}$, and is said to be {\it stably uniform} (in the sense of \cite{BV18} pp.\ 30 -- 31, or {\it locally uniform} in the sense of \cite{Mih16} Definition 4.7) if for any rational domain $U \subset \Spa(R)$, the affinoid ring $(\cO_{\Spa(R)}(U),\cO_{\Spa(R)}^{+}(U))$ is uniform.

\begin{lmm}
\label{stably uniform}
Let $F_1/F_0$ be an extension of valuation fields with non-trivial valuations, and $A$ a Banach $F_0$-algebra. If $F_0$ is spherically complete and $(F_1 \hat{\otimes}_{F_0} A)_{\ad}$ is a stably uniform affinoid ring over $F_1$, then $A_{\ad}$ is a stably uniform affinoid ring over $F_0$.
\end{lmm}

\begin{proof}
Let $U$ be a rational domain of $\Spa(A_{\ad})$. Take an $n \in \N$, a $g \in A$, and an $f \in A^n$ such that $\sum_{i \in n} f(i) A = A$ and $U = \set{x \in \Spa(A_{\ad})}{\forall i \in n, \v{f(i)(x)} \leq \v{g(x)}}$. Put $U_{F_1} \coloneqq \set{x \in \Spa((F_1 \hat{\otimes}_{F_0} A)_{\ad})}{\forall i \in n, \v{(1 \otimes f(i))(x)} \leq \v{(1 \otimes g)(x)}}$. We denote by $T$ the $n$-tuple of indeterminates, and by $(gT -f)$ the closure of the ideal of $F_0 \ens{T}$ (resp.\ $F_1 \ens{T}$) generated by $\set{gT(i) - f(i)}{i \in n}$. We regard $A$ (resp.\ $F_1 \hat{\otimes}_{F_0} A$) as a Banach $F_0 \ens{T}$-algebra (resp.\ Banach $F_1 \ens{T}$-algebra) by the $F_0$-algebra (resp.\ $F_1$-algebra) homomorphism assigning $f(i)$ to $T(i)$ for each $i \in n$. Then, the natural $A$-algebra homomorphisms
\be
F_0 \ens{T}/(gT - f) \hat{\otimes}_{F_0 \ens{T}} A & \to & \cO_{\Spa(A_{\ad})}(U) \\
F_1 \ens{T}/(gT - f) \hat{\otimes}_{F_1 \ens{T}} (F_1 \hat{\otimes}_{F_0} A) & \to & \cO_{\Spa((F_1 \hat{\otimes}_{F_0} A)_{\ad})}(U_{F_1})
\ee
are homeomorphisms by \cite{Mih16} Proposition 2.16 applied to $U$ and $U_{F_1}$. Moreover, the natural $A$-algebra homomorphism
\be
\iota \colon F_0 \ens{T}/(gT - f) \hat{\otimes}_{F_0 \ens{T}} A \to F_1 \ens{T}/(gT - f) \hat{\otimes}_{F_1 \ens{T}} (F_1 \hat{\otimes}_{F_0} A)
\ee
is the composite of the natural $A$-algebra homomorphisms
\be
F_0 \ens{T}/(gT - f) \hat{\otimes}_{F_0 \ens{T}} A & \to & F_1 \hat{\otimes}_{F_0} (F_0 \ens{T}/(gT - f) \hat{\otimes}_{F_0 \ens{T}} A) \\
F_1 \hat{\otimes}_{F_0} (F_0 \ens{T}/(gT - f) \hat{\otimes}_{F_0 \ens{T}} A) & \to & F_1 \ens{T}/(gT - f) \hat{\otimes}_{F_1 \ens{T}} (F_1 \hat{\otimes}_{F_0} A).
\ee
The former one is an isometry by Lemma \ref{scalar extension}, and the latter one is an isometric isomorphism by the universality of push out (cf. \cite{BK20} Proposition 2.16). This implies that $\iota$ is an isometry. By the stable uniformity of $(F_1 \hat{\otimes}_{F_0} A)_{\ad}$, $(\cO_{\Spa((F_1 \hat{\otimes}_{F_0} A)_{\ad}}(U_{F_1}),\cO_{\Spa((F_1 \hat{\otimes}_{F_0} A)_{\ad}}^{+}(U_{F_1}))$ is uniform, and hence the norm of $F_1 \ens{T}/(gT - f) \hat{\otimes}_{F_1 \ens{T}} (F_1 \hat{\otimes}_{F_0} A)$ is equivalent to its spectral seminorm by \cite{Mih16} Proposition 2.6. In addition, since $\iota$ is an isometry, the norm of $F_0 \ens{T}/(gT - f) \hat{\otimes}_{F_0 \ens{T}} A$ is equivalent to its spectral seminorm, too. This implies that $(\cO_{\Spa(A_{\ad})}(U),\cO_{\Spa(A_{\ad})}^{+}(U))$ is uniform again by \cite{Mih16} Proposition 2.6. Thus, $\Spa(A_{\ad})$ is stably uniform.
\end{proof}

\begin{lmm}
\label{descent}
Let $F_1/F_0$ be an extension of complete valuation fields, and $X$ a weighted set. Then the short $F_1$-linear homomorphism $F_1 \hat{\otimes}_{F_0} F_0^{\hat{\oplus} X} \to F_1^{\hat{\oplus} X}$ induced by the inclusion $F_0^{\oplus X} \hookrightarrow F_1^{\oplus X}$ is an isometric $F_1$-linear isomorphism.
\end{lmm}

\begin{proof}
Replacing $X$ by $\set{x \in X}{\n{x}_X > 0}$, we may assume that the seminorms of $F_0^{\oplus X}$ and $F_1^{\oplus X}$ are norms. We denote by $V$ the completion of $F_1 \otimes_{F_0} F_0^{\oplus X}$. By \cite{BGR84} 2.1.7 Proposition 4, the composite of the $F_1$-linear homomorphisms
\be
F_1 \otimes_{F_0} \hat{\id}_{F_0^{\oplus X}} \colon F_1 \otimes_{F_0} F_0^{\oplus X} & \hookrightarrow & F_1 \otimes_{F_0} F_0^{\hat{\oplus} X} \\
\hat{\id}_{F_1 \otimes_{F_0} F_0^{\hat{\oplus} X}} \colon F_1 \otimes_{F_0} F_0^{\hat{\oplus} X} & \to & F_1 \hat{\otimes}_{F_0} F_0^{\hat{\oplus} X}
\ee
induces an isometric $F_1$-linear isomorphism $j \colon V \to F_1 \hat{\otimes}_{F_0} F_0^{\hat{\oplus} X}$. The short $F_1$-linear homomorphism in the assertion is the composite of the isometric $F_1$-linear isomorphism $j^{-1} \colon F_1 \hat{\otimes}_{F_0} F_0^{\hat{\oplus} X} \to V$ and the completion $V \to F_1^{\hat{\oplus} X}$ of the isometric $F_1$-linear isomorphism $F_1 \otimes_{F_0} F_0^{\oplus X} \to F_1^{\oplus X}$, and hence is an isometric $F_1$-linear isomorphism.
\end{proof}

\begin{lmm}
\label{stably uniform 2}
If $k$ is a perfectoid field or spherically complete, then $\cA_k(D)_{\ad}$ is a stably uniform affinoid ring over $k$.
\end{lmm}

\begin{proof}
Every perfectoid affinoid $F_1$-algebra is stably uniform by \cite{Sch12} Corollary 6.8. Therefore, the assertion for the first case follows from Theorem \ref{perfectoid} (2). Suppose that $k$ is spherically complete. Let $\cC$ denote the completion of an algebraic closure of $k$. By Lemma \ref{descent} applied to $(F_0,F_1,X) = (k,\cC,D)$, the short $\cC$-algebra homomorphism $\cC \hat{\otimes}_k \cA_k(D) \to \cA_{\cC}(D)$ induced by the inclusion $k[D] \hookrightarrow \cC[D]$ is an isometric $\cC$-algebra isomorphism. By the assertion for the first case, $\cA_{\cC}(D)_{\ad}$ is a stably uniform affinoid ring over $\cC$, and hence so is $(\cC \hat{\otimes}_k \cA_k(D))_{\ad}$. Therefore, $\cA_k(D)_{\ad}$ is a stably uniform affinoid ring over $k$ by Lemma \ref{stably uniform}.
\end{proof}

\begin{proof}[Proof of Theorem \ref{sheafy group}]
The assertion (1) immediately follows from Theorem \ref{perfectoid} (2) and \cite{Sch12} Theorem 6.3 (iii). The assertion (2) immediately follows from Lemma \ref{stably uniform 2} and \cite{BV18} Theorem 7 (or independently \cite{Mih16} Theorem 4.9).
\end{proof}

Let $B$ be a Banach $k$-algebra, and $\chi$ a short group homomorphism $D \to B^{\times}$. By Proposition \ref{universality of completed group algebra}, $\chi$ induces a short $k$-algebra homomorphism $\cA_k(\chi) \colon \cA_k(D) \to B$. We denote by $\cA_k(D)_{\chi}$ the closure of the image of $\cA_k(\chi)$, and by $V(\chi)$ the affinoid pre-adic space associated to $(\cA_k(D)_{\chi})_{\ad}$.

\begin{thm}
\label{sheafy quotient}
If $k$ is a perfectoid field and $B$ is a uniform Banach $k$-algebra, then $V(\chi)$ and $\Spa(\cA_k(D)/\ker(\cA_k(\chi)))$ are affinoid perfectoid adic spaces over $\Spa(k_{\ad})$.
\end{thm}

\begin{proof}
By Theorem \ref{perfectoid} (ii), $\cA_k(D)$ is a perfectoid $k$-algebra. By Proposition \ref{perfectoid Banach vs perfectoid affinoid} and \cite{KL19} Theorem 3.3.18 (iii) (b), $\cA_k(D)/\ker(\cA_k(\chi))$ is a perfectoid Banach $k$-algebra. By Proposition \ref{perfectoid Banach vs perfectoid affinoid} and \cite{KL19} Theorem 3.3.18 (ii) (b), $\cA_k(D)_{\chi}$ is a perfectoid Banach $k$-algebra. Thus the assertion follows from and \cite{Sch12} Theorem 6.3 (iii).
\end{proof}

%% file: Polytope.tex
\section{Perfectoid Algebra Associated to Polytope}
\label{Perfectoid Algebra Associated to Polytope}

\subsection{Affine Germ}
\label{Affine Germ}

In order to define a cosimplicial analytic space, we introduce a functorial way to construct an affinoid adic space associated to combinatorial data analogous to a polytope such as a standard simplex. Let $\Lambda \subset \R$ be a subring, e.g.\ $\Z[p^{-1}]$ or $\R$.

\begin{dfn}
A {\it $\Lambda$-affine germ} is a set $X$ equipped with an $n_X \in \N$ and an injective map $e_X \colon X \hookrightarrow \Lambda^{n_X}$ such that $\im(e_X)$ is bounded with respect to the Euclidean norm $\n{\bullet}_{\R^{n_X}}$.
\end{dfn}

We introduce the most important example of a $\Lambda$-affine polytope as a generalisation of a standard simplex, and will use it in the formulation of the cosimplicial analytic space.

\begin{exm}
\label{standard simplex}
Let $n \in \N \cup \ens{-1}$ and $\ell \in \Lambda \cap (0,\infty)$. We recall that $[n]$ is the set of all natural numbers smaller than or equal to $n$. We put
\be
\ell \Delta_{\Lambda}^n \coloneqq \set{\ell t}{t \in (\Lambda \cap [0,1])^{[n]} \land \sum_{i \in [n]} t(i) = 1},
\ee
and equip it with the inclusion $e_{\ell \Delta_{\Lambda}^n} \colon \ell \Delta_{\Lambda}^n \hookrightarrow \Lambda^{[n]}$. When $\ell = 1$, we abbreviate $\ell \Delta_{\Lambda}^n$ to $\Delta_{\Lambda}^n$. The $\Lambda$-affine germ $\ell \Delta_{\Lambda}^n$ will be used in the formulation of an analytic counterpart of the standard simplex in Definition \ref{analytic standard simplex}.
\end{exm}

For a $\Lambda$-affine germ $X$ and its subset $S$, we regard $S$ as a $\Lambda$-affine germ with respect to $e_S \coloneqq e_X |_S$ unless otherwise specified. We introduce several operations on a $\Lambda$-affine germ.

\begin{dfn}
Let $X$ be a $\Lambda$-affine germ $X$, and $A \subset \R$ a $\Lambda$-subalgebra. We denote by $A \otimes_{\Lambda} X$ the $A$-affine germ which shares the underlying set with $X$ and is equipped with the composite $e_{A \otimes_{\Lambda} X}$ of $e_X$ and the inclusion $\Lambda^{e_X} \hookrightarrow A^{e_X}$. We define
\be
\Conv_A(X) \coloneqq \set{\sum_{i \in [m]} t(i) x(i)}{m \in \N \land t \in \Delta_A^m \land x \in \im(e_X)},
\ee
and call it {\it the $A$-linear convex hull of $X$}. We denote by $\ol{\Conv}_{\Lambda}(X)$ the closure of $\Conv_{\Lambda}(X)$ in $\R^{n_X}$. We regard $\Conv_A(X)$ and $\ol{\Conv}_{\Lambda}(X)$ as an $A$-affine germ with respect to the inclusion $e_{\Conv_A(X)} \colon \Conv_A(X) \hookrightarrow A^{n_X}$ and an $\R$-affine germ with respect to the inclusion $e_{\ol{\Conv}_{\Lambda}(X)} \colon \ol{\Conv}_{\Lambda}(X) \hookrightarrow \R^{n_X}$ respectively.
\end{dfn}

Using the operations, we introduce the notions of a linearly convex $\Lambda$-affine germ and a $\Lambda$-affine polytope.

\begin{dfn}
Let $X$ be a $\Lambda$-affine germ.
\bi
\item[(1)] We say that a subset $S \subset X$ is a {\it convex generator of $X$} if it satisfies $\Conv_{\Lambda}(S) = \im(e_X)$. 
\item[(2)] We say that a convex generator $S$ of $X$ is {\it fine} if it is a finite set and the equality $\Lambda^{e_X} \cap \Conv_{\Frac(\Lambda)}(S) = \Conv_{\Lambda}(S)$ holds.
\item[(3)] We say that $X$ is {\it linearly convex} if $X$ is a convex generator of $X$.
\item[(4)] We say that $X$ is an {\it $\Lambda$-affine polytope} if it admits a fine convex generator of $X$.
\ei
\end{dfn}

For example, for any $\ell \in \Lambda \cap (0,\infty)$, $\ell \Delta_{\Lambda}^n$ is a $\Lambda$-affine polytope, because the canonical basis of $\Lambda^{[n]}$ rescaled by $\ell$ forms its fine convex generator. The linear convexity for the case where $\Lambda = \R$ and $e_X$ is an inclusion map is precisely the same as the classical convexity. The following criterion of the linear convexity is well-known for the classical convexity, and the typical proof works for our general setting:

\begin{prp}
\label{criterion of convexity}
For any $\Lambda$-affine germ $X$, the following are equivalent:
\bi
\item[(1)] The $\Lambda$-affine germ $X$ is linearly convex, i.e.\ the equality $\Conv_{\Lambda}(X) = \im(e_X)$ holds.
\item[(2)] There exists a convex generator of $X$, i.e.\ the equality $\Conv_{\Lambda}(S) = \im(e_X)$ holds for some subset $S \subset X$.
\ei
\end{prp}

In particular, we immediately obtain the following from Proposition \ref{criterion of convexity} and the definition of a $\Lambda$-affine polytope:

\begin{crl}
\label{convexity of polytope}
Every $\Lambda$-affine polytope is linearly convex.
\end{crl}

We formalise the notion of a linear inequality with coefficients in $\Lambda$ as the tuple of the coefficients itself so that the finiteness of a system of linear inequalities with coefficients in $\Lambda$ formally makes sense. For subsets $S_0$ and $S_1$ of $\R^n$ with $n \in \N$, a system $F$ of linear inequalities with coefficients in $\Lambda$ is said to {\it define $S_0$ in $S_1$} if $F$ is a system of $n$-ary relations and $S_0$ coincides with the set of solutions of $F$ in $S_1$.

\vs
When $\Lambda = \R$, then an $\R$-affine germ $X$ is a $\Lambda$-affine polytope if and only if there exists a finite system $F$ of linear inequalities with coefficients in $\R$ such that $F$ defines $\im(e_X)$ in $\Lambda^{n_X}$ by a classical result called Minkowski--Weyl Theorem. However, such a reformulation does not work for a general setting, as we see in the following counterexample:

\begin{exm}
We consider the case $\Lambda = \Z[p^{-1}]$. Then, the finite system
\be
\left\{
\begin{array}{rcl}
x & \geq & 0 \\
(p+1)x & \leq & 1
\end{array}
\right.
\ee
of linear inequalities with coefficients in $\Lambda$ defines a $\Lambda$-affine germ
\be
X \coloneqq \set{x \in \Lambda}{0 \leq x \land (p+1)x \leq 1}
\ee
with respect to the inclusion $e_X \colon X \hookrightarrow \Lambda$. The closure of $\im(e_X)$ in $\R$ coincides with $[0,(p+1)^{-1}]$, which is a closed subset of $\R$ defined by the same system of linear inequalities and forms an $\R$-affine polytope with respect to the inclusion $e_{[0,(p+1)^{-1}]}  \colon [0,(p+1)^{-1}] \hookrightarrow \R$. On the other hand, $X$ itself is not a $\Lambda$-affine polytope. Indeed, for any non-empty finite subset $S \subset X$, $\max(S)$ is the maximum of $\Conv_{\Lambda}(S)$, while $X$ itself does not have the maximum. In this way, a finite system of linear inequalities with coefficients in $\Lambda$ defines a $\Lambda$-affine germ $X$ which is not necessarily a $\Lambda$-affine polytope, even if it defines an $\R$-affine polytope which coincides with the closure of $\im(e_X)$.
\end{exm}

Despite of the counterexample above, it might look more natural to formulate a $\Lambda$-affine polytope in terms of a system of linear inequalities with coefficients in $\Lambda$. However, such a formulation is incompatible with the closure in $\R^{n_X}$ in the following sense:

\begin{exm}
We consider the case $\Lambda = \Z[p^{-1}]$. Then, the finite system
\be
\left\{
\begin{array}{rcl}
(p+1)x & \leq & 1 \\
1 & \leq & (p+1)x
\end{array}
\right.
\ee
of linear inequalities with coefficients in $\Lambda$ defines $\emptyset$ in $\Lambda$, and its closure $\emptyset$ in $\R$ does not coincide with the subset $\ens{(p+1)^{-1}} \subset \R$ defined by the same system of linear inequalities.
\end{exm}

In contrast to the pathological phenomena related to the formulation using a finite set of linear inequalities with coefficients in $\Lambda$, the next two propositions explain the main benefits of our formulation using the linear convex hull, which will be used in in the proof of Lemma \ref{polytope and linear inequality 2}, and will be used in the proof of Lemma \ref{convex hull 3}, which will be used in the proof of Theorem \ref{Weierstrass localisation} describing a specific Weierstrass localisation of the dual $D_k^{\vee}$ of a $p$-divisible normed group $D$ introduced in Definition \ref{perfectoid dual group}:

\begin{prp}[Compatibility of the expression by the convex hull and the closure]
\label{polytope and closure}
Suppose that $\Lambda$ is a dense subring of $\R$. Let $S$ be a $\Lambda$-affine germ whose underlying set is finite, e.g.\ a fine convex generator of a $\Lambda$-affine polytope. Then the equality $\Conv_{\R}(S) = \ol{\Conv}_{\Lambda}(S)$ holds.
\end{prp}

\begin{prp}[Compatibility of the expression by a finite set of linear inequalities and the closure]
\label{polytope and linear inequality}
Let $X$ be a $\Lambda$-affine polytope. Then there exists a finite system $F$ of linear inequalities with coefficients in $\Lambda$ such that the following hold:
\bi
\item[(1)] The system $F$ defines $\im(e_X)$ in $\Lambda^{n_X}$.
\item[(2)] For any system $F'$ of linear inequalities with coefficients in $\R$ such that $F'$ defines $\im(e_X)$ in $\Lambda^{n_X}$, if $F \subset F'$, then $F'$ also defines $\ol{\Conv}_{\Lambda}(X)$ in $\R^{n_X}$.
\ei
\end{prp}

\begin{proof}[Proof of Proposition \ref{polytope and closure}]
It suffices to show that $\Conv_{\Lambda}(S)$ is a dense subset of $\Conv_{\R}(S_{\R})$ and $\Conv_{\R}(S_{\R})$ is closed in $\R^{n_S}$. We denote by $c \in \N$ the cardinality of $S$, and fix a bijective map $s \colon c \to S$. For an $A \in \ens{\Lambda,\R}$, we put $T_A \coloneqq \set{t \in (A \cap [0,1])^{c-1}}{\sum_{i \in c-1} t(i) \leq 1}$.

\vs
Since $\Lambda$ is dense in $\R$, $T_{\Lambda}$ is dense in $T_{\R}$. Since $\Delta_A^{c-1}$ is the graph of the continuous map $T_A \to A, \ t \mapsto 1 - \sum_{i \in c-1} t(i)$ for any $A \in \ens{\Lambda,\R}$, $\Delta_{\Lambda}^{c-1}$ is dense in $\Delta_{\R}^{c-1}$. Since $\Conv_A(S)$ coincides with the image of the continuous map $\Delta_A^{c-1} \to A^{n_S}, \ t \mapsto \sum_{i \in c} t(i) e_S(s(i))$ for any $A \in \ens{\Lambda,\R}$, $\Conv_{\Lambda}(S)$ is dense in $\Conv_{\R}(S)$. Since $\Delta_{\R}^{c-1}$ is compact and $\R^{n_S}$ is Hausdorff, $\Conv_{\R}(S)$ is closed in $\R^{n_S}$.
\end{proof}

\begin{proof}[Proof of Proposition \ref{polytope and linear inequality}]
We take a fine convex generator $S$ of $X$, denote by $c \in \N$ the cardinality of $S$, and fix a bijective map $s \colon c \to S$. Then $\Conv_{\Frac(\Lambda)}(S)$ coincides with
\be
& & \set{y \in (\Frac(\Lambda) \cap [0,1])^{c}}{\sum_{i \in c} y(i) = 1} \\
& = & \set{x \in \Frac(\Lambda)^{n_X}}{\exists y \in \Frac(\Lambda)^{c}, \sum_{i \in c} y(i) = 1 \land \sum_{i \in c} y(i) s(i) = x \land \forall i \in c, y(i) \geq 0}.
\ee
Since $\Frac(\Lambda)$ is a field, Fourier--Motzkin elimination of $y$ applied to the finite system
\be
F_0 \coloneqq 
\left\{
\begin{array}{rcl}
\sum_{i \in c} y(i) & \leq & 1 \\
\sum_{i \in c} y(i) & \geq & 1 \\
\sum_{i \in c} y(i) s(i) & \leq & x \\
\sum_{i \in c} y(i) s(i) & \geq & x \\
\forall i \in c, y(i) & \geq & 0
\end{array}
\right.
\ee
of linear inequalities with coefficients in $\Lambda$ returns a finite system $F$ of linear inequalities with coefficients in $\Frac(\Lambda)$ such that it defines $\Conv_{\Frac(\Lambda)}(S)$ in $\Frac(\Lambda)^{n_X}$. Multiplying $F$ by a constant in $\Lambda \setminus \ens{0}$, we may assume that coefficients of $F$ belong to $\Lambda$. Since $S$ is fine, we have $\im(e_X) = \Lambda^{n_X} \cap \Conv_{\Frac(\Lambda)}(S)$. This implies that $F$ defines $\im(e_X)$ in $\Lambda^{n_X}$. 

\vs
Let $F'$ be a system of linear inequalities with coefficients in $\R$ such that $F'$ defines $\im(e_X)$ in $\Lambda^{n_X}$. Suppose $F \subset F'$. We denote by $Y \subset \R^{n_X}$ the subset defined by $F'$ in $\R^{n_X}$. We show $Y = \ol{\Conv}_{\Lambda}(X)$. By Proposition \ref{polytope and closure}, we have $\Conv_{\R}(S) = \ol{\Conv}_{\Lambda}(S)$. Therefore, it suffices to show $Y = \Conv_{\R}(S)$. Since $F'$ defines $Y$ in $\R^{n_X}$, $F'$ defines $\Lambda^{n_X} \cap Y$ in $\Lambda^{n_X}$. This implies $\im(e_X) = \Lambda^{n_X} \cap Y$. By $e_{\R \otimes_{\Lambda} X}(S) = e_X(S) \subset \im(e_X) = \Lambda^{n_X} \cap Y \subset Y$, we obtain $\Conv_{\R}(S) \subset Y$, because $Y$ is convex in the classical sense. Since we have the equivalence $F \Leftrightarrow \exists y \in A^c, F_0$ as predicates on $x \in A^{n_X}$ for any $\Lambda$-subalgebra $A \subset \R$ which is a field by the definition of $F$ using the quantifier elimination, $F$ defines $\Conv_{\R}(S)$ in $\R^{n_X}$. By $F \subset F'$, we obtain $Y \subset \Conv_{\R}(S)$. Thus, we obtain $Y = \Conv_{\R}(S)$.
\end{proof}

Let $X$ and $Y$ be $\Lambda$-affine germs, and $\varphi \colon X \to Y$ a map. We say that $\varphi$ is {\it $\Lambda$-affine} if there exists a pair $(v,\phi)$ of a $v \in \Lambda^{n_Y}$ and a $\Lambda$-linear homomorphism $\phi \colon \Lambda^{n_X} \to \Lambda^{n_Y}$ such that for any $x \in X$, the equality $\phi(e_X(x)) - e_Y(\varphi(x)) = v$ holds.

\vs
We say that $\varphi$ is {\it $\Lambda$-linear} if such a $v$ can be chosen to be zero. We note that when $\Lambda$ is $\Z[p^{-1}]$, then a map between $\Lambda$-modules is a $\Lambda$-linear homomorphism if and only if it is a group homomorphism, and hence we simply deal with group homomorphisms in that case. The notion of a $\Lambda$-affine map is useful when we introduce the dimension of a $\Lambda$-affine germ.

\subsection{Exponential Map}
\label{Exponential Map}

Henceforth, we only consider the case $\Lambda = \Z[p^{-1}]$, fix an algebraic closure $k^{\alg}$ of $k$, and denote by $\cC$ the completion of $k^{\alg}$. Let $X$ be a $\Lambda$-affine germ. For an $n \in \N$, we denote by $(\cC/k)^{n \flat \times}$ the Abelian group of group homomorphisms $\chi \colon \Lambda^n \to \cC^{\times}$ for which for any $i \in n$, there exists an $h \in \N$ with $\chi(p^h \delta_{i/n}) \in k^{\times}$, where the characteristic function $\delta_{i/n}$ of $\ens{i} \subset n$ is regarded as the $\Lambda$-valued one. As the convention shows, this is an analogue of $(\cC^n)^{\flat \times}$, i.e.\ the multiplicative group of the tilt of the perfectoid $\cC$-algebra $\cC^n$.

\vs
Let $K$ be a closed subextension of $\cC/k$. For a $\chi \in (\cC/k)^{n_X \flat \times}$, we denote by $\chi_X$ the map $X \to \cC$ assigning $\chi(e_X(x))$ to each $x \in X$. A map $f \colon X \to \cC$ is said to be {\it $k$-exponential} if there exists a $\chi \in (\cC/k)^{n \flat \times}$ with $f = \chi_X$, and is said to be {\it $(K,k)$-analytic} if for any $\epsilon \in (0,\infty)$, there exists a $K$-linear combination $g$ of $k$-exponential maps $X \to \cC$ such that $\n{f - g}_{\sup} < \epsilon$.

\begin{dfn}
\label{analytic convex hull}
We denote by $D_{X,k} \subset (\cC^X)^{\times}$ the subgroup of $k$-exponential maps, and by $A_{X,K/k} \subset \cC^X$ the $K$-subalgebra of $(K,k)$-analytic maps. We equip them with the restrictions $\n{\bullet}_{D_{X,k}}$ and $\n{\bullet}_{A_{X,K/k}}$ of the supremum norm on $X$. We denote by $X_{K/k}$ the affinoid pre-adic space $\Spa((A_{X,K/k})_{\ad})$ over $\Spa(K_{\ad})$.
\end{dfn}

We note that the supremum norm of a $k$-exponential map $X \to \cC$ is finite by the boundedness of $\im(e_X)$, and hence $\n{\bullet}_{D_{X,k}}$ and $\n{\bullet}_{A_{X,K/k}}$ are actually norms, i.e.\ their images are contained in $[0,\infty)$. In particular, $D_{X,k}$ forms a $p$-divisible normed group by Example \ref{p-divisible normed group} (1) and (2). The affinoid pre-adic space $X_{K/k}$ plays a role of a convex hull of $X$, and will be used in the formulation of analytic singular homology in \S \ref{Analytic Singular Homology}.

\begin{thm}
\label{sheafy germ}
If $K$ is a perfectoid field, then $(D_{X,k})_K^{\vee}$ (cf.\ Definition \ref{perfectoid dual group}) and $X_{K/k}$ are affinoid perfectoid adic spaces over $\Spa(K_{\ad})$.
\end{thm}

\begin{proof}
The assertion for $(D_{X,k})_K^{\vee}$ follows from Theorem \ref{sheafy group}. We denote by $B \subset \cC^X$ the $K$-subalgebra of bounded maps, and equip it with the supremum norm, for which $B$ is uniform. By the definition, $A_{X,K/k}$ is the closure of the image of the short $K$-algebra homomorphism $\cA_K(D_{X,k}) \to B$ (cf.\ Definition \ref{completed group algebra}) associated to the inclusion $\iota \colon D_{X,k} \hookrightarrow B^{\times}$ by the universality of the completed group algebra (cf.\ Proposition \ref{universality of completed group algebra}), and hence $X_{K/k}$ coincides with the affinoid perfectoid $V(\iota)$ in Theorem \ref{sheafy quotient}.
\end{proof}

Since we have a natural morphism $X_{K/k} \to (D_{X,k})_K^{\vee}$ associated to the canonical short $K$-algebra homomorphism $\cA_K(D_{X,k}) \to A_{X,K/k}$, a morphism from $(D_{X,k})_K^{\vee}$ to an adic space induces a morphism from $X_{K/k}$. Therefore, study of $(D_{X,k})_K^{\vee}$ is helpful when one needs a non-trivial example of a morphism from $X_{K/k}$. We note that the morphism $X_{K/k} \to (D_{X,k})_K^{\vee}$ is not necessarily an isomorphism, due to $\sum_{i=0}^{p-2} \mu^i = 0$ for any non-trivial $p-1$-torsion element $\mu \in k_{\leq 1}^{\times}$.

\vs
First, we show the functoriality of the assignments $X \mapsto (D_{X,k})_K^{\vee}$ and $X \mapsto X_{K/k}$. For a subring $L \subset \cC$, we denote by $L^{\times} D_{X,k} \subset (\cC^X)^{\times}$ the subgroup generated by $L^{\times} \cup D_{X,k}$.

\vs
Let $Y$ be a $\Lambda$-affine germs, and $\varphi \colon X \to Y$ a $\Lambda$-affine map. For a map $f \colon Y \to \cC$, we denote by $\varphi^* f \colon X \to \cC$ the composite $f \circ \varphi$.

\begin{prp}
\label{functoriality of the analytification}
\bi
\item[(1)] If $\varphi$ is a $\Lambda$-linear map, then $\varphi^* d \in D_{X,k}$ for any $d \in D_{Y,k}$, and $\varphi^* f \in A_{X,K/k}$ for any $f \in A_{Y,K/k}$.
\item[(2)] Let $L$ be a subring of $\cC$. For any $d \in D_{Y,k}$ with $\im(d) \subset L^{\times}$, $\varphi^* d \in L^{\times} D_{X,k}$.
\item[(3)] Let $Z$ be a $\Lambda$-affine germ, and $\psi$ a $\Lambda$-linear map $Y \to Z$. Then the equality $\varphi^* \psi^* f = (\psi \circ \varphi)^* f$ holds for any $f \in A_{Z,K/k}$.
\ei
\end{prp}

\begin{proof}
The assertion (3) follows from the definition. The assertion (1) for $f$ follows from that for $d$ since the composition of $\varphi$ is short. We show the assertions (1) for $d$ and (2). Since $d$ is exponential, there exists a $\phi_0 \in (\cC/k)^{n_Y \flat \times}$ such that for any $y \in Y$, the equality $d(y) = \phi_0(e_Y(y))$ holds. In both of (1) and (2), the assignment $d \mapsto \varphi^* d$ (resp.\ $f \mapsto \varphi^* f$) gives a short group homomorphism $D_{\varphi,k} \colon D_{Y,k} \to D_{X,k}$ (resp.\ $K$-algebra homomorphism $A_{\varphi,K/k} \colon A_{Y,K/k} \to A_{X,K/k}$). By the functoriality of the group algebra, $D_{\varphi,k}$ induces a $K$-algebra homomorphism $K[D_{\varphi,k}] \colon K[D_{Y,k}] \to K[D_{X,k}]$. Since $D_{\varphi,k}$ is short, $K[D_{\varphi,k}]$ induces a short $K$-algebra homomorphism $\cA_K(D_{\varphi,k}) \colon \cA_K(D_{Y,k}) \to \cA_K(D_{X,k})$ (cf.\ Definition \ref{completed group algebra}). We denote by $(D_{\varphi,k})_K^{\vee}$ (resp.\ $\varphi_{K/k}$) the morphism in $(V)$ associated to $\cA_K(D_{\varphi,k})$ (resp.\ $A_{\varphi,K/k}$).
\end{proof}

We denote by $\Germ_p$ the big category of $\Lambda$-affine germs and linear maps, by $(V)$ Huber's big category of triples in the sense of \cite{Hub96} 1.1.2, and by $\Adic \subset (V)$ the full subcategory of adic spaces. Now we are ready for formulating the functoriality of $(D_{X,k})_K^{\vee}$ and $X_{K/k}$ on $X$. We note that the naturality of the canonical morphism $X_{K/k} \to (D_{X,k})_K^{\vee}$ with respect to the functoriality immediately follows from the definition.

\begin{crl}
The assignment $X \mapsto (D_{X,k})_K^{\vee}$ (resp.\ $X_{K/k}$) forms a functor $\Germ_p \to (V)/\Spa(K_{\ad})$ with respect to the assignment $\varphi \mapsto (D_{\varphi,k})_K^{\vee}$ (resp.\ $\varphi_{K/k}$). If $K$ is a perfectoid field, then it forms a functor $\Germ_p \to \Adic/\Spa(K_{\ad})$.
\end{crl}

\begin{proof}
The first statement immediately follows from  Proposition \ref{functoriality of the analytification} (3). The second assignment follows from Theorem \ref{sheafy germ} (resp.\ Theorem \ref{sheafy quotient}).
\end{proof}

Let $X$ be a $\Lambda$-affine germ. For an $x \in X$, the evaluation map $D_{X,k} \to \cC^{\times}$ at $x$ is a short group homomorphism, and hence extends to a unique $K$-algebra homomorphism $\ev_{X,K/k,x} \colon \cA_K(D_{X,k}) \to \cC$ by Proposition \ref{universality of completed group algebra}.

\begin{prp}
\label{supremum of evaluation}
For any $f \in K^{\times} D_{X,k}$, the following equality holds:
\be
\n{f}_{\cA_K(D_{X,k})} = \sup_{x \in X} \v{\ev_{X,K/k,x}(f)}
\ee
\end{prp}

\begin{proof}
By $f \in K^{\times} D_{X,k}$, there exists a $(c,\chi) \in K^{\times} \times (\cC/k)^{n_X \flat \times}$ with $f = c \chi_X$. We have
\be
\n{f}_{\cA_K(D_{X,k})} = \n{c \chi_X}_{\cA_K(D_{X,k})} = \v{c} \ \n{\chi_X}_{D_{X,k}} = \v{c} \sup_{x \in X} \v{\chi_X(x)} = \sup_{x \in X} \v{c \chi_X(x))} = \sup_{x \in X} \v{\ev_{X,K/k,x}(f)}.
\ee
\end{proof}

For a (possibly infinite but mainly finite) subset $F \subset \cA_K(D_{X,k})$, we put
\be
X F & \coloneqq & \set{x \in X}{\forall f \in F, \v{\ev_{X,K/k,x}(f)} \leq 1} \\
(D_{X,k})_K^{\vee} F & \coloneqq & \set{x \in (D_{X,k})_K^{\vee}}{\forall f \in F, x(f) \leq 1}.
\ee
We regard $X F$ as a $\Lambda$-affine germ with respect to the restriction $e_{X F}$ of $e_X$ to $X F$. We denote by $K^{\times} D_{X,k}^1 \subset K^{\times} D_{X,k}$ the subset of maps $f \colon X \to \cC$ satisfying that there exists a tuple $(c,m,\theta,\phi,r,a,b)$ of a $c \in K^{\times}$, an $m \in \N$, a $\theta \in (\cC/k)^{1 \flat \times}$, a $\phi \in \Hom_{\Ab}(\Lambda^{n_X},\Lambda)$, an $r \in (0,1)$, an $a \in \Lambda$, and a $b \in \Lambda^m$ such that $\v{c} = r^a$, $\v{\theta(i)(1)} = r^{b(i)}$ for any $i \in m$, and $f(x) = c \prod_{i \in m} \theta(i)(\phi(i)(e_X(x)))$ for any $x \in X$. We only consider the case where the parameter set $F$ is contained in $K^{\times} D_{X,k}^1$, and hence in $\cA_K(D_{X,k})^{\times}$. Although a similar construction works for rational domains, we do not need to deal with it because we only consider parameters in $\cA_K(D_{X,k})^{\times}$.

\vs
Now we give a presentation of a Weierstrass domain of $(D_{X,k})_K^{\vee}$ for a specific case.

\begin{thm}
\label{Weierstrass localisation}
Let $F \subset K^{\times} D_{X,k}^1$ be a finite subset, and $n \in \N$ denote the cardinality of $F$. Assume that the following hold:
\bi
\item[(1)] There exists a bijective map $\iota \colon n \to F$ such that $X \set{\iota(i)}{i \in m}$ is a $\Lambda$-affine polytope for any $m \in n$.
\item[(2)] The complete valuation field $K$ is a perfectoid field or spherically complete.
\ei
Then the inclusion $X F \hookrightarrow X$ induces an isomorphism $(X F)_{K/k} \to X_{K/k} F$ in $\Adic$.
\end{thm}

In order to prove Theorem \ref{Weierstrass localisation}, we prepare several lemmata. Before that, we give an example of an application of Theorem \ref{Weierstrass localisation}. For a $(\theta,i) \in (\cC/k)^{1 \flat \times} \times n_X$, we denote by $\theta_{X,i} \in D_{X,k}$ the composite of $e_X$, the $i$-th canonical projection $\pr_{\Lambda,i/n_X} \colon \Lambda^{n_X} \twoheadrightarrow \Lambda$, and $\theta$.

\vs
Regarding $\Lambda$ itself as a $\Lambda$-affine germ with respect to the identity map $\Lambda \to \Lambda^1$, we fix an exponential map $\pi \colon \Lambda \to \cC$ with $\pi(1) \in k_{<1} \setminus \ens{0}$, which exists by $\Hom_{\Ab}(\Lambda,\cC^{\times}) = \Hom_{\Ab}(\Z[p^{-1}],\cC^{\times}) \cong (\cC^{\flat})^{\times}$ because $\cC$ is algebraically closed.

\begin{exm}
\label{Subdivision}
Let $n \in \N$, $m \in \N$, and $\ell \in \Lambda \cap (0,\infty)$. For an $a \colon [n] \to p^m$, we put
\be
Y \coloneqq \ell \Delta_{\Lambda}^n \cap \prod_{i \in [n]} [p^{-m} \ell a(i), p^{-m} \ell(a(i)+1)]
\ee
and regard it as a $\Lambda$-affine germ with respect to the inclusion $e_Y \colon Y \hookrightarrow \Lambda^{[n]}$. If $K$ is a perfectoid field or spherically complete , then the morphism $(D_{Y,k})_K^{\vee} \to (D_{\ell \Delta_{\Lambda}^n,k})_K^{\vee}$ associated to the inclusion $Y \hookrightarrow \ell \Delta_{\Lambda}^n$ satisfies the universality of a Weierstrass domain. Indeed, the finite subset $F \subset K^{\times} D_{\ell \Delta_{\Lambda}^n,k}^1$ given as
\be
\set{\pi(1)^{-\ell a(i)} \pi_{\ell \Delta_{\Lambda}^n,i}^{p^m}}{i \in [n]} \cup \set{\pi(1)^{\ell(a(i)+1)} \pi_{\ell \Delta_{\Lambda}^n,i}^{-p^m}}{i \in [n]}
\ee
satisfies $Y = \ell \Delta_{\Lambda}^n F$, and every bijective map $\iota \colon 2n+2 \to F$ satisfies the condition in Theorem \ref{Weierstrass localisation}.
\end{exm}

For a $\phi \in \Hom_{\Ab}(\Lambda^{n_X},\Lambda)$, we denote by $\phi_{\R} \in \Hom_{\Ab}(\R^{n_X},\R)$ the unique $\R$-linear extension of $\phi$. We interpret the linear convexity into the convexity with respect to $D_{X,k}$.

\begin{lmm}
\label{infimum and convext hull}
\bi
\item[(1)] For any $\chi \in (\cC/k)^{n_X \flat \times}$, the following equalities hold:
\be
\sup_{x \in X} \v{\chi_X(x)} & = & \sup_{x \in \Conv_{\Lambda}(X)} \v{\chi(x)} \\
\inf_{x \in X} \v{\chi_X(x)} & = & \inf_{x \in \Conv_{\Lambda}(X)} \v{\chi(x)}
\ee
\item[(2)] For any $\phi \in \Hom_{\Ab}(\Lambda^{n_X},\Lambda)$, the following equalities hold:
\be
\sup_{x \in X} \phi(e_X(x)) & = & \sup_{x \in \ol{\Conv}_{\Lambda}(X)} \phi_{\R}(x) \\
\inf_{x \in X} \phi(e_X(x)) & = & \inf_{x \in \ol{\Conv}_{\Lambda}(X)} \phi_{\R}(x)
\ee
\ei
\end{lmm}

\begin{proof}
(1) By
\be
\sup_{x \in X} \v{\chi_X(x)} & = & \sup_{x \in X} \v{\chi(e_X(x))} = \left( \inf_{x \in X} \chi^{-1}(e_X(x)) \right)^{-1} = \left( \inf_{x \in X} (\chi^{-1})_X(x) \right)^{-1} \\
\sup_{x \in \Conv_{\Lambda}(X)} \v{\chi(x)} & = & \left( \inf_{x \in \Conv_{\Lambda}(X)} \chi^{-1}(x) \right)^{-1},
\ee
it suffices to show the assertion for the infimum.

\vs
For any $x_0 \in \Conv_{\Lambda}(X)$, there exists an $m \in \N$, a $t \in (\Lambda \cap [0,1])^m$, and an $x' \in \im(e_X)^m$ such that $\sum_{i \in m} t(i) = 1$ and $\sum_{i \in m} t(i) x'(i) = x_0$ by the definition of $\Conv_{\Lambda}(X)$, and hence
\be
\v{\chi(x_0)} = \v{\chi \left( \sum_{i \in m} t(i) x'(i) \right)} = \prod_{i \in m} \v{\chi(x'(i))}^{t(i)} \geq \inf_{i \in m} \v{\chi(x'(i))} \geq \inf_{x \in \im(e_X)} \v{\chi(x)}.
\ee
This implies the desired equality.

\vs
(2) By
\be
\sup_{x \in X} \phi(e_X(x)) & = & - \inf_{x \in X} (- \phi)(e_X(x)) \\
\sup_{x \in \ol{\Conv}_{\Lambda}(X)} \phi_{\R}(x) & = & - \inf_{x \in \ol{\Conv}_{\Lambda}(X)} (- \phi)_{\R}(x),
\ee
it suffices to show the assertion for the infimum.

\vs
By (1) applied to $\pi \circ \phi$, we have
\be
& & \v{\pi(1)}^{\inf_{x \in X} \phi(e_X(x))} = \sup_{x \in X} \v{(\pi \circ \phi)(e_X(x))} = \sup_{x \in \Conv_{\Lambda}(X)} \v{(\pi \circ \phi)(x)} \\
& = & \v{\pi(1)}^{\inf_{x \in \Conv_{\Lambda}(X)} \phi(x)} = \v{\pi(1)}^{\inf_{x \in \Conv_{\Lambda}(X)} \phi_{\R}(x)}
\ee
and hence
\be
\inf_{x \in X} \phi(e_X(x)) = \inf_{x \in \Conv_{\Lambda}(X)} \phi_{\R}(x).
\ee
Since $\phi_{\R}$ is continuous and $\Conv_{\Lambda}(X)$ is dense in $\ol{\Conv}_{\Lambda}(X)$, we obtain the desired equality.
\end{proof}

We show that $\ol{\Conv}_{\Lambda}(X)$ is separated from a point in its complement by a group homomorphism $\Lambda^{n_X} \to \Lambda$.

\begin{lmm}
\label{linear convexity}
For any $(x_0,\theta) \in \R^{n_X} \times (\cC/k)^{1 \flat \times}$ with $\theta(1) \in k_{<1} \setminus \ens{0}$, $x_0$ belongs to $\R^{n_X} \setminus \ol{\Conv}_{\Lambda}(X)$ if and only if there exists a $\phi \in \Hom_{\Ab}(\Lambda^{n_X},\Lambda)$ satisfying the following:
\be
\n{(\theta \circ \phi)_X}_{D_{X,k}} < \v{\theta(1)}^{\phi_{\R}(x_0)}
\ee
\end{lmm}

\begin{proof}
Since $\Conv_{\Lambda}(X)$ is linearly convex and $\Conv_{\Lambda}(X)^2 \times (\Lambda \cap [0,1])$ is dense in $\ol{\Conv}_{\Lambda}(X)^2 \times [0,1]$, $\ol{\Conv}_{\Lambda}(X)$ is a convex subset of $\R^{n_X}$ in the classical sense by the continuity of the map
\be
(\R^{n_X})^2 \times [0,1] \to \R^{n_X}, \ (x_0,x_1,t) \mapsto t x_0 + (1 - t) x_1.
\ee
First, suppose $x_0 \in \ol{\Conv}_{\Lambda}(X)$. Let $\phi \in \Hom_{\Ab}(\Lambda^{n_X},\Lambda)$. For any $x \in X$, we have
\be
\v{(\theta \circ \phi)_X(x)} = \v{\theta(\phi(e_X(x)))} = \v{\theta(1)}^{\phi(e_X(x))} = \v{\theta(1)}^{\phi_{\R}(e_X(x))}.
\ee
This implies
\be
\n{(\theta \circ \phi)_X}_{D_{X,k}} = \sup_{x \in X} \v{(\theta \circ \phi)_X(x)} = \v{\theta(1)}^{\inf_{x \in X} \phi(e_X(x))} = \v{\theta(1)}^{\inf_{x \in \ol{\Conv}_{\Lambda}(X)} \phi_{\R}(x)} = \sup_{x \in \ol{\Conv}_{\Lambda}(X)} \v{\theta(1)}^{\phi_{\R}(x)}
\ee
by Lemma \ref{infimum and convext hull} (2) applied to $\phi$, and hence
\be
\v{\theta(1)}^{\phi_{\R}(x_0)} \leq \n{(\theta \circ \phi)_X}_{D_{X,k}}
\ee
by the assumption $x_0 \in \ol{\Conv}_{\Lambda}(X)$.

\vs
Next, suppose $x_0 \notin \ol{\Conv}_{\Lambda}(X)$. By the boundedness of $\im(e_X)$, $\Conv_{\Lambda}(X)$ is bounded, and hence $\ol{\Conv}_{\Lambda}(X)$ is compact. Therefore, there exists an $x_1 \in \ol{\Conv}_{\Lambda}(X)$ such that the equality
\be
\n{x_1 - x_0}_{\R^{n_X}} = \sup_{x \in \ol{\Conv}_{\Lambda}(X)} \n{x - x_0}_{\R^{n_X}}
\ee
holds. By $x_0 \notin \ol{\Conv}_{\Lambda}(X)$ and $x_1 \in \ol{\Conv}_{\Lambda}(X)$, we have $x_1 \neq x_0$. For an $x \in \R^{n_X}$, we denote by $\theta_x$ the group homomorphism
\be
\R^{n_X} \to \R, \ y \mapsto x \cdot y.
\ee
We show that for any $x \in \ol{\Conv}_{\Lambda}(X)$, there exists a $c \in [1,\infty)$ such that $(x - x_0) - c(x_1 - x_0) \in \ker(\theta_{x_1 - x_0})$.

\vs
By $\R^{n_X} = \R(x_1 - x_0) + \ker(\theta_{x_1 - x_0})$, there exists a $c \in \R$ such that $(x - x_0) - c(x_1 - x_0) \in \ker(\theta_{x_1 - x_0})$. It suffices to show $c \geq 1$. Assume $c < 1$. We put $y \coloneqq (x - x_0) - c(x_1 - x_0) \in \ker(\theta_{x_1 - x_0})$, and denote by $f$ the function
\be
\R \to \R, \ t \mapsto \n{(t x_1 + (1-t) x) - x_0}_{\R^{n_X}}.
\ee
By $f(1) = \n{x_1 - x_0}_{\R^{n_X}} \neq 0$, $f$ is differentiable at a neighbourhood of $1$. For any $t \in \R$, we have
\be
& & f(t) = \n{(t x_1 + (1-t) x) - x_0}_{\R^{n_X}} = \n{t(x_1 - x_0) + (1-t)(x - x_0)}_{\R^{n_X}} \\
& = & \n{t(x_1 - x_0) + (1-t)(c(x_1 - x_0) + y)}_{\R^{n_X}} = \n{(t + (1-t)c)(x_1 - x_0) + (1-t) y}_{\R^{n_X}} \\
& = & \left( (t + (1-t)c)^2 \n{x_1 - x_0}_{\R^{n_X}}^2 + (1-t)^2 \n{y}_{\R^{n_X}}^2 \right)^{\frac{1}{2}} \\
& = & \left( (-(1-t)(1-c)+1)^2 \n{x_1 - x_0}_{\R^{n_X}}^2 + (1-t)^2 \n{y}_{\R^{n_X}}^2 \right)^{\frac{1}{2}}
\ee
and hence
\be
& & \frac{df}{dt} = \frac{d}{dt} \left( (-(1-t)(1-c)+1)^2 \n{x_1 - x_0}_{\R^{n_X}}^2 + (1-t)^2 \n{y}_{\R^{n_X}}^2 \right)^{\frac{1}{2}} \\
& = & \frac{1}{2} \left( (-(1-t)(1-c)+1)^2 \n{x_1 - x_0}_{\R^{n_X}}^2 + (1-t)^2 \n{y}_{\R^{n_X}}^2 \right)^{- \frac{1}{2}} \\
& & \left( 2(1-c)(-(1-t)(1-c)+1) \n{x_1 - x_0}_{\R^{n_X}}^2 - 2(1-t) \n{y}_{\R^{n_X}}^2 \right)
\ee
at a neighbourhood of $1$. In particular, we obtain
\be
\frac{df}{dt}(1) = \frac{1}{2} \n{x_1 - x_0}_{\R^{n_X}}^{-1} \cdot 2(1-c) \n{x_1 - x_0}_{\R^{n_X}}^2 = (1-c) \n{x_1 - x_0}_{\R^{n_X}} > 0.
\ee
This implies $f(t_0) < f(1)$ for some $t_0 \in [0,1)$. By the convexity of $\ol{\Conv}_{\Lambda}(X)$, we have $t_0 x_1 + (1-t_0) x \in \ol{\Conv}_{\Lambda}(X)$, and hence
\be
f(t_0) = \n{(t_0 x_1 + (1-t_0) x) - x_0}_{\R^{n_X}} \geq \inf_{x' \in \ol{\Conv}_{\Lambda}(X)} \n{x' - x_0}_{\R^{n_X}} = \n{x_1 - x_0}_{\R^{n_X}} = f(1).
\ee
This contradicts $f(t_0) < f(1)$. Thus, we obtain $c \geq 1$.

\vs
We put $\phi_0 \coloneqq \theta_{x_1 - x_0} \in \Hom_{\Ab}(\R^{n_X},\R)$. For any $x \in \ol{\Conv}_{\Lambda}(X)$, there exists a $c \in [1,\infty)$ such that $(x - x_0) - c(x_1 - x_0) \in \ker(\theta_{x_1 - x_0})$ by the argument above, and hence
\be
\phi_0(x) & = & \theta_{x_1 - x_0}(x) = \theta_{x_1 - x_0}(x_0 + c(x_1 - x_0)) = \theta_{x_1 - x_0}(x_0) + c \n{x_1 - x_0}_{\R^{n_X}}^2 \\
& = & \phi_0(x_0) + c \n{x_1 - x_0}_{\R^{n_X}}^2 \geq \phi_0(x_0) + \n{x_1 - x_0}_{\R^{n_X}}^2
\ee
This implies
\be
\inf_{x \in \ol{\Conv}_{\Lambda}(X)} \phi_0(x) \leq \phi_0(x_0) + \n{x_1 - x_0}_{\R^{n_X}}^2 > \phi_0(x_0).
\ee
Since $\Lambda^{n_X}$ is dense in $\R^{n_X}$, there exists an there exists an $x_2 \in \Lambda^{n_X}$ such that
\be
\sup_{x \in \ol{\Conv}_{\Lambda}(X)} \theta_{x_2}(x) > \theta_{x_2}(x_0).
\ee
by the continuity of the maps
\be
\begin{array}{rclrcl}
\R^{n_X} & \to & \R, & x' & \mapsto & \theta_{x'}(x_0) \\
\R^{n_X} & \to & \R, & x' & \mapsto & \inf_{x \in \ol{\Conv}_{\Lambda}(X)} \theta_{x'}(x).
\end{array}
\ee
We denote by $\phi$ the group homomorphism
\be
\Lambda^{n_X} \to \Lambda, \ y \mapsto &x_2 \cdot y.
\ee
Then we have $\phi_{\R} = \theta_{x_2}$, and hence
\be
\inf_{x \in \ol{\Conv}_{\Lambda}(X)} \phi_{\R}(x) > \phi_{\R}(x_0).
\ee
We obtain
\be
& & \n{(\theta \circ \phi)_X}_{D_{X,k}} = \sup_{x \in X} \v{(\theta \circ \phi)_X(x)} = \sup_{x \in X} \v{\theta(1)}^{\phi(e_X(x))} \\
& = & \sup_{x \in \im(e_X)} \v{\theta(1)}^{\phi(x)} = \sup_{x \in \ol{\Conv}_{\Lambda}(X)} \v{\theta(1)}^{\phi_{\R}(x)} = \v{\theta(1)}^{\inf_{x \in \ol{\Conv}_{\Lambda}(X)} \phi_{\R}(x)}
\ee
by the continuity of the map
\be
\R^{n_X} \to [0,\infty), \ x \mapsto \v{\theta(1)}^{\phi_{\R}(x)},
\ee
and hence
\be
\n{(\theta \circ \phi)_X}_{D_{X,k}} < \v{\theta(1)}^{\phi_{\R}(x_0)}.
\ee
by $\inf_{x \in \ol{\Conv}_{\Lambda}(X)} \phi_{\R}(x) > \phi_{\R}(x_0)$.
\end{proof}

For a $\xi \in \cM_K(\cA_K(D_{X,k}))$, we put
\be
\Re_X(\xi) \coloneqq (\log_{\v{\pi(1)}} \v{\pi_{X,i}(\xi)})_{i \in n_X} \in \R^{n_X}.
\ee
We show that $\Re_X(\xi)$ does not depend on the choice of $\pi$.

\begin{lmm}
\label{independence of Re}
For any $(\xi,\theta) \in \cM_K(\cA_K(D_{X,k})) \times (\cC/k)^{1 \flat \times}$ with $\theta(1) \in k_{<1} \setminus \ens{0}$, the following equality holds:
\be
\Re_X(\xi) = (\log_{\v{\theta(1)}} \v{\theta_{X,i}(\xi)})_{i \in n_X}
\ee
\end{lmm}

\begin{proof}
It suffices to show that for any $(\epsilon,i) \in n_X \times (0,\infty)$, the inequality
\be
\v{\log_{\v{\theta(1)}} \v{\theta_{X,i}(\xi)} - \log_{\v{\pi(1)}} \v{\pi_{X,i}(\xi)}} < \epsilon
\ee
holds. If $X = \emptyset$, then we have $\theta_{X,i} = \pi_{X,i}$. Therefore, we may assume $X \neq \emptyset$. We denote by $m_i \in \R$ (resp.\ $M_i \in \R$) the infimum (resp.\ supremum) of the $i$-th entry of points in $X$. We put
\be
\rho \coloneqq \log_{\v{\pi(1)}} \v{\theta(1)} \in (0,\infty).
\ee
Since $\Lambda$ is dense in $\R$, there exists an $r \in \Lambda$ such that the inequality
\be
\v{r - \rho} \rho^{-1} \left( \max \ens{-m_i, M_i} + \v{- \log_{\v{\pi(1)}} \n{\pi_{X,i}}_{D_{X,k}}} \right) < \epsilon
\ee
holds. We put $r = p^{-d} m$ with $(d,m) \in \Z^2$, and denote by $\pi^r$ the group homomorphism
\be
\Lambda \to \cC^{\times}, \ a \mapsto \pi(ra).
\ee
Then, we have $(\pi^r)^{p^d} = \pi^m$ with respect to the usual action of $\Z$ on the multiplicative Abelian group $(\cC/k)^{1 \flat \times}$, and hence
\be
\v{\pi^r_{X,i}(\xi)} = \v{(\pi^r_{X,i})^{p^d}(\xi)}^{p^{-d}} = \v{(\pi^r)^{p^d}_{X,i}(\xi)}^{p^{-d}} = \v{\pi^m_{X,i}(\xi)}^{p^{-d}} = \v{\pi_{X,i}(\xi)}^{p^{-d} m} = \v{\pi_{X,i}(\xi)}^r
\ee
by the multiplicativity of $\xi$. For any $x \in X$, we have
\be
\v{\frac{\theta_{X,i}(x)}{\pi^r_{X,i}(x)}} = \v{\frac{\theta(x(i))}{\pi(r x(i))}} = \frac{\v{\theta(x(i))}}{\v{\pi(r x(i))}} = \frac{\v{\theta(1)}^{x(i)}}{\v{\pi(1)}^{r x(i)}} = \frac{\v{\pi(1)}^{\rho x(i)}}{\v{\pi(1)}^{r x(i)}} = \v{\pi(1)}^{(\rho - r) x(i)}
\ee
This implies
\be
& & \v{\log_{\v{\pi(1)}} \v{\frac{\theta_{X,i}}{\pi^r_{X,i}}(\xi)}} = \max \set{\sigma \log_{\v{\pi(1)}} \v{\frac{\theta_{X,i}}{\pi^r_{X,i}}(\xi)}}{\sigma \in \ens{-1,1}} \\
& = & - \log_{\v{\pi(1)}} \max \ens{\v{\frac{\theta_{X,i}}{\pi^r_{X,i}}(\xi)}, \v{\frac{\pi^r_{X,i}}{\theta_{X,i}}(\xi)}} \\
& \leq & - \log_{\v{\pi(1)}} \max \ens{\n{\frac{\theta_{X,i}}{\pi^r_{X,i}}}_{D_{X,k}}, \n{\frac{\pi^r_{X,i}}{\theta_{X,i}}}_{D_{X,k}}} \\
& = & - \log_{\v{\pi(1)}} \max \ens{\sup_{x \in X} \v{\frac{\theta_{X,i}(x)}{\pi^r_{X,i}(x)}}, \sup_{x \in X} \v{\frac{\pi^r_{X,i}(x)}{\theta_{X,i}(x)}}} \\
& = & - \log_{\v{\pi(1)}} \max \ens{\sup_{x \in X} \v{\pi(1)}^{(\rho - r) x(i)}, \sup_{x \in X} \v{\pi(1)}^{-(\rho - r) x(i)}} \\
& = & - \log_{\v{\pi(1)}} \max \ens{\v{\pi(1)}^{(\rho - r) m_i}, \v{\pi(1)}^{(\rho - r) M_i}, \v{\pi(1)}^{-(\rho - r) m_i}, \v{\pi(1)}^{-(\rho - r) M_i}} \\
& = & - \log_{\v{\pi(1)}} \v{\pi(1)}^{\min \ens{(\rho - r) m_i, (\rho - r) M_i}, -(\rho - r) m_i, -(\rho - r) M_i} \\
& = & - \v{\rho - r} \min \ens{m_i, -M_i} = \v{\rho - r} \max \ens{-m_i, M_i},
\ee
and hence
\be
& & \v{\log_{\v{\theta(1)}} \v{\theta_{X,i}(\xi)} - \log_{\v{\pi(1)}} \v{\pi_{X,i}(\xi)}} = \v{\log_{\v{\pi(1)}^{\rho}} \v{\left( \frac{\theta_{X,i}}{\pi^r_{X,i}} \pi^r_{X,i} \right)(\xi)} - \log_{\v{\pi(1)}} \v{\pi_{X,i}(\xi)}} \\
& = & \v{\log_{\v{\pi(1)}} \left( \v{\frac{\theta_{X,i}}{\pi^r_{X,i}}(\xi)} \v{\pi^r_{X,i}(\xi)} \right)^{\rho^{-1}} - \log_{\v{\pi(1)}} \v{\pi_{X,i}(\xi)}} \\
& = & \v{\log_{\v{\pi(1)}} \left( \v{\frac{\theta_{X,i}}{\pi^r_{X,i}}(\xi)}^{\rho^{-1}} \v{\pi_{X,i}(\xi)}^{\rho^{-1} r - 1} \right)} \\
& \leq & \v{\log_{\v{\pi(1)}} \v{\frac{\theta_{X,i}}{\pi^r_{X,i}}(\xi)}^{\rho^{-1}}} + \v{\log_{\v{\pi(1)}} \v{\pi_{X,i}(\xi)}^{\rho^{-1} r - 1}} \\
& = & \rho^{-1} \v{\log_{\v{\pi(1)}} \v{\frac{\theta_{X,i}}{\pi^r_{X,i}}(\xi)}} + \v{\rho^{-1} r - 1} \ \v{- \log_{\v{\pi(1)}} \v{\pi_{X,i}(\xi)}} \\
& \leq & \rho^{-1} \v{\rho - r} \max \ens{-m_i, M_i} + \v{\rho^{-1} r - 1} \ \v{- \log_{\v{\pi(1)}} \n{\pi_{X,i}}_{D_{X,k}}} \\
& = & \v{\rho - r} \rho^{-1} \left( \max \ens{-m_i, M_i} + \v{- \log_{\v{\pi(1)}} \n{\pi_{X,i}}_{D_{X,k}}} \right) \\
& < & \epsilon.
\ee
\end{proof}

By the definition of the topology of $\cM_K(\cA_K(D_{X,k}))$, the map
\be
\Re_X \colon \cM_K(\cA_K(D_{X,k})) \to \R^{n_X}, \ \xi \mapsto \Re_X(\xi)
\ee
is continuous. Although $\Re(\xi)$ does not possess the information of the evaluation of a non-trivial linear combination of elements in $D_{X,k}$ at $\xi$, we show that it possess the information of the evaluation of an element of $K^{\times} D_{X,k}$ at $\xi$.

\begin{lmm}
\label{evaluation and realisation}
Let $f \in K^{\times} D_{X,k}$ and $\xi \in \cM_K(\cA_K(D_{X,k}))$.
\bi
\item[(1)] For any $(\theta,\phi) \in (\cC/k)^{1 \flat \times} \times \Hom_{\Ab}(\Lambda^{n_X},\Lambda)$ with $f = (\theta \circ \phi)_X$, the following equality holds:
\be
\v{f(\xi)} = \v{\theta(1)}^{\phi_{\R}(\Re_X(\xi))}
\ee
\item[(2)] For any $x \in X$ with $\Re_{\R}(\xi) = e_X(x)$, the following equality holds:
\be
\v{f(\xi)} = \v{\ev_{X,K/k,x}(f)}
\ee
\ei
\end{lmm}

\begin{proof}
(1) For any $x \in X$, we have 
\be
f(x) = \theta(\phi(e_X(x))) = \prod_{i \in n_X} \theta_{X,i}^{\phi(\delta_{i/n_X})}(x),
\ee
and hence $f = \prod_{i \in n_X} \theta_{X,i}^{\phi(\delta_{i/n_X})}$. This implies
\be
\v{f(\xi)} = \v{\prod_{i \in n_X} \theta_{X,i}^{\phi(\delta_{i/n_X})}(\xi)} = \prod_{i \in n_X} \v{\theta_{X,i}(\xi)}^{\phi(\delta_{i/n_X})} = \v{\theta(1)}^{\sum_{i \in n_X} \phi(\delta_{i/n_X}) \log_{\theta(1)} \v{\theta_{X,i}(\xi)}}.
\ee
In addition, by Lemma \ref{independence of Re}, we have
\be
\v{\theta(1)}^{\sum_{i \in n_X} \phi(\delta_{i/n_X}) \log_{\theta(1)} \v{\theta_{X,i}(\xi)}} = \v{\theta(1)}^{\sum_{i \in n_X} \phi(\delta_{i/n_X}) \Re_X(\xi)(i)} = \v{\theta(1)}^{\phi_{\R} \left( \sum_{i \in n_X} \Re_x(\xi)(i) \delta_{i/n_X} \right)} = \v{\theta(1)}^{\phi_{\R}(\Re_X(\xi))}.
\ee
Therefore, we obtain
\be
\v{f(\xi)} = \v{\theta(1)}^{\sum_{i \in n_X} \phi(\delta_{i/n_X}) \log_{\theta(1)} \v{\theta_{X,i}(\xi)}} = \v{\theta(1)}^{\phi_{\R}(\Re_X(\xi))}.
\ee
(2) By $f \in K^{\times} D_{X,k}$, there exists a $(c,\chi) \in K^{\times} \times (\cC/k)^{n_X \flat \times}$ with $f = c \chi_X$ holds. For each $i \in n_X$, we denote by $\theta_i$ the group homomorphism
\be
\Lambda \to \cC^{\times}, a \mapsto \chi(a \delta_{i/n_X}),
\ee
and by $f_i \in D_{X,k}$ the map $X \to \cC$ assigning $\theta_i(\pr_{\Lambda,i/n_X}(e_X(x)))$ to each $x \in X$, where $\pr_{\Lambda,i/n_X}$ denotes the $i$-th canonical projection $\Lambda^{n_X} \to \Lambda$. By
\be
\chi = \prod_{i \in n_X} \theta_i \circ \pr_{\Lambda,i/n_X},
\ee
we have
\be
f(x) = c \chi_X(x) = c \chi(e_X(x)) = c \prod_{i \in n_X} \theta_i(\pr_{\Lambda,i/n_X}(e_X(x))) = c \prod_{i \in n_X} f_i(x)
\ee
for any $x \in x$, and hence $f = c \prod_{i \in n_X} f_i$. Therefore, by (1) applied to $(f_i,\theta_i,\pr_{\Lambda,i/n_X})$ for each $i \in n_X$, we have
\be
& & \v{f(\xi)} = \v{\left( c \prod_{i \in n_X} f_i \right)(\xi)} = \v{c} \prod_{i \in n_X} \v{f_i(\xi)} \\
& = & \v{c} \prod_{i \in n_X} \v{\theta_i(1)}^{(\pr_{\Lambda,i/n_X})_{\R}(\Re_X(\xi))} = \v{c} \prod_{i \in n_X} \v{\chi(\delta_{i/n_X})}^{\Re_X(\xi)(i)} = \v{c} \prod_{i \in n_X} \v{\chi(\Re_X(\xi)(i) \delta_{i/n_X})} \\
& = & \v{c} \v{\chi \left( \sum_{i \in n_X} \Re_X(\xi)(i) \delta_{i/n_X} \right)} = \v{c} \v{\chi(\Re_X(\xi))} = \v{c \chi(e_X(x))} = \v{\ev_{X,K/k,x}(f)}.
\ee
\end{proof}

\begin{lmm}
\label{evaluation and realisation 2}
Let $f \in K^{\times} D_{X,k}$, $m \in \N$, and $(c,\theta,\phi) \in K^{\times} \times ((\cC/k)^{1 \flat \times})^m \times \Hom_{\Ab}(\Lambda^{n_X},\Lambda)^m$ satisfying $f = c \prod_{i \in m} (\theta(i) \circ \phi(i))_X$.
\bi
\item[(1)] For any $\xi \in \cM_K(\cA_K(D_{X,k}))$, the following equality holds:
\be
\v{f(\xi)} = \v{c} \prod_{i \in m} \v{\theta(i)(1)}^{\phi(i)_{\R}(\Re_X(\xi))}
\ee
\item[(2)] For any $x \in X$, the following equality holds:
\be
\v{\ev_{X,K/k,x}(f)} = \v{c} \prod_{i \in m} \v{\theta(i)(1)}^{\phi(i)(e_X(x))}
\ee
\ei
\end{lmm}

\begin{proof}
The assertion (1) immediately follows from Lemma \ref{evaluation and realisation} (1) applied to $(\theta_i,\phi_i)$ for each $i \in m$, and the assertion (2) immediately follows from Lemma \ref{evaluation and realisation} (2) applied to $\theta_i \circ \phi_i$ for each $i \in m$ and the assertion (1).
\end{proof}

We describe $\im(\Re_X)$ in terms of linear convexity.

\begin{lmm}
\label{convex hull}
The image of $\Re_X$ coincides with $\ol{\Conv}_{\Lambda}(X)$.
\end{lmm}

\begin{proof}
Since $\cM_K(\cA_K(D_{X,k}))$ is compact by \cite{Ber90} 1.2.1 Theorem and $\R^{n_X}$ is Hausdorff, $\im(\Re_X)$ is closed. Therefore it suffices to show $\Conv_{\Lambda}(X) \subset \im(\Re_X) \subset \ol{\Conv}_{\Lambda}(X)$.
\vs
First, let $x_0 \in \Conv_{\Lambda}(X)$. in order to show $x_0 \in \im(\Re_X)$, we show that the group homomorphism
\be
e \colon (\cC/k)^{n_X \flat \times} \to \cC^{\times}, \ \chi \mapsto \chi(x)
\ee
factors through the map
\be
(\bullet)_X \colon (\cC/k)^{n_X \flat \times} \twoheadrightarrow D_{X,k}, \ \chi \mapsto \chi_X.
\ee
Let $\chi \in ((\cC/k)^{n_X \flat \times})^2$ with $\chi(0)_X = \chi(1)_X$. By $x_0 \in \Conv_{\Lambda}(X)$, there exist an $m \in \N$, a $t \in (\Lambda \cap [0,1])^m$, and an $x \in X^m$ such that $\sum_{i \in m} t(i) = 1$ and $\sum_{i \in m} t(i) e_X(x(i)) = x_0$. For any $j \in \ens{0,1}$, we have
\be
e(\chi(j)) = \chi(j)(x_0) = \chi(j) \left( \sum_{i \in m} t(i) e_X(x(i)) \right) = \sum_{i \in m} t(i) \chi(j)(e_X(x(i))) = \sum_{i \in m} t(i) \chi(j)_X(x(i)).
\ee
This implies $e(\chi(0)) = e(\chi(1))$. Therefore, $e$ factors through $r$, i.e.\ there exists a unique $e' \in \Hom_{\Ab}(D_{X,k},\cC^{\times})$ such that $e' \circ (\bullet)_X = e$.

\vs
For any $d \in D_{X,k}$, we have $\v{e'(d)} \leq \n{d}_{D_{X,k}}$. Indeed, for any $\chi \in (\cC/k)^{n_X \flat \times}$ with $\chi_X = d$, we have
\be
& & \v{e'(d)} = \v{e'(\chi_X)} = \v{e(\chi)} = \v{\chi(x_0)} \leq \sup_{x \in \Conv_{\Lambda}(X)} \v{\chi(x)} = \sup_{x \in X} \v{\chi(e_X(x))} \\
& = & \sup_{x \in X} \v{\chi_X(x)} = \sup_{x \in X} \v{d(x)} = \n{d}_{D_{X,k}} 
\ee
by Lemma \ref{infimum and convext hull} (1) applied to $\chi$. This implies that $e'$ extends to a unique short $K$-algebra homomorphism $f \colon \cA_K(D_{X,k}) \to \cC$. We denote by $\xi \in \cM_K(\cA_K(D_{X,k}))$ the bounded multiplicative seminorm on $\cA_K(D_{X,k})$ corresponding to the character $f$, i.e.\ the composite of $f$ and the valuation of $\cC$. Then we obtain
\be
& & \Re_X(\xi) = (\log_{\v{\pi(1)}} \v{\pi_{X,i}(\xi)})_{i \in n_X} = (\log_{\v{\pi(1)}} \v{f(\pi_{X,i})})_{i \in n_X} = (\log_{\v{\pi(1)}} \v{e'(\pi_{X,i})})_{i \in n_X} \\
& = & (\log_{\v{\pi(1)}} \v{\pi(x_0(i))})_{i \in n_X} = (\log_{\v{\pi(1)}} \v{\pi(1)}^{x_0(i)})_{i \in n_X} = (x_0(i))_{i \in n_X} = x_0,
\ee
and hence $x \in \im(\Re_X)$.

\vs
Next, let $x_0 \in \im(\Re_X)$. We show $x_0 \in \ol{\Conv}_{\Lambda}(X)$. Assume $x_0 \notin \ol{\Conv}_{\Lambda}(X)$. Then by Lemma \ref{linear convexity}, there exists a $\phi \in \Hom_{\Ab}(\Lambda^{n_X},\Lambda)$ with $\n{(\pi \circ \phi)_X}_{D_{X,k}} < \v{\pi(1)}^{\phi_{\R}(x_0)}$. By $x_0 \in \im(\Re_X)$, there exists a $\xi \in \cM_K(\cA_K(D_{X,k}))$ such that $\Re_X(\xi) = x_0$, i.e.\ $\log_{\pi(1)} \v{\pi_{X,i}(\xi)} = x_0(i)$ for any $i \in n_X$. By Lemma \ref{evaluation and realisation} (1) applied to $(\pi,\phi)$, we have
\be
\v{(\pi \circ \phi)_X(\xi)} = \v{\pi(1)}^{\phi_{\R}(\Re_X(\xi)} = \v{\pi(1)}^{\phi_{\R}(x_0)}.
\ee
We obtain
\be
\n{(\pi \circ \phi)_X}_{\cA_K(D_{X,k})} = \n{(\pi \circ \phi)_X}_{D_{X,k}} < \v{\pi(1)}^{\phi_{\R}(x_0)} = \v{(\pi \circ \phi)_X(\xi)},
\ee
but this contradicts that $\xi$ is a bounded multiplicative seminorm on $\cA_K(D_{X,k})$. We conclude $x_0 \in \ol{\Conv}_{\Lambda}(X)$.
\end{proof}

\begin{lmm}
\label{polytope and linear inequality 2}
Let $X$ be a $\Lambda$-affine polytope. For any $f \in K^{\times} D_{X,k}$, if $X \ens{f}$ is a $\Lambda$-affine polytope, then there exists a finite system $F$ of linear inequalities with coefficients in $\Lambda$ such that the following hold:
\bi
\item[(1)] The system $F$ defines $\im(e_{X \ens{f}})$ in $\Lambda^{n_X}$.
\item[(2)] For any system $F'$ of linear inequalities with coefficients in $\R$ such that $F'$ defines $\im(e_{X \ens{f}})$ in $\Lambda^{n_X}$, if $F \subset F'$, then $F'$ also defines $\ol{\Conv}_{\Lambda}(X \ens{f})$ in $\R^{n_X}$ and for any $\xi \in \cM_K(\cA_K(D_{X,k}))$, $\v{f(\xi)} \leq 1$ is equivalent to that $\Re_X(\xi)$ is a solution of $F'$.
\item[(3)] A subsystem of $F$ defines $\ol{\Conv}_{\Lambda}(X)$ in $\R^{n_X}$.
\ei
\end{lmm}

Henceforth in this subsection, for each $(\chi,i) \in (\cC/k)^{n_X \flat \times} \times n_X$, we denote by $\theta_{\chi,i}$ the group homomorphism
\be
\Lambda \to \cC^{\times}, \ a \mapsto \chi(a \delta_{i/n_X}),
\ee
and put $\rho_{\chi,i} \coloneqq \log_{\v{\pi(1)}} \v{\theta_{\chi,i}(1)} \in \R$.

\begin{proof}
For any $\chi \in (\cC/k)^{n_X \flat \times}$, we have
\be
\chi = \prod_{i \in n_X} (\theta_{\chi,i} \circ \pr_{\Lambda,i/n_X})_X.
\ee
For any $(\chi,\xi) \in (\cC/k)^{n_X \flat \times} \times \cM_K(\cA_K(D_{X,k}))$, we obtain
\be
\v{\chi_X(\xi)} = \v{1} \prod_{i \in n_X} \v{\theta_{\chi,i}(1)}^{(\pr_{\Lambda,i/n_X})_{\R}(\Re_X(\xi))} = \prod_{i \in n_X} \v{\theta_{\chi,i}(1)}^{\Re_X(\xi)(i)} = \prod_{i \in n_X} \v{\pi(1)}^{\rho_{\chi,i} \Re_X(\xi)(i)} = \v{\pi(1)}^{(\rho_{\chi,i})_{i \in n_X} \cdot \Re_X(\xi)}
\ee
by Lemma \ref{evaluation and realisation 2} (1) applied to $\chi$, $n_X$, and $(1,(\theta_{\chi,i})_{i \in n_X},(\pr_{\Lambda,i/n_X})_{i \in n_X})$. For any $(\chi,x) \in (\cC/k)^{n_X \flat \times} \times X$, we obtain
\be
\v{\chi_X(x)} = \v{1} \prod_{i \in n_X} \v{\theta_{\chi,i}(1)}^{\pr_{\Lambda,i/n_X}(e_X(x))} = \prod_{i \in n_X} \v{\theta_{\chi,i}(1)}^{e_X(x)(i)} = \prod_{i \in n_X} \v{\pi(1)}^{\rho_{\chi,i} e_X(x)(i)} = \v{\pi(1)}^{(\rho_{\chi,i})_{i \in n_X} \cdot e_X(x)}
\ee
by Lemma \ref{evaluation and realisation 2} (2) applied to the same data.

\vs
We put $f = c_f (\chi_f)_X$ with $(c_f,\chi_f) \in K^{\times} \times (\cC/k)^{n_X \flat \times}$ and $r_f \coloneqq \log_{\pi(1)} \v{c_f}$. For any $\xi \in \cM_K(\cA_K(D_{X,k}))$, we have
\be
\v{f(\xi)} & = & \v{c_f} \ \v{(\chi_f)_X(\xi)} = \v{\pi(1)}^{r_f} \v{\pi(1)}^{(\rho_{\chi_f,i})_{i \in n_X} \cdot \Re_X(\xi)} = \v{\pi(1)}^{(\rho_{\chi_f,i})_{i \in n_X} \cdot \Re_X(\xi) + r_f},
\ee
and hence $\v{f(\xi)} \leq 1$ is equivalent to $(\rho_{\chi_f,i})_{i \in n_X} \cdot \Re_X(\xi) + r_f \geq 0$. Similarly, for any $x \in X$, we have
\be
\v{\ev_{X,K/k,x}(f)} & = & \v{c_f} \ \v{(\chi_f)_X(x)} = \v{\pi(1)}^{r_f} \v{\pi(1)}^{(\rho_{\chi_f,i})_{i \in n_X} \cdot e_X(x)} = \v{\pi(1)}^{(\rho_{\chi_f,i})_{i \in n_X} \cdot e_X(x) + r_f},
\ee
and hence $\v{\ev_{X,K/k,x}(f)} \leq 1$ is equivalent to $(\rho_{\chi_f,i})_{i \in n_X} \cdot e_X(x) + r_f \geq 0$.

\vs
By Proposition \ref{polytope and linear inequality} applied to $X$, there exists a finite system $F_0$ of linear inequalities with coefficients in $\Lambda$ such that the following hold:
\bi
\item[(1)'] The system $F_0$ defines $\im(e_X)$ in $\Lambda^{n_X}$.
\item[(2)'] For any system $F'$ of linear inequalities with coefficients in $\R$ such that $F'$ defines $\im(e_X)$ in $\Lambda^{n_X}$, if $F_0 \subset F'$, then $F'$ also defines $\ol{\Conv}_{\Lambda}(X)$ in $\R^{n_X}$.
\ei
By the conditions (1)' and (2)', $F_0$ defines $\ol{\Conv}_{\Lambda}(X)$ in $\R^{n_X}$. We denote by $F$ the system of linear inequalities with coefficients in $\Lambda$ given as the conjunction of $F_0$ and the linear inequality $(\rho_{\chi_f,i})_{i \in n_X} \cdot x + r_f \geq 0$ on $x \in \R^{n_X}$ with coefficients in $\R$. By
\be
\im(e_{X \ens{f}}) & = & \set{e_{X \ens{f}}(x)}{x \in X \ens{f}} = \set{e_X(x)}{x \in X \land \v{\ev_{X,K/k,x}(f)} \leq 1} \\
& = & \set{e_X(x)}{x \in X \land (\rho_{\chi_f,i})_{i \in n_X} \cdot e_X(x) + r_f \geq 0} \\
& = & \set{x \in \im(e_X)}{(\rho_{\chi_f,i})_{i \in n_X} \cdot x + r_f \geq 0},
\ee
$F$ defines $\im(e_{X \ens{f}})$ in $\Lambda^{n_X}$. Let $F'$ be a system of linear inequalities with coefficients in $\R$ with $F \subset F'$ such that $F'$ defines $\im(e_{X \ens{f}})$ in $\Lambda^{n_X}$. By $F_0 \subset F \subset F'$ and the condition (2)', $F'$ also defines $\ol{\Conv}_{\Lambda}(X \ens{f})$ in $\R^{n_X}$. Therefore, it suffices to show that for any $\xi \in \cM_K(\cA_K(D_{X,k}))$, $\v{f(\xi)} \leq 1$ is equivalent to that $\Re_X(\xi)$ is a solution of $F'$.

\vs
First, suppose that $\xi$ is a solution of $F'$. Then, we have $(\rho_{\chi_f,i})_{i \in n_X} \cdot \Re_X(\xi) + r_f \geq 0$, and hence $\v{f(\xi)} \leq 1$. Secondly, suppose $\v{f(\xi)} \leq 1$. By Lemma \ref{convex hull}, we have $\Re_X(\xi) \in \ol{\Conv}_{\Lambda}(X)$. Since $F_0$ defines $\ol{\Conv}_{\Lambda}(X)$ by the conditions (1)' and (2)', $\Re_X(\xi)$ is a solution of $F_0$. By $\v{f(\xi)} \leq 1$, we have $(\rho_{\chi_f,i})_{i \in n_X} \cdot \Re_X(\xi) + r_f \geq 0$, i.e.\ $\Re_X(\xi)$ is a solution of $(\rho_{\chi_f,i})_{i \in n_X} \cdot x + r_f \geq 0$. Therefore, $\Re_X(\xi)$ is a solution of $F$. Since $F$ and $F'$ define the common subset $\ol{\Conv}_{\Lambda}(X \ens{f})$ in $\R^{n_X}$, $\Re_X(\xi)$ is a solution of $F'$.
\end{proof}

\begin{lmm}
\label{convex hull 3}
Let $f \in K^{\times} D_{X,k}^1$. If $X$ and $X \ens{f}$ are $\Lambda$-affine polytopes, then the image of $\set{\xi \in \cM_K(\cA_K(D_{X,k}))}{\v{f(\xi)} \leq 1}$ by $\Re_X$ coincides with the closure of $\im(e_{X \ens{f}})$.
\end{lmm}

\begin{proof}
By Corollary \ref{convexity of polytope} applied to $X \ens{f}$, $X \ens{f}$ is linearly convex, i.e.\ the equality $\im(e_{X \ens{f}}) = \Conv_{\Lambda}(X \ens{f})$ holds. By Lemma \ref{polytope and linear inequality 2}, there exists a finite system $F$ of linear inequalities with coefficients in $\Lambda$ such that the following hold:
\bi
\item[(1)] The system $F$ defines $\im(e_{X \ens{f}})$ in $\Lambda^{n_X}$.
\item[(2)] For any system $F'$ of linear inequalities with coefficients in $\R$ such that $F'$ defines $\im(e_{X \ens{f}})$ in $\Lambda^{n_X}$, if $F \subset F'$, then $F'$ also defines $\ol{\Conv}_{\Lambda}(X \ens{f})$ in $\R^{n_X}$ and for any $\xi \in \cM_K(\cA_K(D_{X,k}))$, $\v{f(\xi)} \leq 1$ is equivalent to that $\Re_X(\xi)$ is a solution of $F'$.
\item[(3)] A subsystem of $F$ defines $\ol{\Conv}_{\Lambda}(X)$ in $\R^{n_X}$.
\ei
We denote by $Y$ the image of $\set{\xi \in \cM_K(\cA_K(D_{X,k}))}{\v{f(\xi)} \leq 1}$ by $\Re_X$. By the conditions (1) and (2), $F$ defines $\ol{\Conv}_{\Lambda}(X \ens{f})$ in $\R^{n_X}$ and $Y$ in $\im(\Re_X)$. By Lemma \ref{convex hull}, $\im(\Re_X)$ coincides with $\ol{\Conv}_{\Lambda}(X)$. By the condition (3), a subsystem of $F$ defines $\ol{\Conv}_{\Lambda}(X)$ in $\R^{n_X}$. Therefore, $F$ defines the common subset in $\R^{n_X}$ and $\im(\Re_X)$, i.e.\ $\ol{\Conv}_{\Lambda}(X \ens{f})$ coincides with $Y$.
\end{proof}

\begin{proof}[Proof of Theorem \ref{Weierstrass localisation}]
If $n = 0$, then the assertion is trivial by $X F = X$. Suppose $n > 0$. By $X F = X \set{\iota(i)}{i \in n-1} \ens{\iota(n-1)}$, we may reduce the assertion to the case $n = 1$, i.e.\ $F = \ens{f}$ for an $f \in K^{\times} D_{X,k}^1$. Since $X$ and $X \ens{f}$ are $\Lambda$-affine polytopes, they are linearly convex by Corollary \ref{convexity of polytope}.

\vs
We denote by $B$ the affinoid ring given as the Weierstrass localisation of $\cA_K(D_{X,k})_{\ad}$ corresponding to the Weierstrass domain $(D_{X,k})_K^{\vee} \ens{f} \hookrightarrow (D_{X,k})_K^{\vee}$, and by $\Psi$ the canonical morphism $\cA_K(D_{X,k})_{\ad} \to B$ of affinoid rings over $K_{\ad}$. It suffices to show that the short $K$-algebra homomorphism $\cA_K(D_{X,k}) \to \cA_K(D_{X \ens{f},k})$ associated to the inclusion $X \ens{f} \hookrightarrow X$ by the functoriality of $\cA_K$ induces an isomorphism $\Phi \colon B \to \cA_K(D_{X \ens{f},k})_{\ad}$ of affinoid rings over $K_{\ad}$.

\vs
We denote by $f' \in \cA_K(D_{X \ens{f},k})$ the image of $f$. By $f \in K^{\times} D_{X,k}^1 \subset K^{\times} D_{X,k}$, there exists a $(c_f,\chi_f) \in K^{\times} \times \Hom_{\Ab}(\Lambda^{n_X},\cC)$ with $f = c_f (\chi_f)_X$. We have $f' = c_f (\chi_f)_{X \ens{f}}$. We have $\chi_f = \prod_{i \in n_X} \theta_{\chi_f,i} \circ \pr_{\Lambda,i/n_X}$ and $f' = c_f \prod_{i \in n_X} (\theta_{\chi_f,i} \circ \pr_{\Lambda,i/n_X})_{X \ens{f}}$. For any $\xi \in \cM_K(\cA_K(D_{X \ens{f},k}))$, by Lemma \ref{evaluation and realisation 2} (1) applied to $((\theta_{\chi_f,i})_{i \in n_X},(\pr_{\Lambda,i/n_X})_{i \in n_X})$, we have
\be
\v{f'(\xi)} = \v{c_f} \prod_{i \in n_X} \v{\theta_{\chi_f,i}(1)}^{(\pr_{\Lambda,i/n_X})_{\R}(\Re_X(\xi))} = \v{c_f} \prod_{i \in n_X} \v{\theta_{\chi_f,i}(1)}^{\Re_X(\xi)(i)}.
\ee
Similarly, for any $x \in X \ens{f}$, by Lemma \ref{evaluation and realisation 2} (2) applied to the same data, we have
\be
\v{\ev_{X \ens{f},k,x}(f')} = \v{c_f} \prod_{i \in n_X} \v{\theta_{\chi_f,i}(1)}^{(\pr_{\Lambda,i/n_X})_{\R}(e_X(x))} = \v{c_f} \prod_{i \in n_X} \v{\theta_{\chi_f,i}(1)}^{e_X(x)(i)}.
\ee
By Lemma \ref{convex hull} applied to $X \ens{f}$, the linear convexity of $X \ens{f}$, and the continuity of the map
\be
\R^{n_X} \to (0,\infty), \ x \to \v{c_f} \prod_{i \in n_X} \v{\theta_{\chi,i}(1)}^{x(i)},
\ee
we obtain
\be
& & \sup_{\xi \in \cM_K(\cA_K(D_{X \ens{f},k}))} \v{f'(\xi)} = \sup_{\xi \in \cM_K(\cA_K(D_{X \ens{f},k}))} \v{c_f} \prod_{i \in n_X} \v{\theta_{\chi_f,i}(1)}^{\Re_X(\xi)(i)} = \sup_{x \in \im(\Re_{X \ens{f}})} \v{c_f} \prod_{i \in n_X} \v{\theta_{\chi_f,i}(1)}^{x(i)} \\
& = & \sup_{x \in X \ens{f}} \v{c_f} \prod_{i \in n_X} \v{\theta_{\chi_f,i}(1)}^{e_X(x)(i)} = \sup_{x \in X \ens{f}} \v{\ev_{X \ens{f},k,x}(f')} = \sup_{x \in X \ens{f}} \v{\ev_{X,K/k,x}(f)} \leq 1.
\ee
By Theorem \ref{perfectoid}, $\n{\bullet}_{\cA_K(D_{X \ens{f},k})}$ is power-multiplicative. Therefore, by \cite{Ber90} 1.3.2 Corollary (ii), we obtain
\be
\n{f'}_{\cA_K(D_{X \ens{f},k})} = \sup_{\xi \in \cM_K(\cA_K(D_{X \ens{f},k}))} \v{f'(\xi)} \leq 1.
\ee
This implies
\be
f' \in \cA_K(D_{X \ens{f},k})_{\leq 1} \subset \Int(\cA_K(D_{X \ens{f},k})_{\leq 1}) = \cA_K(D_{X \ens{f},k})_{\ad}^{+},
\ee
and hence $\xi(f') \leq 1$ for any $\xi \in \Spa(\cA_K(D_{X \ens{f},k})_{\ad})$. Therefore, by the universality of the Weierstrass localisation $B$, there exists a unique morphism $\Phi \colon B \to \cA_K(D_{X \ens{f},k})_{\ad}$ of affinoid rings over $K_{\ad}$ such that the short $K$-algebra homomorphism $\cA_K(D_{X,k}) \to \cA_K(D_{X \ens{f},k})$ coincides with $\Phi \circ \Psi$.

\vs
By Lemma \ref{stably uniform 2}, $\cA_K(D_{X,k})_{\ad}$ is stably uniform, and hence $B$ is uniform. We equip $B^{\triangleleft}$ with the unique power-multiplicative norm, i.e.\ the norm given by the spectral radius. By $n_X = n_{X \ens{f}}$, the restriction map $D_{X,k} \to D_{X \ens{f},k}$ is surjective. Therefore, in order to show that $\Phi$ is an isomorphism, it suffices to show $\n{\Psi(\chi_X)}_B \leq \n{\chi_{X \ens{f}}}_{D_{X \ens{f},k}}$ for any $\chi \in (\cC/k)^{n_X \flat \times}$ by the homomorphism theorem and the universality of the completed group algebra $\cA_K(D_{X \ens{f},k})$ (cf.\ Proposition \ref{universality of completed group algebra}). We have $\chi = \prod_{i \in n_X} \theta_{\chi,i} \circ \pr_{\Lambda,i/n_X}$ and $\chi_X = \prod_{i \in n_X} (\theta_{\chi,i} \circ \pr_{\Lambda,i/n_X})_X$. By the power-multiplicativity of the norm of $B$ and \cite{Ber90} 1.3.2 Corollary (ii), we have
\be
\n{\Psi(\chi_X)}_B = \sup_{\xi \in \cM_K(B)} \v{\Phi(\chi_X)(\xi)} = \sup \set{\v{\chi_X(\xi)}}{\xi \in \cM_K(\cA_K(D_{X,k})) \land \v{f(\xi)} \leq 1}.
\ee
By Lemma \ref{evaluation and realisation 2} (1) applied to $(\theta_{\chi,i},\pr_{\Lambda,i/n_X})$ for each $i \in n_X$, the right hand side coincides with
\be
& & \sup \set{\v{\left( \prod_{i \in n_X} (\theta_{\chi,i} \circ \pr_{\Lambda,i/n_X})_{X \ens{f}} \right)(\xi)}}{\xi \in \cM_K(\cA_K(D_{X,k})) \land \v{f(\xi)} \leq 1} \\
& = & \sup \set{\prod_{i \in n_X} \v{\theta_{\chi,i}(1)}^{\Re_{X \ens{f}}(\xi)(i)}}{\xi \in \cM_K(\cA_K(D_{X,k})) \land \v{f(\xi)} \leq 1}.
\ee
By Lemma \ref{convex hull 3}, Lemma \ref{evaluation and realisation 2} (2) applied to the same data, and the continuity of the map
\be
\R^{n_X} \to (0,\infty), \ x \mapsto \prod_{i \in n_X} \v{\theta_{\chi,i}(1)}^{x(i)},
\ee
the right hand side coincides with
\be
\sup_{x \in \im(e_{X \ens{f}})} \prod_{i \in n_X} \v{\theta_{\chi,i}(1)}^{x(i)} = \sup_{x \in X \ens{f}} \prod_{i \in n_X} \v{\theta_{\chi,i}(1)}^{e_{X \ens{f}}(x)(i)} = \sup_{x \in X \ens{f}} \ev_{X \ens{f},k,x}(\chi_{X \ens{f}}) = \n{\chi_{X \ens{f}}}_{D_{X \ens{f},k}}.
\ee
\end{proof}

%% file: Cycle.tex
\section{Analytic Cycle}
\label{Analytic Cycle}

In this section, we introduce an analytic counterpart of cycles by formulating singular homology of an adic space. We continue to only consider the case $\Lambda = \Z[p^{-1}]$. Throughout this section, we suppose that the residue field $k_{\leq 1}/k_{< 1}$ of $k$ is a finite field. We denote by $q_k$ its cardinality, and put $\ell_k \coloneqq q_k - 1$. Let $n \in \N$.

\subsection{Analytic Standard Simplex}
\label{Analytic Standard Simplex}

We first introduce an analytic counterpart of the standard simplex.

\begin{dfn}
\label{analytic standard simplex}
We put $\Lambda^n_k \coloneqq \ell_k \Delta_{\Lambda}^n$, $D^n_k \coloneqq D_{\Lambda^n_k,k}$, $\cA_{K/k}^n = A_{\Lambda^n_k,K/k}$, and $\Delta_{K/k}^n = (\Lambda^n_k)_{K/k}$ (cf.\ Definition \ref{completed group algebra}, Example \ref{standard simplex}, and Definition \ref{analytic convex hull}). We denote by $\pr^n_{K/k}$ the canonical projection $\cA_K(D^n_k) \twoheadrightarrow \cA_{K/k}^n$ (cf.\ Definition \ref{completed group algebra}).
\end{dfn}

The scale $\ell_k$ will play an important role in the definition of a character in Definition \ref{integration of a character}. We promise to consider $X_{K/k}$ only for a $\Lambda$-affine germ $X$ rather than an $\R$-affine germ so that the notation $\Delta_{K/k}^n$ will not be confounded with $X_{K/k}$ with $X = \Delta^n$, which is an $\R$-affine germ rather than a $\Lambda$-affine germ. Since $\Delta_{\Lambda}^n$ is linearly convex, $\Delta_{K/k}^n$ plays a role analogous to $\Delta^n$.

\vs
We note that $\cA_K(D^n_k)$ is simpler and hence is more useful in a concrete computation than $\cA_{K/k}^n$. On the other hand, $\cA_{K/k}^n$ is more suitable for the purpose to formulate analytic singular homology because a variant defined by using $\cA_K(D^n_k)$ does not have sufficiently many cycles. Therefore, we prepare a technical proposition to help us to reduce explicit computations on $\cA_{K/k}^n$ to those on $\cA_K(D^n_k)$.

\begin{prp}
\label{kernel of comparison}
Suppose $\set{x \in \cC}{\exists r \in \N, x^{p^r} \in k} \subset K$ so that $d(x) \in K$ holds for any $(d,x) \in D_{\Delta^n_k	,k} \times \Delta^n_k$. For any $i \in [n]$, $\ker(\pr^n_{K/k})$ is the closure of the $K$-vector subspace generated by $\set{[d] - d(\ell_k \delta_{i/[n]})[d(\ell_k \delta_{i/[n]})^{-1}d]}{d \in D_{\Delta^n_k,k}}$, where the characteristic function $\delta_{i/[n]}$ is regarded as the $\Lambda$-valued one.
\end{prp}

We will use Proposition \ref{kernel of comparison} in the formulation of the integral in \S \ref{Integration along Rigid Analytic Cycles}. In order to show Proposition \ref{analytic standard simplex}, we prepare lemmata on the invertibility of a certain multivariable extension of Vandermond matrix. We note that the invertibility of the general multivariable extension of Vandermond is a difficult problem (cf.\ \cite{Olv06}), and we only deal with a quite specific case.

\begin{lmm}
\label{Vandermond}
Let $F$ be a field, $\ell$ and $m$ positive integers, $\zeta \in F^{\times}$ a primitive $\ell$-th root of $1$, and $\iota$ a bijective map $\ell^m \to (\Z/\ell \Z)^{\oplus m}$. Then the square matrix $(\omega^{\iota(i) \cdot \iota(j)})_{i,j \in \ell^m}$ of size $\ell^m$ with coefficients in $F$ is invertible, where $\cdot$ denotes the standard inner product $((\Z/\ell \Z)^{\oplus m})^2 \to \Z/\ell \Z$ and the power $\Z/\ell \Z \to F^{\times}, \ a \mapsto \omega^a$ is the one induced by the usual power $\Z \to F^{\times}, \ a \mapsto \omega^a$.
\end{lmm}

\begin{proof}
We denote by $f \colon F^{\ell^m} \to F^{\ell^m}$ the $F$-linear homomorphism given by the square matrix in the assertion. We denote by $j \colon F[(\Z/\ell \Z)^{\oplus m}] \to F^{\ell^m}$ the $F$-linear algebra isomorphism induced by the group homomorphism
\be
(\Z/\ell \Z)^{\oplus m} \to (F^{\ell^m})^{\times}, \ v \mapsto (\omega^{v(i)})_{i \in \ell^m}.
\ee
Then the composite $f \circ j$ coincides with the Fourier transform over $F$ on the discrete Abelian group $(\Z/\ell \Z)^{\oplus m}$ with respect to $\omega$. Since the charasteristic of $F$ is coprime to the order $\ell^m$ of $(\Z/\ell \Z)^{\oplus m}$ due to the existence of a primitive $\ell$-th root of unity, the Fourier transform $f \circ j$ is an $F$-linear isomorphism. This implies $f = (f \circ j) \circ j^{-1}$ is an $F$-linear isomorphism, and hence its matrix presentation is invertible.
\end{proof}

\begin{lmm}
\label{coefficient reconstruction}
Let $f \in K[D_{\Delta^n_k,k}]$ and $i \in [n]$.
\bi
\item[(1)] The inequality $\n{f}_{\cA_K(D_{\Delta^n_k,k})} \geq \n{\pr^n_{K/k}(f)}_{\cA_{K/k}^n}$ holds.
\item[(2)] If $f \in \ker(\pr^n_{K/k})$, then $f$ belongs to the $K$-vector subspace generated by the subset $\set{[d] - d(\ell_k \delta_{i/[n]})[d(\ell_k \delta_{i/[n]})^{-1}d]}{d \in D_{\Delta^n_k,k}}$.
\item[(3)] If $f(d) = 0$ for any $d \in D_{\Delta^n_k,k}$ with $d(\ell_k \delta_{i/[n]}) \neq 1$, then the equality $\n{f}_{\cA_K(D_{\Delta^n_k,k})} = \n{\pr^n_{K/k}(f)}_{\cA_{K/k}^n}$ holds.
\ei
\end{lmm}

\begin{proof}
The assertion (1) follows from the strong triangular inequality. The assertion (2) follows from (3), since $f$ replaced by $f - \sum_{d \in D_{\Delta^n_k,k}} f(d)([d] - d(\ell_k \delta_{i/[n]})[d(\ell_k \delta_{i/[n]})^{-1}d])$ satisfies the condition in (3). Therefore, it suffices to show $\n{f}_{\cA_K(D_{\Delta^n_k,k})} < (1 + \epsilon) \n{\pr^n_{K/k}(f)}_{\cA_{K/k}^n}$ for any $\epsilon \in (0,\infty)$ under the assumption $f(d) = 0$ for any $d \in D_{\Delta^n_k,k}$ with $d(\ell_k \delta_{i/[n]}) \neq 1$. By the functoriality in Proposition \ref{functoriality of the analytification} (1) and (3), it suffices to consider only the case $i = 0$.

\vs
Put $D \coloneqq \set{d \in D_{\Delta^n_k,k}}{d(\ell_k \delta_{i/[n]}) = 1}$. We denote by $D[\ell_k] \subset D$ the subgroup of torsion elements. By the definition, $D[\ell_k]$ coincides with the finite group of exponential maps $d \colon \Delta^n_k \to k$ with $d(\ell_k \delta_{0/[n]}) = 1$ whose image is contained in the subgroup of $k^{\times}$ consisting of $\ell_k$-th root of $1$.

\vs
Put $S \coloneqq \set{d \in D}{\exists e \in D[\ell_k], f(de) \neq 0}$. Since $D[\ell_k]$ is a finite group, $S$ is a finite set. We put $[n]_{>0} \coloneqq [n] \setminus \ens{0}$ and $X \coloneqq \set{x \in (\Lambda \cap [0,\infty))^{[n]_{>0}}}{\sum_{i=1}^{n} x(i) \leq 1}$. We denote by $s$ the embedding
\be
\Q^{[n]_{>0}} \hookrightarrow \set{y \in \Q^{[n]}}{\sum_{i \in [n]} y(i) = 1}, \ x \mapsto 
\left(
\begin{array}{rcl}
[n] \to \Q, \ i \mapsto 
\left\{
\begin{array}{ll}
1 - \sum_{j=1}^{n} x(j) & (i = 0) \\
x(i) & (i \neq 0)
\end{array}
\right.
\end{array}
\right).
\ee
We take a complete system $S_0$ of representatives of $S$ with respect to the action of $D[\ell_k]$ given by the multiplication. For each $d \in S_0 \cup D[\ell_k]$, we fix a $\chi_d \in (\cC/k)^{n \flat \times}$ with $d = (\chi_d)_{\Delta^n_k}$. Let $d \in S_0$. We denote by $d \otimes \Q$ the $\Q$-linear homomorphism $\Q^{[n]_{>0}} \cong \Lambda^{[n]_{>0}} \otimes_{\Z} \Q \to \cC^{\times} \otimes_{\Z} \Q$ associated to the group homomorphism
\be
\Lambda^{[n]_{>0}} \to \cC^{\times}, \ x \mapsto \chi_d(s(x)).
\ee
By the definition, we have $d(x) \otimes 1 = (d \otimes \Q)((x(i))_{i \in [n]_{>0}})$. In this sense, $d \otimes \Q$ possesses all the information of $d$ modulo $\ell_k$-torsion.

\vs
We denote by $s_0 \colon \Q^{[n]_{>0}} \hookrightarrow \Q^{[n]}$ the zero extension, which is $\Q$-linear unlike $s$. For any $(d,x) \in (S_0 \cup D[\ell_k]) \times \Q^{[n]_{>0}}$, we have $\chi_d(s(x)) = \chi_d(s_0(x))$ by $\chi_d(\ell_k \delta_{0/[n]}) = 1$. Take a sufficiently large $N_0 \in \N$ with $(q_k-2) n < 2^{-1} q_k^{N_0}$. Fix a primitive $\ell_k$-th root $\zeta \in k^{\times}$ of $1$.  For each $(e,i) \in D[\ell_k] \times [n]_{>0}$, we denote by $\ell_{e,i} \in \ell_k$ the discrete logarithm of $e(\ell_k \delta_{i/[n]})$ with base $\zeta$. By Lemma \ref{Vandermond} applied to $\zeta$, the extended Vandermond matrix $(\zeta^{\sum_{i=1}^{n} \ell_{e',i} \ell_{e,i}})_{(e,e') \in D[\ell_k]^2}$ admits an inverse $B_{\zeta}$.

\vs
Let $e \in D[\ell_k]$, and put $y_e \coloneqq (q_k^{-N_0} \ell_{e,i})_{i \in [n]_{>0}} \in (\Q \cap [0,(2n)^{-1}))^{[n]_{>0}}$. For any $(d,m,v) \in S_0 \times \N \times X$, we have
\be
& & \v{\chi_d(s(p^{-m} v))} = \v{\chi_d(s_0(p^{-m} v))} = \v{\chi_d(p^{-m} s_0(v))} \\
& = & \left( \v{\chi_d(p^{-m} s_0(v))}^{p^m} \right)^{p^{-m}} = \v{\chi_d(s_0(v))}^{p^{-m}} = \v{\chi_d(s(v))}^{p^{-m}}.
\ee
Therefore, replacing $y_e$ by $p^{-m} y_e$ for a sufficiently large $m \in \N$ by replacing $N_0$ by sufficiently larger one and putting $C \coloneqq (1+\epsilon)^{2(3 \# S_0 - 1)\# S_0)^{-1}} \in (1,\infty)$, we may assume $\v{\chi_d(s(y_e))} > C^{-1}$.

\vs
For any $(d_0,d_1) \in S_0^2$ with $d_0 \neq d_1$, since $d_0 d_1^{-1} \notin D[\ell_k]$, we have $d_0 \otimes \Q \neq d_1 \otimes \Q$. Therefore, the set
\be
\set{x \in \Q^{[n]_{>0}}}{\exists (d_0,d_1) \in S_0^2, d_0 \neq d_1 \land (d_0 \otimes \Q)(x) = (d_1 \otimes \Q)(x)}
\ee
is a finite union of proper $\Q$-linear subspaces, which is a nowhere dense closed subset of $\Q^{[n]_{>0}}$ with respect to the Euclidean topology. In particular, there exists an $x_e \in (\Q \cap [0,\infty))^{[n]_{>0}}$ such that $(d_0 \otimes \Q)(m x_e + y_e) \neq (d_1 \otimes \Q)(m x_e + y_e)$ for any $(d_0,d_1,m) \in S_0^2 \times \Q^{\times}$ with $d_0 \neq d_1$. Replacing $x_e$ by $m! x_e$ for a sufficiently large $m \in \N$, we may assume $x_e \in (\Lambda \cap [0,\infty))^{[n]_{>0}}$. Replacing $x_e$ by $p^{-m} x_e$ for a sufficiently large $m \in \N$, we may assume $\sum_{i=1}^{n} x_e(i) < 2^{-1}\ell_k^{-1}$. In particular, we have
\be
\sum_{i=1}^{n} (x_e + y_e)(i) = \sum_{i=1}^{n} x_e(i) + q_k^{-N_0} \sum_{i=1}^{n} \ell_{e,i} < 2^{-1}\ell_k^{-1} + q_k^{-N_0}(q_k-2) n < 2^{-1}(\ell_k^{-1} + 1) \leq 1.
\ee
This implies $x_e + y_e \in X$, and hence $d(s(x_e + y_e))$ makes sense for any $d \in D_{\Delta^n_k,k}$. For any $(d_0,d_1) \in S_0^2$ with $d_0 \neq d_1$, $(d_0 \otimes \Q)(x_e + y_e) \neq (d_1 \otimes \Q)(x_e + y_e)$ implies $d_0(s(x_e + y_e)) d_1(s(x_e + y_e))^{-1}$ is not an $\ell_k$-torsion, and hence $d_0(s(x_e + y_e)) \neq d_1(s(x_e + y_e))$.

\vs
For any $m \in \N \setminus \ens{0}$, we have
\be
& & \v{d_0(s(p^{-m} x_e + y_e)) - d_1(s(p^{-m} x_e + y_e))} \\
& = & \v{\chi_{d_0}(s(p^{-m} x_e + y_e)) - \chi_{d_1}(s(p^{-m} x_e))} = \v{\chi_{d_0}(s_0(p^{-m} x_e)) - \chi_{d_1}(s_0(p^{-m} x_e))} \\
& = & \v{\chi_{d_0}(p^{-m} s_0(x_e)) - \chi_{d_1}(p^{-m} s_0(x_e))} = \v{\chi_{d_1}(p^{-m} s_0(x_e))} \v{\chi_{d_0} \chi_{d_1}^{-1}(p^{-m} s_0(x_e)) - 1} \\
& = & \v{\chi_{d_1}(s_0(x_e))}^{p^{-m}} \v{\frac{\chi_{d_0} \chi_{d_1}^{-1}(p^{-m} s_0(x_e))^{p^m} - 1}{\sum_{i \in p^m} \chi_{d_0} \chi_{d_1}^{-1}(p^{-m} s_0(x_e))^i}} \\
& = & \v{\chi_{d_1}(s_0(x_e))}^{p^{-m}} \v{\frac{\chi_{d_0} \chi_{d_1}^{-1}(s_0(x_e)) - 1}{\sum_{i \in p^m} \left( 1 + (\chi_{d_0} \chi_{d_1}^{-1}(p^{-m} s_0(x_e)) - 1) \right)^i}} \\
& = & \v{\chi_{d_1}(s_0(x_e))}^{p^{-m}} \v{\frac{\chi_{d_0} \chi_{d_1}^{-1}(s_0(x_e)) - 1}{p^m + \sum_{i \in p^m} \sum_{j=1}^{i} \binom{i}{j}(\chi_{d_0} \chi_{d_1}^{-1}(p^{-m} s_0(x_e)) - 1)^j}} \\
& = & \frac{\v{\chi_{d_1}(s_0(x_e))}^{p^{-m}}}{\Bigv{\frac{p^m}{\chi_{d_0} \chi_{d_1}^{-1}(s_0(x_e)) - 1} + 1 + (\chi_{d_0} \chi_{d_1}^{-1}(p^{-m} s_0(x_e)) - 1) \sum_{j=1}^{p^m-2} \left( \sum_{i=j+1}^{p^m-1} \binom{i}{j} \right) (\chi_{d_0} \chi_{d_1}^{-1}(p^{-m} s_0(x_e)) - 1)^{j-1}}}.
\ee
This implies that if $\v{\chi_{d_0} \chi_{d_1}^{-1}(s_0(x_e)) - 1} < 1$, we have
\be
& & \lim_{m \to \infty} \v{d_0(s(p^{-m} x_e + y_e)) - d_1(s(p^{-m} x_e))} \\
& = & \lim_{m \to \infty} \frac{\v{\chi_{d_1}(s_0(x_e))}^{p^{-m}}}{\Bigv{\frac{p^m}{\chi_{d_0} \chi_{d_1}^{-1}(s_0(x_e)) - 1} + 1 + (\chi_{d_0} \chi_{d_1}^{-1}(p^{-m} s_0(x_e)) - 1) \sum_{j=1}^{p^m-2} \left( \sum_{i=j+1}^{p^m-1} \binom{i}{j} \right) (\chi_{d_0} \chi_{d_1}^{-1}(p^{-m} s_0(x_e)) - 1)^{j-1}}} \\
& = & \lim_{m \to \infty} \v{\chi_{d_1}(s_0(x_e))}^{p^{-m}} = 1.
\ee
On the other hand, if $\v{\chi_{d_0} \chi_{d_1}^{-1}(s_0(x_e)) - 1} \geq 1$, then we already have
\be
\v{\chi_{d_0} \chi_{d_1}^{-1}(s_0(p^{-m} x_e)) - 1} \geq 1 > 2^{-((\# S_0 - 1)\# S_0)^{-1}}
\ee
for any $m \in \N$. Therefore, replacing $x_e$ by $p^{-m} x_e$ for a sufficiently large $m \in \N$, we may assume $C^{-1} < \v{\chi_{d_1}(s_0(p^{-m} x_e))} < C$ and $C^{-1} < \v{\chi_{d_0} \chi_{d_1}^{-1}(s_0(p^{-m} x_e)) - 1}$ for any $m \in \N$.

\vs
We denote by $N_1 \in \N$ the minimum of an $m \in \N$ satisfying $\# S_0 \leq p^{m}$. Let $c \in \# S_0$. Put $v_{e,c} \coloneqq c p^{-N_1} \ell_k x_e + y_e \in (\Lambda \cap [0,\infty))^{[n]_{>0}}$. We have
\be
\sum_{i=1}^{n} v_{e,c}(i) & = & c p^{-N_1} \ell_k \sum_{i=1}^{n} x_e(i) + q_k^{-N_0} \sum_{i=1}^{n} \ell_{e,i} < c p^{-N_1} \ell_k \times 2^{-1}\ell_k^{-1} + q_k^{-N_0} (q_k-2) n \\
& < & \# S_0 p^{-N_1} 2^{-1} + 2^{-1} \leq 2^{-1} + 2^{-1} = 1.
\ee
This implies $v_{e,c} \in X$, and hence $d(s(v_{e,c}))$ makes sense for any $d \in D_{\Delta^n_k,k}$. we have
\be
& & f(s(v_{e,c})) = \sum_{d \in D_{\Delta^n_k,k}} f(d) d(s(v_{e,c})) = \sum_{d \in S} f(d) d(s(v_{e,c})) \\
& = & \sum_{d \in S_0} \sum_{e' \in D[\ell_k]} f(de') (de')(s(v_{e,c})) = \sum_{d \in S_0} \sum_{e' \in D[\ell_k]} f(de') d(s(v_{e,c})) e'(s(v_{e,c})) \\
& = & \sum_{d \in S_0} \sum_{e' \in D[\ell_k]} f(de') \chi_d(s(v_{e,c})) \chi_{e'}(s(v_{e,c})) \\
& = & \sum_{d \in S_0} \sum_{e' \in D[\ell_k]} f(de') \chi_d(s_0(c p^{-N_1} \ell_k x_e + y_e)) \chi_{e'}(s_0(c p^{-N_1} \ell_k x_e + y_e)) \\
& = & \sum_{d \in S_0} \sum_{e' \in D[\ell_k]} f(de') \chi_d(s_0(p^{-N_1} x_e))^c \chi_d(s_0(y_e)) \chi_{e'}(s_0(p^{-N_1} x_e))^c \zeta^{\sum_{i=1}^{n} \ell_{e',i} \times q_k^{-N_0} \ell_{e,i}} \\
& = & \sum_{d \in S_0} \sum_{e' \in D[\ell_k]} f(de') \chi_d(s(p^{-N_1} x_e))^c \chi_d(s(y_e)) (\zeta^{(1+\ell_k)^N_0})^{q_k^{-N_0} \sum_{i=1}^{n} \ell_{e',i} \ell_{e,i}} \\
& = & \sum_{d \in S_0} \sum_{e' \in D[\ell_k]} f(de') \chi_d(s(p^{-N_1} x_e))^c \chi_d(s(y_e)) \zeta^{\sum_{i=1}^{n} \ell_{e',i} \ell_{e,i}} \\
& = & \sum_{d \in S_0} \left( \sum_{e' \in D[\ell_k]} f(de') \zeta^{\sum_{i=1}^{n} \ell_{e',i} \ell_{e,i}} \chi_d(s(y_e)) \right) \chi_d(s(p^{-N_1} x_e))^c.
\ee
By $\chi_{d_0}(s(p^{-N_1} x_e)) \neq \chi_{d_1}(s(p^{-N_1} x_e))$ for any $(d_0,d_1) \in S_0^2$ with $d_0 \neq d_1$, Vandermond's matrix for $(\chi_d(s(p^{-N_1} x_e)))_{d \in S_0}$ admits an inverse $A_e$. We denote by $D_e$ the diagonal matrix whose index set is $S_0$ and whose diagonal entries are given by $(\chi_d(s(y_e))^{-1})_{d \in S_0}$. By the computation above, the coefficient vector $(f(de'))_{(d,e') \in S_0 \times }D[\ell_k]$ is expressed as the matrix product of $B_{\zeta}^{\oplus S_0}$, $\bigoplus_{e \in D[\ell_k]} D_e$, and $\bigoplus_{e \in D[\ell_k]} A_e$ applied to $(f(s(v_{e,c}))_{(e,c) \in D[\ell_k] \times \#S_0}$. In particular, we have
\be
& & \n{f}_{\cA_K(D_{\Delta^n_k,k})} = \sup_{d \in S} \v{f(d)} = \sup_{(d,e') \in S_0 \times D[\ell_k]} \v{f(de)} \\
& \leq & \n{B_{\zeta}} \max_{e \in D[\ell_k]} \n{D_e} \max_{e \in D[\ell_k]} \n{A_e} \sup_{(e,c) \in D[\ell_k] \times \#S_0} \v{f(s(v_{e,c}))} \\
& \leq & \n{B_{\zeta}} \max_{e \in D[\ell_k]} \sup_{d \in S_0} \v{\chi_d(s(y_e))^{-1}} \max_{e \in D[\ell_k]} \left( \prod_{(d_0,d_1) \in S_0, d_0 \neq d_1} \v{\chi_{d_0}(s(p^{-N_1} x_e)) - \chi_{d_1}(s(p^{-N_1} x_e))} \right)^{-1} \\
& & \left( \max_{d \in S_0} \v{\chi_d(s(p^{-N_1} x_e))} \right)^{2^{-1}(\# S_0 - 1) \# S_0} \n{\pr^n_{K/k}(f)}_{\cA_{K/k}^n},
\ee
where the norms of matrices denote the operator norms with respect to the $\ell^{\infty}$-norms. The operator norm $\n{B_{\zeta}}$ coincides with $1$, because the reduction of the extended Vandermond matrix is again an extended Vandermond matrix over the residue field $k_{\leq 1}/k_{< 1}$ of $k$, which is still invertible by Lemma \ref{Vandermond} applied to the residue class $\zeta + k_{< 1}$ of $\zeta$. Estimating the operator norm $\n{A_e}$ by the valuations of $\det(A_e)$ and entries of $A_e$, we obtain
\be
& & \n{B_{\zeta}} \max_{e \in D[\ell_k]} \sup_{d \in S_0} \v{\chi_d(s(y_e))^{-1}} \max_{e \in D[\ell_k]} \left( \prod_{(d_0,d_1) \in S_0, d_0 \neq d_1} \v{\chi_{d_0}(s(p^{-N_1} x_e)) - \chi_{d_1}(s(p^{-N_1} x_e))} \right)^{-1} \\
& & \left( \max_{d \in S_0} \v{\chi_d(s(p^{-N_1} x_e))} \right)^{2^{-1}(\# S_0 - 1) \# S_0} \n{\pr^n_{K/k}(f)}_{\cA_{K/k}^n} \\
& = & 1 \max_{e \in D[\ell_k]} \left( \inf_{d \in S_0} \v{\chi_d(s(y_e))} \right)^{-1} \max_{e \in D[\ell_k]} \left( \prod_{(d_0,d_1) \in S_0, d_0 \neq d_1} \v{\chi_{d_0}\chi_{d_1}^{-1}(s(p^{-N_1} x_e)) - 1} \v{\chi_{d_1}(s(p^{-N_1} x_e))} \right)^{-1} \\
& & \left( \max_{d \in S_0} \v{\chi_d(s(p^{-N_1} x_e))} \right)^{2^{-1}(\# S_0 - 1) \# S_0} \n{\pr^n_{K/k}(f)}_{\cA_{K/k}^n} \\
& < & \max_{e \in D[\ell_k]} \left( \inf_{d \in S_0} C^{-1} \right)^{-1} \max_{e \in D[\ell_k]} \left( \prod_{(d_0,d_1) \in S_0, d_0 \neq d_1} C \times C \right)^{-1} \left( \max_{d \in S_0} C \right)^{2^{-1}(\# S_0 - 1) \# S_0} \n{\pr^n_{K/k}(f)}_{\cA_{K/k}^n} \\
& = & C^{2^{-1} 3(\# S_0 - 1) \# S_0} \n{\pr^n_{K/k}(f)}_{\cA_{K/k}^n} = (1 + \epsilon) \n{\pr^n_{K/k}(f)}_{\cA_{K/k}^n}.
\ee
This implies $\n{f}_{\cA_K(D_{\Delta^n_k,k})} < (1 + \epsilon) \n{\pr^n_{K/k}(f)}_{\cA_{K/k}^n}$
\end{proof}

\begin{proof}[Proof of Proposition \ref{kernel of comparison}]
We denote by by $I$ the $K$-vector subspace of $\cA_K(D_{\Delta^n_k,k})$ generated by $\set{[d] - d(\ell_k \delta_{i/[n]})[d(\ell_k \delta_{i/[n]})^{-1}d]}{d \in D_{\Delta^n_k,k}}$. We have $I \subset \ker(\pr^n_{K/k})$ by the definition of $\pr^n_{K/k}$. Let $(f,\epsilon) \in \ker(\pr^n_{K/k}) \times (0,\infty)$. By the definition of $\cA_K(D_{\Delta^n_k,k})$, there exists a $g \in K[D_{\Delta^n_k,k}]$ such that $\n{f-g}_{\cA_K(D_{\Delta^n_k,k})} < \epsilon$. By $f \in \ker(\pr^n_{K/k})$, we have $\pr^n_{K/k}(g) = - \pr^n_{K/k}(f -g)$. Since $\pr^n_{K/k}$ is short, we obtain
\be
\n{\pr^n_{K/k}(g)}_{\cA_{K/k}^n} = \n{- \pr^n_{K/k}(f-g)}_{\cA_{K/k}^n} \leq \n{f-g}_{\cA_K(D_{\Delta^n_k,k})} < \epsilon.
\ee
We denote by $h \in K[D_{\Delta^n_k,k}] \cap I$ the essentially finite sum
\be
\sum_{d \in D_{\Delta^n_k,k}} g(d)([d] - d(\ell_k \delta_{i/[n]})[d(\ell_k \delta_{i/[n]})^{-1}d]).
\ee
For any $d \in D_{\Delta^n_k,k}$, we have
\be
(d(\ell_k \delta_{i/[n]})^{-1}d)[\ell_k \delta_{i/[n]}] = (d(\ell_k \delta_{i/[n]})^{-1}d)(\ell_k \delta_{i/[n]}) = 1.
\ee
Therefore, for any $d \in D_{\Delta^n_k,k}$ with $d(\ell_k \delta_{i/[n]}) \neq 1$, we have $(g-h)(d) = g(d) - g(d) = 0$. Moreover, by Lemma \ref{coefficient reconstruction} (2) applied to $g - h$, we obtain $\n{g - h}_{\cA_K(D_{\Delta^n_k,k})} = \n{\pr^n_{K/k}(g - h)}_{\cA_{K/k}^n}$. This implies
\be
& & \n{f - h}_{\cA_K(D_{\Delta^n_k,k})} \leq \n{f - g}_{\cA_K(D_{\Delta^n_k,k})} + \n{g - h}_{\cA_K(D_{\Delta^n_k,k})} = \n{f - g}_{\cA_K(D_{\Delta^n_k,k})} + \n{\pr^n_{K/k}(g - h)}_{\cA_{K/k}^n} \\
& = & \n{f - g}_{\cA_K(D_{\Delta^n_k,k})} + \n{\pr^n_{K/k}(g)}_{\cA_{K/k}^n} = \n{f - g}_{\cA_K(D_{\Delta^n_k,k})} + \n{- \pr^n_{K/k}(f - g)}_{\cA_{K/k}^n}.
\ee
By Lemma \ref{coefficient reconstruction} (1) applied to $f - g$,
\be
\n{f - g}_{\cA_K(D_{\Delta^n_k,k})} + \n{- \pr^n_{K/k}(f - g)}_{\cA_{K/k}^n} \leq 2 \n{f - g}_{\cA_K(D_{\Delta^n_k,k})} < 2 \epsilon.
\ee
This implies $\n{f - h}_{\cA_K(D_{\Delta^n_k,k})} < 2 \epsilon$. This completes the proof by $h \in I \subset \ker(\pr^n_{K/k})$.
\end{proof}

\subsection{Analytic Singular Homology}
\label{Analytic Singular Homology}

We have fixed an algebraic closure $k^{\alg}$ of $k$, and denote by $\cC$ the completion of $k^{\alg}$. Let $K$ be a perfectoid subfield of $\cC$ containing $\set{x \in \cC}{\exists r \in \N, x^{p^r} \in k}$ as in Proposition \ref{kernel of comparison} which is Galois in the sense that it corresponds to a normal closed (not necessarily open) subgroup $G_K$ of the absolute Galois group $G_k \coloneqq \Gal(k^{\alg}/k)$ of $k$. Put $G_{K/k} \coloneqq G_k/G_K$. The affinoid pre-adic space $\Delta_{K/k}^n$ over $\Spa(K_{\ad})$ defined in Definition \ref{analytic standard simplex} is an affinoid perfectoid space by Theorem \ref{sheafy germ}, and is regarded as an affinoid adic space over $\Spa(k_{\ad})$ through the canonical morphism $\Spa(K_{\ad}) \to \Spa(k_{\ad})$.

\vs
We denote by $\Delta$ the category of non-zero finite ordinals and order-preserving maps, by $\sSet$ the category of simplicial sets and simplicial maps, by $\Ab_p$ the category of $p$-divisible normed groups (cf.\ \S \ref{p-divisible normed Group}) and short group homomorphisms, by $\Ch_{\Ab}$ the category of chain complices and chain homomorphisms, by $\Alg(k)$ the category of Banach $k$-algebras and short $k$-algebra homomorphisms, by $\PreAdic/k \subset (V)$ the subcategory of pre-adic spaces over $\Spa(k_{\ad})$ and morphisms over $\Spa(k_{\ad})$, and by $\Adic/k \subset \PreAdic/k$ the full subcategory of adic spaces over $\Spa(k_{\ad})$.

\vs
The correspondence $[n] \mapsto \Lambda^{\oplus [n]}$ gives a functor $\Delta \to \Ab$, and hence the correspondence $[n] \mapsto \Hom_{\Ab}(\Lambda^{[n]},\cC^{\times})$ gives a functor $\Delta \to \Ab$. The correspondence $[n] \mapsto (\cC/k)^{(n+1) \flat \times}$ gives its subfunctor $\Delta \to \Ab_p^{\op}$. By the definition, the correspondence $[n] \mapsto \cA_k^n$ gives a functor $\Delta \to \Alg(k)^{\op}$. Since the correspondence $A \mapsto \Spa(A_{\ad})$ gives a functor $\Alg(k)^{\op} \to \PreAdic/k$, the correspondence $[n] \mapsto \Delta_{K/k}^n$ gives a functor $\Delta_{K/k}^{\bullet} \colon \Delta \to \PreAdic/k$. Thus we obtained a cosimplicial object $\Delta_{K/k}^{\bullet}$ in $\PreAdic/k$. Since $\Delta_{K/k}^n$ is an affinoid adic space for any $n \in \N$, $\Delta_{K/k}$ forms a cosimplicial object in $\Adic/k$. We denote by $\cN_{K/k}$ the nerve functor $\PreAdic/k \to \sSet$ associated to $\Delta_{K/k}^{\bullet}$.

\vs
Let $M$ be an Abelian group. We denote by $C_{K/k,M} \colon \PreAdic/k \to \Ch_{\Ab}$ the singular complex with coefficients in $M$ associated to $\cN_{K/k}$. For a subfunctor $C \subset C_{K/k,M}$, we denote by $\rH_*(\bullet,C)$ (resp.\ $\rH^*(\bullet,C)$) the composite of $C$ and the homology $\Ch_{\Ab} \to \Ab^{\N}$ (resp.\ the cohomology $\Ch_{\Ab} \to (\Ab^{\op})^{\N}$).

\vs
For example, $\rH_*(\bullet,C_{K/k,M})$ (resp.\ $\rH^*(\bullet,C_{K/k,M})$)  is an analytic analogue of the singular homology (resp.\ the singular cohomology) with coefficients in $M$. The continuous isometric left action of $G_{K/k}$ on $K^{\times}$ induces a left action of $G_{K/k}$ on $D^n_k$ (cf.\ Definition \ref{analytic standard simplex}) preserving the weight by the assumption $\set{x \in \cC}{\exists r \in \N, x^{p^r} \in k} \subset K$. It induces an isometric $k$-algebra left action of $G_{K/k}$ on $K[D^n_k]$, which extends to an isometric $k$-algebra left action of $G_{K/k}$ on $\cA_K(D^n_k)$, which extends to an isometric $k$-algebra left action of $G_{K/k}$ on $\cA_K^n$ by Proposition \ref{kernel of comparison}. It induces a right action of $G_{K/k}$ on $\Delta_{K/k}^n$ over $\Spa(k_{\ad})$, which induces a right action of $G_{K/k}$ on $\Delta_{K/k}^{\bullet}$, which induces a left action of $G_{K/k}$ on $\cN_{K/k}$, which induces a left action of $G_{K/k}$ on $C_{K/k,M}$. If a subfunctor $C \subset C_{K/k,M}$ is stable under the action of $G_{K/k}$, the restriction of the action of $G_{K/k}$ to $C$ induces a left action of $G_{K/k}$ on $\rH_*(\bullet,C)$ (resp.\ a right action of $G_{K/k}$ on $\rH^*(\bullet,C)$). In this way, we obtain Galois representations in functorial ways.

\vs
As the singular homology can be regarded as a group of equivalence classes of continuous cycles, it is natural to regard the analytic analogue of the singular homology as a group of equivalence classes of ``analytic cycles''. One natural direction to formulate an integration of differential forms on $X$ along ``analytic cycles'' is to assume the overconvergence of differential forms by equipping $X$ with the structure as the completion of a $k$-dagger space (cf.\ \cite{Klo00}). However, it works only when $X$ is locally of finite type over $\Spa(k_{\ad})$. Instead, we consider the integrability of a given differential form along a given ``analytic cycle'' in \S \ref{Integration along Rigid Analytic Cycles}.

%% file: Integration.tex
\section{$p$-adic Period and Integration}
\label{p-adic Period and Integration}

We formulate an integration along an ``analytic cycle'' in terms of $p$-adic periods. Like the requirement of the overconvergence in the classical integrability, we require technical conditions of ``rigidity'' on an ``analytic cycle'' in the formulation. In the following, we assume that $k$ is a finite dimensional extension of $\Qp$, and denote by $\kappa$ the residue field $k_{\leq 1}/k_{< 1}$.

\subsection{Rigid Analytic Singular Simplex}
\label{Rigid Analytic Singular Simplex}

In order to formulate the integration, we introduce the notion of a rigid analytic singular simplex. Before that, we briefly recall the $p$-adic period ring $\BdR$, because we will refer to it as the codomain of the integration along a rigid analytic singular simplex.

\vs
We denote by $R$ the $\kappa$-algebra $\varprojlim_{x \mapsto x^p} \cC_{\leq 1}/p \cC_{\leq 1}$, and identify the underlying multiplicative monoids of $R$ and $\Frac(R)$ with the multiplicative monoids $\varinjlim_{x \mapsto x^p} \cC_{\leq 1}$ and $\varinjlim_{x \mapsto x^p} \cC$ respectively in the usual way. For each $(\alpha,\nu) \in \Frac(R) \times \N$, we denote by $\alpha(\nu) \in \cC$ the $\nu$-th entry of $\alpha$ regarded as an element of $\varprojlim_{x \mapsto x^p} \cC \subset \cC^{\N}$. The convention is ambiguous when $\alpha \in R$ because it is an element of $\varprojlim_{x \mapsto x^p} \cC_{\leq 1}/p \cC_{\leq 1} \subset (\cC_{\leq 1}/p \cC_{\leq 1})^{\N}$, for which the substitution makes sense, but we solve the ambiguity by clarifying that $\alpha(\nu)$ always means the evaluation in the introduced sense.

\vs
We put $W(R)_k \coloneqq k_{\leq 1} \otimes_{W(\kappa)} W(R)$. We equip $R$ with the inverse limit topology of $\varprojlim_{x \mapsto x^p} \cC_{\leq 1}/p \cC_{\leq 1}$, $W(R)$ with the direct product topology of $R$, $W(R)[p^{-1}]$ with the strongest topology $\tau_{\can}$ as a topological $W(R)$-algebra, $W(R)_k[p^{-1}]$ with the strongest topology as a topological $W(R)[p^{-1}]$-module (cf.\ \cite{War93} 2.19 Corollary), $\BdR^{+}$ the inverse limit topology of $W(R)_k[p^{-1}]$, and $\BdR$ with the strongest topology as a topological $\BdR^{+}$-algebra. See \S \ref{Appendix} for details on the topologies. We note that there are several other choices of a topology of $W(R)[p^{-1}]$ analogous to that of Laurent power series.

\begin{dfn}
\label{rigid datum}
A {\it $(k,n)$-exponential datum} is a finite subset $F \subset D^n_k$. We denote by $\ED(k,n)$ the set of $(k,n)$-exponential data.
\end{dfn}

Let $F$ be a $(k,n)$-exponential datum. We denote by $K \ens{F}$ the quotient normed $K$-algebra of the Tate algebra $K \set{\n{d}_{D^n_k}^{-1} T_d}{d \in F}$ by the kernel of the unique $K$-algebra morphism to $\cA_{K/k}^n$ assigning $d$ to $T_d$ for each $d \in F$. For an $i \in \N$, we abbreviate $\rH^0(\Spa(K \ens{F}_{\ad}),\Omega^i_{\Spa(K \ens{F}_{\ad})})$ to $\Omega^i_F$, and equip it with the canonical topology as a topological $K \ens{F}$-module. By definition, $K \ens{F}$ is an affinoid $K$-algebra. The correspondence $F \rightsquigarrow K \ens{F}$ forms a functor from the directed set of $(k,n)$-exponential data to the category of affinoid $K$-algebras and admissible injective $K$-algebra homomorphisms.

\vs
The set $\ED(k,n)$ is directed with respect to the inclusion. Let $S$ be a (possibly non-directed) subset of $\ED(k,n)$. We will explicitly consider the case where $S$ is the set $\ID(K/k,n)$, which will be defined in \S \ref{Integration along Rigid Analytic Cycles}.

\begin{dfn}
\label{analytic differenial form}
For an $i \in \N$, we denote by $\Omega^i_S$ the colimit
\be
\colim_{F \in S} \Omega^i_F,
\ee
as a topological space
\end{dfn}

We formally regard $\Omega^i_S$ as the space of ``differential $i$-forms'' on $\Delta_{K/k}^n$ in order to avoid the difficulty based on the fact that $\cA_{K/k}^n$ is not an affinoid algebra when $n > 0$.

\vs
Let $\gamma \colon \Delta_{K/k}^n \to X$ be a morphism of pre-adic spaces over $\Spa(k_{\ad})$. An {\it $S$-rigid analytic factorisation} of $\gamma$ is a pair $(F,g)$ of an $F \in I$ and a morphism $g \colon \Spa(K \ens{F}_{\ad}) \to X$ of pre-adic spaces over $\Spa(k)$ such that the diagram
\be
\xymatrix{
\Delta_{K/k}^n \ar[r] \ar[rd]_-{\gamma} & \Spa(K \ens{F}_{\ad}) \ar[d]^-{g} \\
& X
}
\ee
commutes. The set of $S$-rigid analytic factorisations of $f$ is partially ordered with respect to the partial order $(F_0,g_0) \leq (F_1,g_1)$ defined as $F_0 \subset F_1$ and the diagram
\be
\xymatrix{
\Spa(K \ens{F_1}_{\ad}) \ar[r] \ar[rd]_-{g_1} & \Spa(K \ens{F_0}_{\ad}) \ar[d]^-{g_0} \\
& X
}
\ee
commutes. We say that $\gamma$ is an {\it $S$-rigid analytic singular $n$-simplex on $X$ over $K$} if the partially ordered set is a non-empty directed set. This is what we mentioned as a technical condition at the beginning of \S \ref{p-adic Period and Integration}.

\vs
Let $\gamma \colon \Delta_{K/k}^n \to X$ be an $S$-rigid analytic singular $n$-simplex on $X$ over $K$. Then it admits an $S$-rigid analytic factorisation, and induces a system of $k$-linear homomorphisms $\rH^0(X,\Omega^n_X) \to \Omega^n_F$ indexed by $S$-rigid analytic factorisations. The system gives a $k$-linear homomorphism $\rH^0(X,\Omega^n_X) \to \Omega^n_S$.

\begin{dfn}
\label{compatible integration}
Fix a $k$-algebra embedding $K \hookrightarrow \BdR^{+}$ compatible with the reduction $\BdR^{+} \to \cC$, through which we regard $\BdR$ as a $K$-algebra. A family $I$ of $K$-linear homomorphisms $I_F \colon \Omega^n_F \to \BdR$ indexed by $F \in S$ is said to be a {\it $S$-compatible integration} if it forms a compatible system with respect to the partial order on $S$ and $I_F$ is continuous with respect to the canonical topology on $\Omega^n_F$ as a finitely generated $K \ens{F}$-module for any $F \in S$. In particular, $I$ defines a continuous map $\Omega^n_S \to \BdR$ by the universality of colimit.
\end{dfn}

By the definition, if we define an $S$-compatible integration, then we obtain a map $\rH^0(X,\Omega^n_X) \to \Omega^n_S \to \BdR$ associated to $\gamma$. We will construct an $S$-compatible integration in \S \ref{Integration along Rigid Analytic Cycles}, and regard the resulting map $\rH^0(X,\Omega_X)^n \to \BdR$ as the integration along $\gamma$.

\subsection{Integration along Rigid Analytic Cycles}
\label{Integration along Rigid Analytic Cycles}

From now on, we only consider the case where the coefficient group $M$ of the homology is a $\Z$-submodule of $\BdR$. We follow the convention in \S \ref{Analytic Singular Homology}. Let $\theta^k$ denote the reduction map $W(R)_k \twoheadrightarrow R$, and $\theta^{\BdR}$ the reduction map $\BdR^{+} \twoheadrightarrow \cC$. 

\vs
We denote by $[\bullet] \colon R \hookrightarrow W(R)$ the \Teichmuller lifting. For each $\chi \in \Frac(R)^{\times}$, put
\be
\cN_R(\chi) \coloneqq 
\left\{
\begin{array}{ll}
\alpha(0)^{-1} [\alpha] & (\alpha \in R) \\
\alpha(0) [\alpha^{-1}]^{-1} & (\alpha \notin R)
\end{array}
\right.
 \in 1 + \Fil^1 \BdR.
\ee
Since the \Teichmuller lifting is a $G_k$-equivariant monoid homomorphism, the map
\be
\cN_R \colon \Frac(R)^{\times} \to 1 + \Fil^1 \BdR, \ \alpha \mapsto \cN_R(\alpha)
\ee
is a $G_k$-equivariant group homomorphism.

\vs
For any $x \in 1+\Fil^1 \BdR$, the infinite sum $\log x \coloneqq - \sum_{a = 1}^{\infty} \frac{(1-x)^a}{a}$ uniquely converges in $\Fil^1 \BdR$ because it is the limit of a sequence of finite sums whose image in the Hausdorff space $\BdR^{+}/\Fil^{a} \BdR^{+}$ is an eventually constant sequence for any $a \in \N$. The logarithm
\be
\log \colon 1 + \Fil^1 \BdR \to \Fil^1 \BdR, \ x \mapsto - \sum_{a = 1}^{\infty} \frac{(1-x)^a}{a}
\ee
is a $G_k$-equivariant group homomorphism with respect to the multiplication of $1 + \Fil^1 \BdR$ and the addition of $\BdR$, and is continuous with respect to the relative topologies as subsets of $\BdR^{+}$, because the composite of $\log$ and the canonical projection $\BdR^{+} \twoheadrightarrow \BdR^{+}/\Fil^a \BdR^{+}$ is a polynomial function for any $a \in \N$. Since $\cN_R$ and $\log$ are $G_k$-equivariant group homomorphisms, so is the composite
\be
\log_R \coloneqq \log \circ \cN_R \colon R \setminus \ens{0} \to \Fil^1 \BdR.
\ee
For a $\xi \in \BdR$, we formally define $\xi^0 \coloneqq 1 \in \BdR$. In order to use powers of values of $\log_R$ as denominators of fractions which will appear in a formulation of an integration of a differential form, we characterise zeros of them.

\begin{prp}
\label{zero of log}
For any $(\alpha,c) \in \Frac(R)^{\times} \times \N$ with $\alpha(0)^{p^{\nu}} \in k^{\times}$ for some $\nu \in \N$, the following are equivalent:
\begin{itemize}
\item[(1)] The equality $(\log_R \alpha)^c = 0$ holds.
\item[(2)] The equality $(\log_R \alpha^{\ell_k})^c = 0$ holds.
\item[(3)] The inequality $c \neq 0$ and the equality $\alpha^{\ell_k} = 1$ hold.
\end{itemize}
\end{prp}

Although it is not so difficult to prove Proposition \ref{zero of log} directly, we will instead later prove a quite stronger result in Proposition \ref{valuation of log} in \S \ref{Appendix}.

\vs
Let $\rho \in D^n_k$ and Let $(i,j) \in [n]^2$. We denote by $\rho_{i \to j}$ the sequence
\be
(\rho(\ell_k \delta_{i/[n]})^{-1} \rho(\ell_k((1 - p^{- \nu}) \delta_{i/[n]} + p^{- \nu} \delta_{j/[n]})))_{\nu \in \N} \in (\cC^{\times})^{\N}.
\ee
By $\rho \in D^n_k$, there exists a $\chi \in (\cC/k)^{n \flat \times}$ with $\rho = \chi_{\Delta^n_k}$. For any $\nu \in \N$, we have
\be
& & \rho_{i \to j}(\nu + 1)^p = \left( \frac{\rho(\ell_k((1 - p^{- (\nu+1)}) \delta_{i/[n]} + p^{- (\nu+1)} \delta_{j/[n]}))}{\rho(\ell_k \delta_{i/[n]})} \right)^p \\
& = & \left( \frac{\chi(\ell_k((1 - p^{- (\nu+1)}) \delta_{i/[n]} + p^{- (\nu+1)} \delta_{j/[n]}))}{\chi(\ell_k \delta_{i/[n]})} \right)^p = \chi(- p^{- (\nu+1)} \delta_{i/[n]} + p^{- (\nu+1)} \delta_{j/[n]})^{\ell_k p} \\
& = & \chi(- p^{- \nu} \delta_{i/[n]} + p^{- \nu} \delta_{j/[n]})^{\ell_k} = \frac{\chi(\ell_k((1 - p^{- \nu}) \delta_{i/[n]} + p^{- \nu} \delta_{j/[n]}))}{\chi(\ell_k \delta_{i/[n]})} \\
& = & \frac{\rho(\ell_k((1 - p^{- \nu}) \delta_{i/[n]} + p^{- \nu} \delta_{j/[n]}))}{\rho(\ell_k \delta_{i/[n]})} = \rho_{i \to j}(\nu).
\ee
Therefore, we obtain $\rho_{i \to j} \in \Frac(R)^{\times}$. In particular, $\log_R \rho_{i \to j}$ makes sense.

\vs
We denote by $P_{\rho}(n)$ the quotient set of $[n]$ with respect to the equivalence relation $i \sim j$ on $(i,j) \in [n]^2$ defined as $\rho_{i \to j} = 1$, which is equivalent to $(\chi(p^{- \nu} \delta_{i/[n]})^{\ell_k})_{\nu \in \N} = (\chi(p^{- \nu} \delta_{j/[n]})^{\ell_k})_{\nu \in \N}$. For each $S \in P_{\rho}(n)$, we denote by $\rho_S \in K^{\times}$ the value $\rho(\ell_k \delta_{i/[n]}) = \chi(\delta_{i/[n]})^{\ell_k}$ for an $i \in S$, which depends only on $S$ by the equivalence above. For each $(S,S') \in P_{\rho}(n)^2$, we denote by $\rho_{S \to S'} \in \Fil^1 \BdR$ the value $\log_R \rho_{i \to j}$ for an $(i,j) \in S \times S'$, which depends only on $(S,S')$ by Proposition \ref{zero of log}.

\vs
Let $S$ and $T$ be finite totally ordered sets. We call an injective map $S \hookrightarrow T$ an {\it $S$-indexed permutation of $T$}. For an $S$-indexed permutation $\sigma$ of $T$, we define its {\it inversion number} $\inv(\sigma)$ the cardinal of the finite set $\set{s \in S^2}{s(0) < s(1) \land \sigma(s(0)) > \sigma((s(1))}$, its {\it signature} $\sgn(\sigma)$ as $(-1)^{\inv(\sigma)}$, and its {\it complement} $\sigma^{\r{c}}$ the unique element of $T \setminus \im(S)$ when $\# S + 1 = \# T$. We denote by $\cS(S,T)$ the set of $S$-permutations of $T$.

\vs
For a $[n-1]$-permutation $\sigma$ of $[n]$, we define its {\it modified signature} $\ol{\sgn}(\sigma)$ as $(-1)^{\sigma^{\r{c}}} \sgn(\sigma)$. For each $m \in [n]$, we denote by $s_{n,m}$ the $m$-th face map $[n-1] \hookrightarrow [n]$. Then, for any $m \in [n]$, the $m$-th face map $s_{n,m}$ is a $[n-1]$-permutation of $[n]$ of inversion number $0$, of signature $1$, of complement $m$, and of modified signature $(-1)^m$.

\begin{dfn}
\label{integration of a character}
For a $\rho \in D^n_k$, we define
\be
\int_n(\rho) \coloneqq \sum_{S \in P_{\rho}(n)} \frac{\rho_S}{\prod_{S' \in P_{\rho}(n) \setminus \ens{S}} \rho_{S \to S'}^{\# S'}} \in \BdR,
\ee
where the denominator of the fractional expression in the right hand is not $0$ by Proposition \ref{zero of log}. For a $(\rho,\sigma) \in D^n_k \times \cS([n-1],[n])$, we put $\int_n(\rho,\sigma) \coloneqq \ol{\sgn}(\sigma) \int_n(\rho)$.
\end{dfn}

The functional $\int_n(\rho,\sigma)$ plays a role of the integration of ``the differential $n$-form $\rho d x_{\sigma(0)} \wedge \cdots \wedge d x_{\sigma(n-1)}$'', although the expression does not rigorously make sense as a differential $n$-form on $\Delta_{K/k}^n$ due to the lack of the coordinate $x_0, \ldots, x_n$ in the ring $\cA_{K/k}^n$.

\vs
Let $K$ be a closed subextension of $\cC/k$ as in \S \ref{Perfectoid Algebra Associated to Polytope}. In order to apply Proposition \ref{kernel of comparison}, assume $\set{x \in \cC}{\exists r \in \N, x^{p^r} \in k} \subset K$. In order to apply Definition \ref{compatible integration}, fix a $k$-algebra embedding $K \hookrightarrow \BdR^{+}$ compatible with the reduction $\BdR^{+} \to \cC$, through which we regard $\BdR$ as a $K$-algebra.

\vs
For any $\sigma \in \cS([n-1],[n])$, the $K$-linear homomorphism
\be
K[D^n_k] \to \BdR, \ f \mapsto \sum_{\rho \in D^n_k} f(\rho) \int_n(\rho,\sigma)
\ee
uniquely factors through the composite of the completion $K[D^n_k] \to \cA_K(D^n_k)$ and the canonical projection $\pr^n_{K/k} \colon \cA_K(D^n_k) \to A_{K/k}^n$ by Lemma \ref{coefficient reconstruction} (2). We denote by $\ol{K[D^n_k]}$ the image of $K[D^n_k]$ in $A_{K/k}^n$, and by $\ol{\int_n(-,\sigma)}$ the induced $K$-linear homomorphism $\ol{K[D^n_k]} \to \BdR$. Unfortunately, $\ol{\int_n(-,\sigma)}$ is not continuous with respect to the relative topology of $\ol{K[D^n_k]}$ unless $n = 0$, and hence it is unable to continuously extend it to a wider space. Also, its graph is not closable, i.e. the closure of its graph does not define a map. Therefore, we directly make use of $\int_n(-,\sigma)$ without extending $\ol{\int_n(-,\sigma)}$, in the following construction of the integration of a differential form on $\Delta_{K/k}^n$.

\vs
Following the setting in \S \ref{Analytic Standard Simplex}, we set $\Lambda \coloneqq \Z[p^{-1}]$. For a $\tau \in \cS([n-1],[n])$, we denote by $\int d^{\tau}$ the map
\be
(D^n_k)^{[n]} \to \BdR, \ \rho \mapsto n^{-1} \sum_{\sigma \in \cS([n-1],[n])} \ol{\sgn}(\sigma) \prod_{h \in [n-1]} \log_R \rho(\tau(h))_{\sigma^{\r{c}} \to \sigma(h)} \int_n \left( \prod_{h=0}^{n} \rho(h) \right).
\ee
For a $(k,n)$-exponential datum $F$ (cf.\ Definition \ref{rigid datum}) and $(\rho,\tau) \in F^{[n]} \times \cS([n-1],[n])$, we denote by $d^{\tau} \rho$ the differential $n$-form on $\Spa(K \ens{F}_{\ad})$ defined as $\rho(\tau^{\r{c}}) (\bigwedge_{h \in [n-1]} d \rho(\tau(h)))$, where the exterior product is taken in the ascending order on the index $h$. The value $\int d^{\tau}(\rho)$ for a $\rho \in (D^n_k)^{[n]}$ plays a role of the integration of the differential $n$-form $d^{\tau} \rho$ on $\Spa(K \ens{F}_{\ad})$ for a $(k,n)$-exponential datum $F$ containing $\im(\rho)$. Although we only need to consider the case $\tau = s_{n,0}$, the formulation using the modified signature of a general permutation make the argument combinatorially natural.

\vs
For a $\tau \in \cS([n-1],[n])$ and an $\alpha \in ((\Frac(R)^{\times})^{[n]^2}$, we set
\be
D^{\tau}(\alpha) \coloneqq \sum_{\sigma \in \cS([n-1],[n])} \ol{\sgn}(\sigma) \prod_{h \in [n-1]} \log_R \alpha(\tau(h),\sigma(h)).
\ee
Another combinatorial description of $\int d^{\tau}(\rho)$ using $D^{\tau}$ justifies the scaling $n^{-1}$.

\begin{prp}
\label{signature and functional}
Let $(\rho,\tau) \in (D^n_k)^{[n]} \times \cS([n-1],[n])$. Fix a $\chi \in ((\cC/k)^{(n+1) \flat \times})^{[n]}$ such that $\rho(h)$ is given as $\chi(h)_{\Delta^n_{\Lambda}}$ for any $h \in [n]$. Let $\alpha \in (\Frac(R)^{\times})^{[n]^2}$ denote the map assigning to each $(h,j) \in [n]^2$ the element of $\Frac(R)^{\times}$ corresponding to the group homomorphism
\be
\Lambda \to \cC^{\times}, \ a \mapsto \chi(h)(a \delta_{j/[n]}),
\ee
where the characteristic function $\delta_{j/[n]}$ of $\ens{j} \subset [n]$ is regarded as the $\Lambda$-valued one. Then, the following equality holds:
\be
\int d^{\tau}(\rho) = D^{\tau}(\alpha) \int_n(\rho)
\ee
\end{prp}

\begin{proof}
By the definition, $\rho(h)_{i \to j}$ coincides with $\alpha(h,j) \alpha(h,i)^{-1} \in \Frac(R)^{\times}$, and hence $\log_R \rho(h)_{i \to j}$ coincides with $\log_R \alpha(h,j) - \log_R \alpha(h.i)$ for any $(h,i,j) \in [n]^3$. This implies that the sum
\be
\sum_{\sigma \in \cS([n-1],[n])} \ol{\sgn}(\sigma) \prod_{h \in [n-1]} \log_R \rho(\tau(h))_{\sigma^{\r{c}} \to \sigma(h)}
\ee
is a linear combination of $\prod_{h \in [n-1]} \log_R \alpha(\tau(h),\sigma(h))$, where $\sigma$ runs through all maps $[n-1] \to [n]$ given by replacing (possibly no) values of a $[n-1]$-permutation of $[n]$ by its complement.

\vs
If the number $c$ of the occurrence of the complement as values of such a $\sigma$ is greater than $1$, then the coefficient of the corresponding product is $0$, because it is a multiple of $\sum_{g \in \cS([c-1],[c-1])} \sgn(g) = 0$. Therefore, it suffices only to consider the case $0 \leq c \leq 1$. In that case, $\sigma$ is again a $[n-1]$-permutation of $[n]$. Therefore, the sum is a linear combination of $\prod_{h \in [n-1]} \log_R \alpha(\tau(h),\sigma(h))$, where $\sigma$ runs through all $[n-1]$-permutation of $[n]$. The coefficient of such a $\sigma$ is its modified signature minus the sum of the modified signatures of the $[n-1]$-permutations of $[n]$ obtained by replacing one value of $\sigma$ by its complement, and is equal to
\be
& & \ol{\sgn}(\sigma) - \sum_{h \in [n-1]} (-1)^{\sigma(h) + \# \set{h' \in [n-1]}{(\sigma(h') > \sigma(h)) \not\Leftrightarrow (\sigma(h') > \sigma^{\r{c}})}} \sgn(\sigma) \\
& = & \sgn(\sigma) \left( (-1)^{\sigma^{\r{c}}} - \sum_{h \in [n-1]} (-1)^{\sigma(h) + \v{\sigma(h) - \sigma^{\r{c}}} - 1}  \right) = \sgn(\sigma) \left( (-1)^{\sigma^{\r{c}}} + \sum_{h' \in \im(\sigma)} (-1)^{\sigma^{\r{c}}} \right)  = n \ol{\sgn}(\sigma).
\ee
This implies the assertion.
\end{proof}

In order to define the integration of a differential $n$-form on the analytic standard $n$-simplex, we introduce the notion of integrability of an rigid datum.

\begin{dfn}
\label{integrable datum}
A $(k,n)$-exponential datum $F$ is said to be a {\it $(K/k,n)$-integrable datum} if there exists a unique continuous $K$-linear homomorphism
\be
I_F \colon \Omega^n_F \to \BdR
\ee
satisfying
\be
I_F(d^{\tau} \rho) = \int d^{\tau}(\rho)
\ee
for any $(\rho,\tau) \in F^{[n]} \times \cS([n-1],[n])$. We denote by $\ID(K/k,n) \subset \ED(k,n)$ the subset of $(K/k,n)$-integrable data.
\end{dfn}

\begin{exm}
Let $u$ be an non-torsion element of $1 + k_{<1}$, and take a $k$-exponential map $\rho \colon D^1_k \to \cC$ such that $\rho(\ell_k \delta_{0/[1]}) = 1$ and $\rho(\ell_k \delta_{1/[1]}) = u$. For any $(i,j) \in \N^2 \setminus \ens{(0,0)}$, we have $P_{\rho^{i+j}}(1) = \ens{\ens{0},\ens{1}}$, and hence
\be
& & \int d^{s_{1,0}}(\rho^i,\rho^j) = (\log_R \rho) \times (-1)^0 \left( \frac{(\rho^{i+j})_{\ens{0}}}{(\rho^{i+j})_{\ens{0} \to \ens{1}}^1} + \frac{(\rho^{i+j})_{\ens{1}}}{(\rho^{i+j})_{\ens{1} \to \ens{0}}^1} \right) \\
& = & (\log_R \rho_{0 \to 1}) \left( \frac{\rho^{i+j}(\ell_k \delta_{0/[1]})}{\log_R (\rho^{i+j})_{0 \to 1}} + \frac{\rho^{i+j}(\ell_k \delta_{1/[1]})}{\log_R (\rho^{i+j})_{1 \to 0}} \right) \\
& = & (\log_R \rho_{0 \to 1}) \left( \frac{1}{(i+j) \log_R \rho_{0 \to 1}} + \frac{u^{(i+j)\ell_k}}{-(i+j) \log_R \rho_{0 \to 1}} \right) = \frac{1 - u^{(i+j)\ell_k}}{i+j}.
\ee
by $P_{\rho^{i+j}}(1) = \ens{\ens{0},\ens{1}}$. The right hand side converges to $\ell_k \log u$ as $i+j \to \infty$, and hence belongs to a bounded subset of $k$ independent of $(i,j)$. This implies that the $K$-linear homomorphism
\be
K[\rho] d \rho \to \BdR, \ \rho^i d \rho \mapsto \int d^{s_{1,0}}(\rho^i,\rho)
\ee
uniquely continuously extends to $\Omega^n_{\ens{\rho}} \cong K \ens{\rho} d \rho$ (with respect to the topology of $\BdR$ explained in \S \ref{Rigid Analytic Singular Simplex} rather than the valuation topology). Therefore, $\ens{\rho}$ is a $(K/k,n)$-integrable datum.
\end{exm}

Now we are ready to define an integration of a differential form. In Definition \ref{integrable datum}, the condition on an integrable datum guarantees that $(I_F)_{F \in \ID(K/k,n)}$ forms an $\ID(K/k,n)$-compatible integration (Definition \ref{compatible integration}). Indeed, for any $(F_0,F_1) \in \ID(K/k,n)^2$ with $F_0 \subset F_1$, the composite of $I_{F_1}$ and the homomorphism $\Omega^n_{F_0} \to \Omega^n_{F_1}$ associated to the inclusion $F_0 \hookrightarrow F_1$ satisfies the condition on $I_{F_0}$, and hence coincides with $I_{F_0}$ by the uniqueness.

\begin{dfn}
We put $\Omega^n_{\Delta_{K/k}^n} \coloneqq \Omega^n_{\ID(K/k,n)}$ (Definition \ref{analytic differenial form}), and denote by
\be
\int_{\Delta_{K/k}^n} \colon \Omega^n_{\Delta_{K/k}^n} \to \BdR,
\ee
the continuous map induced by the $\ID(K/k,n)$-compatible integration $(I_F)_{F \in \ID(K/k,n)}$.
\end{dfn}

The map $\int_{\Delta_{K/k}^n}$ plays the role of an integration of a differential form on $\Delta_{K/k}^n$. We note that $\ID(K/k,n)$ is not directed unless $n = 0$, and hence the colimit is not replaced by the colimit as topological $K$-modules. This obstruction prevents us to formulate an integration of a differential form on an adic space along a fixed analytic $n$-simplex in a well-defined way even if we assume that the simplex factors through some (not unique) $(K/k,n)$-integrable datum. Therefore, we introduce a narrower class of simplexes in order to ensure the well-definedness of integration.

\begin{dfn}
A partially ordered set $P$ is said to be a {\it minimally directed set} if every pair $(p,q) \in P^2$ of minimal elements admits a common upperbound.
\end{dfn}

Let $X$ be an adic space over $\Spa(k_{\ad})$. Using the notion of minimal directedness, we introduce a narrow class of analytic simplexes on $X$. Let $\gamma$ be an $\ID(K/k,n)$-rigid analytic singular $n$-simplex on $X$ over $K$. We show that the minimal directedness ensures the well-definedness of integration of a differential form.

\begin{prp}
\label{well-definedness of integration}
Let $P$ be a minimally directed subset of $P_{\gamma}$. Let $(F_0,F_1) \in P^2$. For each $i \in 2$, denote by $\gamma^{*,F_i}$ denote the composite $e_i \circ e_i$ of the homomorphism
\be
d_i \colon \rH^0(X,\Omega^n_X) \to \Omega^n_{F_i}
\ee
associated to the factorisation $\Spa(K \ens{F_i}_{\ad}) \to X$ of $\gamma$ and the canonical map
\be
e_i \colon \Omega^n_{F_i} \to \Omega^n_{\Delta_{K/k}^n}.
\ee
Then $\gamma^{*,F_0} = \gamma^{*,F_1}$ holds.
\end{prp}

\begin{proof}
First, assume $F_0 \subset F_1$.  Then the diagram
\be
\xymatrix{
\rH^0(X,\Omega^n_X) \ar[d]^-{d_0} \ar[r]^-{d_1} & \Omega^n_{F_1} \ar[d]^-{e_1} \\
\Omega^n_{F_0} \ar[r]^-{e_0} \ar[ru] & \Omega^n_{\Delta_{K/k}^n}
}
\ee
commutes with respect to the homomorphism $\Omega^n_{F_0} \to \Omega^n_{F_1}$ associated to the inclusion $F_0 \hookrightarrow F_1$. This implies $\gamma^{*,F_0} = \gamma^{*,F_1}$.

\vs
Next, consider the general case. Since $P$ is a set of finite sets partially ordered by the inclusion, every element has a minimal lower bound. Let $F'_0$ and $F'_1$ be minimal lowerbound of $F_0$ and $F_1$ in $P$ respectively. Since $P$ is minimally directed, $F'_0$ and $F'_1$ admit a common upperbound $F$. Therefore, the argument above implies $\gamma^{*,F_0} = \gamma^{*,F'_0} = \gamma^{*,F} = \gamma^{*,F'_1} = \gamma^{*,F_1}$.
\end{proof}

When $n > 0$, for each $m \in [n]$, we abuse the notation $s_{n,m}$ of the $m$-th face map $[n-1] \hookrightarrow [n]$ also for the $\Lambda$-linear morphism $D^{n-1}_k \hookrightarrow D^n_k$ given as the $m$-th face map in the simplicial sense, denote by $\rho_{\check{m}} \in D_{D^{n-1}_k,K}$ the composite $\rho \circ s_{n,m}$ for a $\rho \in D_{\ell_k \Delta^n_{\Lambda},n}$ (cf.\ Proposition \ref{functoriality of the analytification} (1)), by $F_{\check{m}} \subset D_{D^{n-1}_k,K}$ the image $\set{\rho_{\check{m}}}{\rho \in F}$ of $F$ for an $F \in \ED(k,n)$, and by $\partial_m \gamma \colon \Delta_{K/k}^{n-1} \to \Delta_{K/k}^n$ the composite $\gamma \circ (s_{n,m})_{K/k}$ (cf.\ Proposition \ref{functoriality of the analytification}).

\begin{dfn}
\label{integrable analytic simplex}
Recursively on $n$, we define a subset $P'(\gamma)$ of $P(\gamma)$ as $P(\gamma)$ itself when $n = 0$ and otherwise as the set consisting of $F \in P(\gamma)$ such that $F_{\check{m}} \in P'(\partial_m \gamma)$ for any $m \in [n]$, and say that $\gamma$ is an {\it integrable $(K/k,n)$-simplex on $X$} if $P'(\gamma)$ is minimally directed and in addition $\partial_m \gamma$ is an integrable $(K/k,n-1)$-simplex for any $m \in [n]$ when $n > 0$.
\end{dfn}

By the recursive definition, the notion of the integrability of an analytic simplex is preserved under taking a face. Therefore, the set of integrable simplexes defines a subcomplex of the chain complex associated to analytic simplexes.

\begin{dfn}
\label{integrable analytic singular homology}
Let $M$ be an Abelian group. We denote by $C_{K/k,M}^{\dagger}$ the subfunctor of $C_{K/k,M}$ defined by the direct sums of $M$ indexed by the set of integrable simplexes, and abbreviate $\rH_*(\bullet,C_{K/k,M}^{\dagger})$ and $\rH^*(\bullet,C_{K/k,M}^{\dagger})$ to $\rH_*^{\dagger}(\bullet,K/k,M)$ and $\rH^*_{\dagger}(\bullet,K/k,M)$ respectively.
\end{dfn}

By the argument in \S \ref{Analytic Singular Homology}, $\rH_*^{\dagger}(\bullet,K/k,M)$ and $\rH^*_{\dagger}(\bullet,K/k,M)$ form Galois representation. In addition, we will define an integration of a closed form along a cycle. For this purpose, we first define an integration of a differential form along a simplex.

\begin{dfn}
\label{integration along a simplex}
Suppose that $\gamma$ is an integrable $(K/k,n)$-simplex on $X$. Then $\set{\gamma^{*,F}}{F \in P'(\gamma)}$ forms a singleton by Proposition\ref{well-definedness of integration}. We denote by $\gamma^*$ its unique element. For an $\omega \in \rH^0(X,\Omega^n_X)$, we define
\be
\int_{\gamma} \omega \coloneqq \int_{\Delta_{K/k,n}} \gamma^* \omega \in \BdR,
\ee
and call it the integral of $\omega$ along $\gamma$.
\end{dfn}

In order to formulate an integration of a closed form along a cycle, we need to show a counterpart of Stokes theorem.

\begin{thm}
\label{Stokes}
Suppose $n > 0$. Then the equality
\be
\int_{\gamma} d \omega = \sum_{m \in [n]} (-1)^m \int_{\partial_m \gamma} \omega
\ee
holds for any differential $(n-1)$-form $\omega$ on $X$ and any integrable $(K/k,n)$-simplex $\gamma$ on $X$.
\end{thm}

In order to show Theorem \ref{Stokes}, we prepare two lemmata.

\begin{lmm}
\label{log combinatorial}
Suppose $n > 0$. Let $\rho \in D^n_k$ and $(i,j) \in [n]^2$ with $i \neq j$. Then, the following equality holds:
\be
\log_R \rho_{j \to i} \int_n(\rho) = \int_{n-1}(\rho_{\check{j}} - \rho_{\check{i}})
\ee
\end{lmm}

\begin{proof}
Let $m \in [n]$. We denote by $S_m \in P_{\rho}(n)$ the equivalence classes of $m$. Put $Q_m \coloneqq P_{\rho_{\check{m}}}(n) \cup \ens{\emptyset}$, $P_m \coloneqq P_{\rho}(n) \setminus \ens{S_m}$, $S'_m \coloneqq S_m \setminus \ens{m}$. By the definition, we have
\be
P_{\rho}(n) \cap Q_m = P_m, \ P_{\rho}(n) = P_m \cup \ens{S_j}, \ Q_m = P_m \cup \ens{S'_j}.
\ee
We formally define $(\rho_{\check{m}})_{\emptyset \to S'} = (\rho_{\check{m}})_{S' \to \emptyset} \coloneqq 1 \in \BdR$ for any $S' \in Q_m \setminus \ens{\emptyset}$. 

\vs
Put $P' \coloneqq P_i \cap P_j$. We formally define $\rho_{\emptyset} = 0$. First, we assume $S_i = S_j$. Then, by Proposition \ref{zero of log}, we have
\be
\log_R \rho_{j \to i} \int_n(\rho) = 0 \int_n(\rho) = 0, \ \int_{n-1}(\rho_{\check{i}}) = \int_{n-1}(\rho_{\check{j}})
\ee
and hence the equality holds. Next, we assume $S_i \neq S_j$. Since $\log_R$ is a group homomorphism, we have
\be
& & \log_R \rho_{j \to i} \int_n(\rho) = \log_R \rho_{j \to i} \sum_{S \in P_{\rho}(n)} \frac{\rho_S}{\prod_{S' \in P_{\rho}(n) \setminus \ens{S}} \rho_{S \to S'}^{\# S'}} \\
& = & \log_R \rho_{j \to i} \sum_{S \in P'} \frac{\rho_S}{\rho_{S \to S_i}^{\# S_i} \rho_{S \to S_j}^{\# S_j} \prod_{S' \in P' \setminus \ens{S}} \rho_{S \to S'}^{\# S'}} + \log_R \rho_{j \to i} \frac{\ell_k^{\# S_i} \rho_{S_i}}{\rho_{S_i \to S_j}^{\# S_j} \prod_{S' \in P'} \rho_{S_i \to S'}^{\# S'}} \\
& & + \log_R \rho_{j \to i} \frac{\ell_k^{\# S_j} \rho_{S_j}}{\rho_{S_j \to S_i}^{\# S_i} \prod_{S' \in P'} \rho_{S_j \to S'}^{\# S'}} \\
& = & \sum_{S \in P'} \frac{\rho_S(\rho_{S \to S_i} - \rho_{S \to S_j})}{\rho_{S \to S_i}^{\# S_i} \rho_{S \to S_j}^{\# S_j} \prod_{S' \in P' \setminus \ens{S}} \rho_{S \to S'}^{\# S'}} + \frac{\ell_k^{\# S_i} \rho_{S_i}}{\rho_{S_i \to S_j}^{\# S_j - 1} \prod_{S' \in P'} \rho_{S_i \to S'}^{\# S'}} \\
& & + \frac{\ell_k^{\# S_j} \rho_{S_j}}{\rho_{S_j \to S_i}^{\# S_i - 1} \prod_{S' \in P'} \rho_{S_j \to S'}^{\# S'}} \\
& = & \sum_{S \in P'} \frac{\rho_S}{\rho_{S \to S_i}^{\# S_i - 1} \rho_{S \to S_j}^{\# S_j} \prod_{S' \in P' \setminus \ens{S}} \rho_{S \to S'}^{\# S'}} + \sum_{S \in P'} \frac{\rho_S}{\rho_{S \to S_i}^{\# S_i} \rho_{S \to S_j}^{\# S_j - 1} \prod_{S' \in P' \setminus \ens{S}} \rho_{S \to S'}^{\# S'}} \\
& & + \frac{\ell_k^{\# S_i} \rho_{S_i}}{\rho_{S_i \to S_j}^{\# S_j - 1} \prod_{S' \in P'} \rho_{S_i \to S'}^{\# S'}} + \frac{\ell_k^{\# S_j} \rho_{S_j}}{\rho_{S_j \to S_i}^{\# S_i - 1} \prod_{S' \in P'} \rho_{S_j \to S'}^{\# S'}} \\
& = & \sum_{S \in P'} \frac{\rho_S}{\rho_{S \to S_i}^{\# S'_i} \rho_{S \to S_j}^{\# S_j} \prod_{S' \in P' \setminus \ens{S}} \rho_{S \to S'}^{\# S'}} + \frac{\ell_k^{\# S_j} \rho_{S_j}}{\rho_{S_j \to S_i}^{\# S'_i} \prod_{S' \in P'} \rho_{S_j \to S'}^{\# S'}} \\
& & + \sum_{S \in P'} \frac{\rho_S}{\rho_{S \to S_i}^{\# S_i} \rho_{S \to S_j}^{\# S'_j} \prod_{S' \in P' \setminus \ens{S}} \rho_{S \to S'}^{\# S'}} + \frac{\ell_k^{\# S_i} \rho_{S_i}}{\rho_{S_i \to S_j}^{\# S'_j} \prod_{S' \in P'} \rho_{S_i \to S'}^{\# S'}} \\
& = & \sum_{S \in P'} \frac{\rho_S}{\prod_{S' \in Q_i \setminus \ens{S}} \rho_{S \to S'}^{\# S'}} + \frac{\ell_k^{\# S_j} \rho_{S_j}}{\prod_{S' \in Q_i \setminus \ens{S_j}} \rho_{S_j \to S'}^{\# S'}} \\
& & + \sum_{S \in P'} \frac{\rho_S}{\prod_{S' \in Q_j \setminus \ens{S}} \rho_{S \to S'}^{\# S'}} + \frac{\ell_k^{\# S_i} \rho_{S_i}}{\prod_{S' \in Q_j \setminus \ens{S_i}} \rho_{S_i \to S'}^{\# S'}} \\
& = & \sum_{S \in Q_i} \frac{\rho_S}{\prod_{S' \in Q_i \setminus \ens{S}} \rho_{S \to S'}^{\# S'}} + \sum_{S \in Q_j} \frac{\rho_S}{\prod_{S' \in Q_j \setminus \ens{S}} \rho_{S \to S'}^{\# S'}} \\
& = & \sum_{S \in P_{\rho_{\check{i}}}(n-1)} \frac{(\rho_{\check{i}})_S}{\prod_{S' \in P_{\rho_{\check{i}}}(n-1) \setminus \ens{S}} (\rho_{\check{i}})_{S \to S'}^{\# S'}} + \sum_{S \in P_{\rho_{\check{j}}}(n-1)} \frac{(\rho_{\check{j}})_S}{\prod_{S' \in P_{\rho_{\check{j}}}(n-1) \setminus \ens{S}} (\rho_{\check{j}})_{S \to S'}^{\# S'}} \\
& = & \int_{n-1}(\rho_{\check{i}}) + \int_{n-1}(\rho_{\check{j}}).
\ee
\end{proof}

\begin{lmm}
\label{log combinatorial 2}
Suppose $n > 0$. Let $\rho \in (D^n_k)^{[n-1]}$ and $\tau \in \cS([n-2],[n-1])$. Put $\rho^{\Pi} = \prod_{h \in [n-1]} \rho(h)$. Define $\alpha \in (\Frac(R)^{\times})^{[n]^2}$ as in Proposition \ref{signature and functional}. For a $\sigma \in \cS([n-2],[n])$, let $\sigma^{\r{c},-}$ denote the minimum of $[n] \setminus \im(\sigma)$, and $\sigma^{\r{c},+}$ denote the maximum of $[n] \setminus \im(\sigma)$. Then, the following equality holds:
\be
& & \sum_{\sigma \in \cS([n-2],[n])} (-1)^{\sigma^{\r{c},-} + \sigma^{\r{c},+} + n} \sgn(\sigma) \left( \log_R \alpha(\tau^{\r{c}},\sigma^{\r{c},+}) - \log_R \alpha(\tau^{\r{c}},\sigma^{\r{c},-}) \right) \prod_{h \in [n-2]} \log_R \alpha(\tau(h),\sigma(h)) \\
& = & \sum_{\sigma \in \cS([n-2],[n])} (-1)^{\sigma^{\r{c},-} + \sigma^{\r{c},+} + n} \sgn(\sigma) \log_R \rho^{\Pi}_{\sigma^{\r{c},-})  \to \sigma^{\r{c},+}} \prod_{h \in [n-2]} \log_R \alpha(\tau(h),\sigma(1+h))
\ee
\end{lmm}

\begin{proof}
For any $(\beta_0,\sigma) \in (\Frac(R)^{\times})^{[n]} \times \cS([n-2],[n])$, we have
\be
& & (-1)^{\sigma^{\r{c},-} + \sigma^{\r{c},+} + n} \sgn(\sigma) \left( \log_R \beta(\sigma^{\r{c},+}) - \log_R \beta(\sigma^{\r{c},-}) \right) \\
& = & (-1)^{\sigma^{\r{c},-} + n - \sigma^{\r{c},+}} \sgn(\sigma) \log_R \beta(\sigma^{\r{c},+}) + (-1)^{\sigma^{\r{c},-} + n - \sigma^{\r{c},-} - 1} \sgn(\sigma) \log_R \beta(\sigma^{\r{c},-}) \\
& = & (-1)^{\sigma^{\r{c},-} + \# \im(\sigma)_{> \sigma^{\r{c},+}}} \sgn(\sigma) \log_R \beta(\sigma^{\r{c},+}) + (-1)^{\sigma^{\r{c},-} + \# \im(\sigma)_{> \sigma^{\r{c},-}}} \sgn(\sigma) \log_R \beta(\sigma^{\r{c},-}) \\
& = & (-1)^{\sigma^{\r{c},-}} \sgn(\ens{(-1,\sigma^{\r{c},+})} \cup \sigma) \log_R \beta(\sigma^{\r{c},+}) + (-1)^{\sigma^{\r{c},-}} \sgn(\ens{(-1,\sigma^{\r{c},-})} \cup \sigma) \log_R \beta(\sigma^{\r{c},-}).
\ee
Let $\beta \in (\Frac(R)^{\times})^{[n-1] \times [n]}$. Applying the result above to $\beta(\tau^{\r{c}},-)$ and concatenating $(\sigma^{\r{c},-}) \in \cS([0],[n])$ or $(\sigma^{\r{c},+}) \in \cS([0],[n])$ and $\sigma \in \cS([n-2],[n] \setminus \ens{\sigma^{\r{c},-},\sigma^{\r{c},+}})$ to obtain a $[n-1]$-permutation of $[n]$, we obtain
 \be
& & \sum_{\sigma \in \cS([n-2],[n])}  (-1)^{\sigma^{\r{c},-} + \sigma^{\r{c},+} + n} \sgn(\sigma) \left( \log_R \beta(\sigma^{\r{c},+}) - \log_R \beta(\sigma^{\r{c},-}) \right) \prod_{h \in [n-2]} \log_R \beta(\tau(h),\sigma(h)) \\
& = & \sum_{\sigma \in \cS([n-1],[n])} (-1)^{\sigma^{\r{c}}} \sgn(\sigma) \log_R \alpha(\tau^{\r{c}},\sigma(0)) \prod_{h \in [n-2]} \log_R \beta(\tau(h),\sigma(1+h)) \\
& = & \sum_{\sigma \in \cS([n-1],[n])} \ol{\sgn}(\sigma) \log_R \alpha(\tau^{\r{c}},\sigma(0)) \prod_{h \in [n-2]} \log_R \beta(\tau(h),\sigma(1+h)).
\ee
We denote by $\ol{\tau} \in \cS([n-1],[n])$ the concatenation of $(\tau^{\r{c}}) \in \cS([0],[n-1])$ and $\tau \in \cS([n-2],[n-1] \setminus \ens{\tau^{\r{c}}}$, and by $\ol{\beta} \in (\Frac(R)^{\times})^{[n]^2}$ the concatenation of $\beta$ and $(1) \in (\Frac(R)^{\times})^{[0] \times [n]}$. We obtain that the right hand side coincides with
\be
\sum_{\sigma \in \cS([n-1],[n])} \ol{\sgn}(\sigma) \prod_{h \in [n-1]} \log_R \beta(\ol{\tau}(h),\sigma(h)) = \sum_{\sigma \in \cS([n-1],[n])} \ol{\sgn}(\sigma) \prod_{h \in [n-1]} \log_R \ol{\beta}(\ol{\tau}(h),\sigma(h)) = D^{\ol{\tau}}(\ol{\beta}).
\ee
Applying the result above to $\alpha$, we obtain
\be
\sum_{\sigma \in \cS([n-2],[n])} (-1)^{\sigma^{\r{c},-} + \sigma^{\r{c},+} + n} \sgn(\sigma) \left( \log_R \alpha(\tau^{\r{c}},\sigma^{\r{c},+}) - \log_R \alpha(\tau^{\r{c}},\sigma^{\r{c},-}) \right) \prod_{h \in [n-2]} \log_R \alpha(\tau(h),\sigma(h)) = D^{\ol{\tau}}(\ol{\alpha}).
\ee
On the other hand, for any $\sigma \in \cS([n-2],[n])$, since $\log_R$ is a group homomorphism, we have
\be
& & (-1)^{\sigma^{\r{c},-} + \sigma^{\r{c},+} + n} \sgn(\sigma) \log_R \rho^{\Pi}_{\sigma^{\r{c},-})  \to \sigma^{\r{c},+}} = (-1)^{\sigma^{\r{c},-} + \sigma^{\r{c},+} + n} \sgn(\sigma) \sum_{m \in [n-1]} \log_R \rho(m)_{\sigma^{\r{c},-})  \to \sigma^{\r{c},+}} \\
& = & \sum_{m \in [n-1]} (-1)^{\sigma^{\r{c},-} + \sigma^{\r{c},+} + n} \sgn(\sigma) \left( \log_R \alpha(m,\sigma^{\r{c},+}) - \log_R \alpha(m,\sigma^{\r{c},-}) \right).
\ee
Applying the first result to $\alpha(m,-)$, we obtain that the right hand side coincides with
\be
\sum_{m \in [n-1]} \left( (-1)^{\sigma^{\r{c},-}} \sgn(\ens{(-1,\sigma^{\r{c},+})} \cup \sigma) \log_R \alpha(m,\sigma^{\r{c},+}) - \log_R \alpha(m,\sigma^{\r{c},-}) \right).
\ee
Therefore, we obtain
\be
& & \sum_{\sigma \in \cS([n^2],[n])} (-1)^{\sigma^{\r{c},-} + \sigma^{\r{c},+} + n} \sgn(\sigma) \log_R \rho^{\Pi}_{\sigma^{\r{c},-})  \to \sigma^{\r{c},+}} \prod_{h \in [n-2]} \log_R \alpha(\tau(h),\sigma(h)) \\
& = & \sum_{\sigma \in \cS([n^2],[n])} \sum_{m \in [n-1]} \left( (-1)^{\sigma^{\r{c},-}} \sgn(\ens{(-1,\sigma^{\r{c},+})} \cup \sigma) \log_R \alpha(m,\sigma^{\r{c},+}) - \log_R \alpha(m,\sigma^{\r{c},-}) \right) \\
& & \prod_{h \in [n-2]} \log_R \alpha(\tau(h),\sigma(h)) \\
& = & \sum_{m \in [n-1]} \sum_{m_1 \in [n]} \sum_{m_2 \in [n] \setminus \ens{m_1}} \sum_{\sigma \in \cS([n-2],[n] \setminus \ens{m_1,m_2})} (-1)^{m_1} \sgn(\ens{(-1,m_2)} \cup \sigma) \log_R \alpha(m,m_2) \\
& & \prod_{h \in [n-2]} \log_R \alpha(\tau(h),\sigma(h)).
\ee
Concatenating $(m_2) \in \cS([0],[n])$ and $\sigma \in \cS([n-2],[n] \setminus \ens{m_1,m_2})$ to obtain a $[n-1]$-permutation of $[n]$ with complement $m_1$, we obtain that the right hand side coincides with
\be
& & \sum_{m \in [n-1]} \sum_{\sigma \in \cS([n-1],[n])} (-1)^{\sigma^{\r{c}}} \sgn(\sigma) \log_R \alpha(m,\sigma(0)) \prod_{h \in [n-2]} \log_R \alpha(\tau(h),\sigma(1+h)) \\
& = & \sum_{m \in [n-1]} \sum_{\sigma \in \cS([n-1],[n])} \ol{\sgn}(\sigma) \log_R \alpha(m,\sigma(0)) \prod_{h \in [n-2]} \log_R \alpha(\tau(h),\sigma(1+h)).
\ee
Denoting by $\tau_m \in [n]^{[n-1]}$ the concatenation of $(m) \in [n]^{[0]}$ and $\tau \in \cS([n-2],[n]) \subset [n]^{[n-2]}$ for each $m \in [n]$, we obtain that the right hand side coincides with
\be
& & \sum_{m \in [n-1]} \sum_{\sigma \in \cS([n-1],[n])} \ol{\sgn}(\sigma) \prod_{h \in [n-1]} \log_R \alpha(\tau_m(h),\sigma(h)) \\
& = & D^{\ol{\tau}}(\ol{\alpha}) + \sum_{m \in [n-1] \setminus \ens{\tau^{\r{c}}}} \sum_{\sigma \in \cS([n-1],[n])} \ol{\sgn}(\sigma) \prod_{h \in [n-1]} \log_R \alpha(\tau_m(h),\sigma(h)) \\
& = & D^{\ol{\tau}}(\ol{\alpha}) + \sum_{m \in \im(\tau)} \sum_{\sigma \in \cS([n-1],[n])} \ol{\sgn}(\sigma) \prod_{h \in [n-1]} \log_R \alpha(\tau_m(h),\sigma(h)) \\
& = & D^{\ol{\tau}}(\ol{\alpha}) + \sum_{h \in [n-2]} \sum_{\sigma \in \cS([n-1],[n])} \ol{\sgn}(\sigma) \prod_{h' \in [n-1]} \log_R \alpha(\tau_{\tau(h)}(h'),\sigma(h')).
\ee
The second term of the right hand side is $0$, because for any $(h,\sigma) \in [n-2] \times \cS([n-1],[n])$, denoting by $(0,h+1) \in \cS([n-1],[n-1])$ the swap of $0$ and $h+1$, we have
\be
& & \ol{\sgn}(\sigma) \prod_{h' \in [n-1]} \log_R \alpha(\tau_{\tau(h)}(h'),\sigma(h')) \\
& = & - \ol{\sgn}(\sigma \circ (0,h+1)) \prod_{h' \in [n-1]} \log_R \alpha(\tau_{\tau(h)}(h'),(\sigma \circ (0,h+1))(h')).
\ee
Therefore, the both hand side of the assertion coincides with $D^{\ol{\tau}}(\ol{\alpha})$.
\end{proof}

\begin{proof}[Proof of Theorem \ref{Stokes}.]
By the definition of the integration, it suffices to show
\be
I_F(d \omega) = \sum_{m \in [n]} (-1)^m I_{F_{\check{m}}}(\omega_{\check{m}})
\ee
for any $F \in \ID(K/k,n)$ and $\omega \in \Omega^{n-1}_F$, where $\omega_{\check{m}} \in \Omega^{n-1}_{F_{\check{m}}}$ denotes the image of $\omega$ by the homomorphism associated to the restriction map $F \to F_{\check{m}}$ given by $s_{n,m}^*$ (cf.\ Proposition \ref{functoriality of the analytification} (1)). Since inserting $1$ to $F$ does not affect $K \ens{F}$ and $K \ens{F_{\check{m}}}$, we may assume $1 \in F$. Since the image of $K[F]$ is dense in $K \ens{F}$, it is reduced to the case where $\omega$ is given as $d^{\tau} \rho$ for some $(\rho,\tau) \in F^{[n-1]} \times \cS([n-2],[n-1])$, by the continuity of $I_F$ and $I_{F_{\check{m}}}$ for any $m \in [n]$ and the Hausdorffness of $\BdR$.

\vs
We denote by $\ol{\rho} \in F^{[n]}$ the concatenation of $\rho \in F^{[n-1]}$ and $(1) \in F^{[0]}$, and put $\rho^{\Pi} \coloneqq \prod_{h=0}^{n-1} \rho(h)$. We define $\ol{\alpha} \in (\Frac(R)^{\times})^{[n]^2}$ by the construction of $\alpha$ in the assertion of Proposition \ref{signature and functional} applied to $\ol{\rho}$. We denote by $\alpha \in (\Frac(R))^{\times})^{[n-1] \times [n]}$ the $[n-1] \times [n]$-matrix given by removing the last row $(1)_{j \in [n]} \in (\Frac(R)^{\times})^{[n]}$ from $\ol{\alpha}$. By the defining property of $I_F$ and Proposition \ref{signature and functional}, we have
\be
& & I_F(d \omega) = I_F(d(d^{\tau} \rho)) = I_F \left( d \rho(\tau^{\r{c}}) \wedge \left( \bigwedge_{h \in [n-2]} d \rho(\tau(h)) \right) \right) = I_F \left( \bigwedge_{h \in [n-1]} d \rho(\ol{\tau}(h)) \right) \\
& = & I_F(d^{\ol{\tau}} \ol{\rho}) = \int d^{\ol{\tau}} (\ol{\rho}) = D^{\ol{\tau}}(\ol{\alpha}) \int_n \left( \prod_{h=0}^{n} \ol{\rho}(h) \right) \\
& = & \sum_{\sigma \in \cS([n-1],[n])} \ol{\sgn}(\sigma) \prod_{h \in [n-1]} \log_R \ol{\alpha}(\ol{\tau}(h),\sigma(h)) \int_n \left( \rho^{\Pi} \right) \\
& = & \sum_{\sigma \in \cS([n-1],[n])} \ol{\sgn}(\sigma) \prod_{h \in [n-1]} \log_R \alpha(\ol{\tau}(h),\sigma(h)) \int_n \left( \rho^{\Pi} \right) \\
& = & \sum_{\sigma \in \cS([n-1],[n])} \ol{\sgn}(\sigma) \left( \prod_{h \in [n-2]} \log_R \alpha(\tau(h),\sigma(1+h)) \right) \log_R \alpha(\tau^{\r{c}},\sigma(0)) \int_n \left( \rho^{\Pi} \right).
\ee
For an $(m_1,m_2) \in [n]^2$, we define $\delta_{m_1 > m_2}$ as $1$ if $m_1 > m_2$ and as $0$ otherwise. Reordering the sum by $(m_1,m_2) = (\sigma^{\r{c}},\sigma(n-1))$, we obtain that the right hand side coincides with
\be
& & \sum_{m_1 \in [n]} \sum_{m_2 \in [n] \setminus \ens{m_1}} \sum_{\sigma \in \cS([n-2],[n] \setminus \ens{m_1,m_2})} (-1)^{m_1 + \# \set{h \in [n^2]}{m_2 > \sigma(h)}} \sgn(\sigma) \\
& & \left( \prod_{h \in [n-2]} \log_R \alpha(\tau(h),\sigma(h)) \right) \log_R \alpha(\tau^{\r{c}},m_2) \int_n \left( \rho^{\Pi} \right) \\
& = & \sum_{m_1 \in [n]} \sum_{m_2 \in [n] \setminus \ens{m_1}} \sum_{\sigma \in \cS([n-2],[n] \setminus \ens{m_1,m_2})} (-1)^{m_1 + \im(\sigma)_{> m_2}} \sgn(\sigma) \\
& & \left( \prod_{h \in [n-2]} \log_R \alpha(\tau(h),\sigma(h)) \right) \log_R \alpha(\tau^{\r{c}},m_2) \int_n \left( \rho^{\Pi} \right) \\
& = & \sum_{m_1 \in [n]} \sum_{m_2 \in [n] \setminus \ens{m_1}} \sum_{\sigma \in \cS([n-2],[n] \setminus \ens{m_1,m_2})} (-1)^{m_1 + (n - m_2 - \delta_{m_1 > m_2})} \sgn(\sigma) \\
& & \left( \prod_{h \in [n-2]} \log_R \alpha(\tau(h),\sigma(h)) \right) \log_R \alpha(\tau^{\r{c}},m_2) \int_n \left( \rho^{\Pi} \right).
\ee
For a $\sigma \in \cS([n-2],[n])$, we denote by $\sigma^{\r{c},-}$ the minimum of $[n] \setminus \im(\sigma)$, and $\sigma^{\r{c},+}$ the maximum of $[n] \setminus \im(\sigma)$. Summing up the terms in the sum with the same complement $\ens{m_1,m_2}$, we obtain that the right hand side coincides with
\be
& & \sum_{\sigma \in \cS([n-2],[n])} (-1)^{\sigma^{\r{c},-} + \sigma^{\r{c},+} + n} \sgn(\sigma) \left( \prod_{h \in [n-2]} \log_R \alpha(\tau(h),\sigma(h)) \right) \\
& & \left( \log_R \alpha(\tau^{\r{c}},\sigma^{\r{c},+}) - \log_R \alpha(\tau^{\r{c}},\sigma^{\r{c},-} \right) \int_n \left( \rho^{\Pi} \right) \\
& = & \sum_{\sigma \in \cS([n-2],[n])} (-1)^{\sigma^{\r{c},-} + \sigma^{\r{c},+} + n} \sgn(\sigma) \left( \prod_{h \in [n-2]} \log_R \alpha(\tau(h),\sigma(h)) \right) \log_R \rho((\tau^{\r{c}})_{\sigma^{\r{c},-} \to \sigma^{\r{c},+}}) \int_n \left( \rho^{\Pi} \right).
\ee
By Lemma \ref{log combinatorial} and Lemma \ref{log combinatorial 2}, we obtain that the right hand side coincides with
\be
& & \sum_{\sigma \in \cS([n-2],[n])} (-1)^{\sigma^{\r{c},-} + \sigma^{\r{c},+} + n} \sgn(\sigma) \left( \prod_{h \in [n-2]} \log_R \alpha(\tau(h),\sigma(h)) \right) \log_R (\rho^{\Pi})_{\sigma^{\r{c},-} \to \sigma^{\r{c},+}} \int_n \left( \rho^{\Pi} \right) \\
& = & \sum_{\sigma \in \cS([n-2],[n])} (-1)^{\sigma^{\r{c},-} + \sigma^{\r{c},+} + n} \sgn(\sigma) \left( \prod_{h \in [n-2]} \log_R \alpha(\tau(h),\sigma(h)) \right) \int_{n-1} \left( \left( \rho^{\Pi} \right)_{\check{\sigma^{\r{c},+}}} - \left( \rho^{\Pi} \right)_{\check{\sigma^{\r{c},-}}} \right).
\ee
Again, reordering the sum by $\ens{m_1,m_2} = [n] \setminus \im(\sigma)$, we obtain that the right hand side coincides with
\be
& & \sum_{m_1 \in [n]} \sum_{m_2 \in [n] \setminus [m_1]} \sum_{\sigma \in \cS([n-2],[n] \setminus \ens{m_1,m_2})} (-1)^{m_1 + m_2 + n} \sgn(\sigma) \left( \prod_{h \in [n-2]} \log_R \alpha(\tau(h),\sigma(h)) \right) \\
& & \int_{n-1} \left( \left( \rho^{\Pi} \right)_{\check{m_2}} - \left( \rho^{\Pi} \right)_{\check{m_1}} \right) \\
& = & \sum_{m_1 \in [n]} \sum_{m_2 \in [n] \setminus \ens{m_1}} \sum_{\sigma \in \cS([n-2],[n] \setminus \ens{m_1,m_2})} (-1)^{m_1 + m_2 + n} \sgn(\sigma) \left( \prod_{h \in [n-2]} \log_R \alpha(\tau(h),\sigma(h)) \right) \\
& & (-1)^{\delta_{m_1 > m_2}} \int_{n-1} \left( \left( \rho^{\Pi} \right)_{\check{m_2}} \right) \\
& = & \sum_{m_1 \in [n]} \sum_{m_2 \in [n] \setminus \ens{m_1}} \sum_{\sigma \in \cS([n-2],[n] \setminus \ens{m_1,m_2})} (-1)^{m_2 + n - m_1 - \delta_{m_1 > m_2}} \sgn(\sigma) \left( \prod_{h \in [n-2]} \log_R \alpha(\tau(h),\sigma(h)) \right) \\
& & \int_{n-1} \left( \left( \rho^{\Pi} \right)_{\check{m_2}} \right) \\
& = & \sum_{m_1 \in [n]} \sum_{m_2 \in [n] \setminus \ens{m_1}} \sum_{\sigma \in \cS([n-2],[n] \setminus \ens{m_1,m_2})} (-1)^{m_2 + \# ([n] \setminus \ens{m_2})_{< m_1}} \sgn(\sigma) \left( \prod_{h \in [n-2]} \log_R \alpha(\tau(h),\sigma(h)) \right) \\
& & \int_{n-1} \left( \left( \rho^{\Pi} \right)_{\check{m_2}} \right).
\ee
For an $m \in \N$, we denote by $d_{n,m}$ the $m$-th degeneration map. Reordering the double sum by $(m_2,m_1)$ and converting a $\sigma \in \cS([n-2],[n] \setminus \ens{m_1,m_2})$ by $d_{n,m_2} \circ \sigma \in \cS([n-2],[n-1] \setminus \ens{d_{n,m_2}(m_1)})$, we obtain that the right hand side coincides with
\be
& & \sum_{m_2 \in [n]} \sum_{m_1 \in [n] \setminus \ens{m_2}} \sum_{\sigma \in \cS([n-2],[n-1] \setminus \ens{d_{n,m_2}(m_1)})} (-1)^{m_2 + \# ([n] \setminus \ens{m_2})_{< m_1}} \sgn(s_{n,m_2} \circ \sigma) \\
& & \left( \prod_{h \in [n-2]} \log_R \alpha(\tau(h),s_{n,m_2}(\sigma(h))) \right) \int_{n-1} \left( \left( \rho^{\Pi} \right)_{\check{m_2}} \right) \\
& = & \sum_{m_2 \in [n]} \sum_{m_1 \in [n-1]} \sum_{\sigma \in \cS([n-2],[n-1] \setminus \ens{m_1})} (-1)^{m_2 + m_1} \sgn(s_{n,m_2} \circ \sigma) \\
& & \left( \prod_{h \in [n-2]} \log_R \alpha(\tau(h),s_{n,m_2}(\sigma(h))) \right) \int_{n-1} \left( \left( \rho^{\Pi} \right)_{\check{m_2}} \right) \\
& = & \sum_{m \in [n]} \sum_{\sigma \in \cS([n-2],[n-1])} (-1)^{m + \sigma^{\r{c}}} \sgn(\sigma) \left( \prod_{h \in [n-2]} \log_R \alpha(\tau(h),s_{n.],m_2}(\sigma(h))) \right) \int_{n-1} \left( \left( \rho^{\Pi} \right)_{\check{m_2}} \right) \\
& = & \sum_{m \in [n]} (-1)^m \sum_{\sigma \in \cS([n-2],[n-1])} \ol{\sgn}(\sigma) \left( \prod_{h \in [n-2]} \log_R \alpha(\tau(h),s_{n,m_2}(\sigma(h))) \right) \int_{n-1} \left( \left( \rho^{\Pi} \right)_{\check{m_2}} \right) \\
& = & \sum_{m \in [n]} (-1)^m D^{\tau}((\alpha(h,s_{n,m}(j)))_{(h,j) \in [n-1]^2})) \int_{n-1} \left( \left( \rho^{\Pi} \right)_{\check{m_2}} \right).
\ee
Since $(\alpha(h,s_{n,m}(j)))_{(h,j) in [n-1]^2}$ represents $(\rho(h)_{\check{m}})_{h \in [n]}$, we obtain that the right hand side coincides with
\be
\sum_{m \in [n]} (-1)^m \int d^{\tau}((\rho(h)_{\check{m}})_{h \in [n-1]}) = \sum_{m \in [n]} (-1)^m I_{F_{\check{m}}}(d^{\tau}((\rho(h)_{\check{m}})_{h \in [n-1]})) = \sum_{m \in [n]} (-1)^m I_{F_{\check{m}}}(\omega_{\check{m}}).
\ee
\end{proof}

Let $X$ be an adic space over $\Spa(k_{\ad})$, and $M$ a subgroup of $\BdR$. For an integrable $(K/k,n)$-simplex $\gamma$ on $X$, the functional $\int_{\gamma} \colon \rH^0(X,\Omega^n_X) \to \BdR$ is the composite of $k$-linear homomorphisms by the definition of the set-theoretic map $\rH^0(X,\Omega^n_X) \to \Omega^n_{\Delta_{K/k}^n}$ and the continuous map $\Omega^n_{\Delta_{K/k}^n} \to \BdR$, and hence is $k$-linear. By Theorem \ref{Stokes}, the map
\be
\int \colon \rH^0(X,\Omega^n_X) \times C_{K/k,M,n}^{\dagger}(X) \to \rH^0(X,\Omega^n_X) \times (C_{K/k,\Z,n}^{\dagger}(X) \otimes M) \to \BdR
\ee
defined as the $(k,\Z)$-bilinear extension of the integration induces a $(k,\Z)$-bilinear form
\be
\int \colon \rZ^n(\rH^0(X,\Omega_X^{\bullet})) \times \rH_n^{\dagger}(X,K/k,M) \to \BdR,
\ee
where $\rZ^n(\rH^0(X,\Omega_X^{\bullet}))$ denotes the $n$-th cocycle subspace of the cochain complex $\rH^0(X,\Omega_X^{\bullet})$.
Indeed, for any $(\omega,\gamma) \in \rZ^n(\rH^0(X,\Omega_X^{\bullet}) \times C_{K/k,M,n-1}^{\dagger}(X)$, we have
\be
\int_{\partial \gamma} \omega = \int_{\gamma} d \omega = \int_{\gamma} 0 = 0,
\ee
where the $(-1)$-st component of a chain complex is defined as $\ens{0}$. It further induces a $(k,\Z)$-bilinear form
\be
\int \colon \rH^n(\rH^0(X,\Omega_X^{\bullet})) \times \rH_n^{\dagger}(X,K/k,M) \to \BdR,
\ee
Indeed, for any $(\omega,\gamma) \in \rH^0(X,\Omega_X^{n-1}) \times \rH_n^{\dagger}(X,K/k,M)$, we have
\be
\int_{\gamma} d \omega = \int_{\partial \gamma} \omega = \int_{0} \omega = 0,
\ee
where $\Omega_X^{-1}$ is defined as $\ens{0}$. In particular, for any $n$-cycle of $C_{K/k,M}^{\dagger}(X)$, if there exists an $\omega \in \rH^n(\rH^0(X,\Omega_X^{\bullet}))$ such that $\int_{\gamma} \omega \neq 0$, then $\gamma$ does not vanish at $\rH_n^{\dagger}(X,K/k,M)$. This provides a way to show the non-triviality of $\rH_n^{\dagger}(X,K/k,M)$, as desired.

\begin{exm}
\label{Tate curve}
Consider the closed $1$-form $\frac{dT}{T}$ on the Tate curve $\Gm/q^{\Z}$ with $q \in k^{\times}$, and the $1$-cycle $\gamma$ given by the composite of the canonical projection $\Gm \twoheadrightarrow \Gm/q^{\Z}$ and the analytic $1$-simplex on $\Gm$ associated to the global section of $\bA_{K/k}^1$ represented by the $k$-exponential map $\rho$ given by $(1,\ul{q}) \in \Hom_{\Ab}(\Lambda,\cC^{\times})^2 \cong (\cC/k)^{2 \flat \times}$, where $\ul{q} \colon \Lambda \to K^{\times}$ denotes the character given as a system of primitive $p$-power roots of $q$. Then, by the definition of the integration, we have
\be
& & \int_{\gamma} \frac{dT}{T} = \int_{\Delta_{K/k,n}} \gamma^* \frac{dT}{T} = \int_{\Delta_{K/k,n}} \rho^{-1} d \rho = \int d^{s_{1,0}}(\rho^{-1},\rho) \\
& = & 1^{-1} \sum_{\sigma \in \cS([0],[1])} \ol{\sgn}(\sigma) \prod_{h \in [0]} \log_R \rho^{s_{1,0}(h)}_{\sigma^{\r{c}} \to \sigma(h)} \int_1 (\rho^{-1} \rho) \\
& = & \sum_{i=0}^{1} \ol{\sgn}(s_{1,i}) \log_R \rho^{s_{1,0}(0)}_{s_{1,i}^{\r{c}} \to s_{1,i}(h)} \int_1 (1) = \sum_{i=0}^{1} (-1)^i \log_R \rho_{i \to 1-i} \\
& = & (-1)^0 \log_R \ul{q}^{\ell_k} + (-1)^1 \log_R \ul{q}^{-\ell_k} = 2 \ell_k \log_R \ul{q} 
\ee
by $P_{1}(1) = \ens{\ens{0,1}}$. In particular, since $q$ cannot be torsion, Proposition \ref{zero of log} implies $2 \ell_k \log_R \ul{q} \neq 0$. Therefore, the cycle $\gamma$ is not null-homologous in the sense that it is non-trivial as an element of $\rH_1^{\dagger}(\Gm/q^{\Z},K/k,M)$ with $1 \in M$, and especially $\rH_1^{\dagger}(\Gm/q^{\Z},K/k,k)$ is a non-zero representation of $G_k$ admitting a $G_k$-equivariant homomorphism to $k \log \ul{q}$.
\end{exm}

Through the observation in Example \ref{Tate curve}, we expect a good description of the integration using $p$-adic periods in some sense.

\subsection{Appendix}
\label{Appendix}

For an actual application of our integration, we need to check the rigidity of an explicit $(k,n)$-exponential datum. For this purpose, we need to study the topological ring structures of $W(R)$, $W(R)[p^{-1}]$, $W(R)_k[p^{-1}]$, $\BdR^{+}$, and $\BdR$.

\vs
We equip $\kappa$ with the trivial valuation, and $R$ with the map
\be
\n{\bullet}_R \colon R \to [0,\infty), \ \alpha \mapsto \v{\alpha(0)},
\ee
for which $R$ forms a Banach $\kappa$-algebra. We note that we do not use the valuation of $\kappa$ except for the formulation of the notion of a Banach $\kappa$-algebra, and hence the valuation symbol $\v{\bullet}$ always denotes the valuation of $\cC$. The topology of $R$ coincides with the norm topology, and hence with the inverse limit topology of $\varprojlim_{x \mapsto x^p} \cC_{\leq 1}$. This implies that the topology of $W(R)$ coincides with the topology as a topological Abelian group generated by $\set{I_h + p^a W(R)}{(h,a) \in \N^2}$, where $I_h \subset W(R)$ denotes the ideal generated by $\set{[\alpha]}{\alpha \in R \land \n{\alpha} \leq \v{p}^h}$. 

\begin{prp}
\label{topology of W(R) 1/p}
The topology of $W(R)[p^{-1}]$ coincides with the topology as a topological Abelian group generated by $\set{I_{h,a}}{(h,a) \in \N^{\N} \times \N}$, where $I_{h,a} \subset W(R)[p^{-1}]$ denotes the $W(R)$-submodule $(\sum_{b \in \N} p^{-b} I_{h(b)}) + p^a W(R)$, and with the direct limit topology of $\varinjlim_{x \mapsto px} W(R)$.
\end{prp}

\begin{proof}
We denote by $\tau$ the topology of $W(R)[p^{-1}]$ as a topological Abelian group generated by $\set{I_{h,a}}{(h,a) \in \N^{\N} \times \N}$. Then $\tau$ coincides with the direct limit topology of $\varinjlim_{x \mapsto px} W(R)$. Since the multiplication of $p$ is a homeomorphism on any topology of $W(R)[p^{-1}]$ as a topological ring by $p \in W(R)[p^{-1}]^{\times}$, every topology of $W(R)[p^{-1}]$ as a topological $W(R)$-algebra is weaker than or equal to $\tau$. Therefore it suffices to show that $\tau$ is a topology as a topological $W(R)$-algebra. For this purpose, it suffices to show that for any $(x_i)_{i \in 2} \in W(R)[p^{-1}]^2$ and $(h,a) \in \N^{\N} \times \N$, there exists an $(h_i,a_i)_{i \in 2} \in (\N^{\N} \times \N)^2$ such that $(x_0 + I_{h_0,a_0})(x_1 + I_{h_1,a_1}) \subset x_0 x_1 + I_{h,a}$. For each $i \in 2$, denote by $B_i \in \N$ the least natural number satisfying $p^{B_i} x_i \in W(R)$, and put $a_i \coloneqq a+B_{1-i}$. For each $(i,b) \in 2 \times \N$, set $h_i(b) \coloneqq \sum_{b' = 0}^{B_{1-i}+2b} h(b') \in \N$. Then we have
\be
& & \set{x - x_0 x_1}{x \in (x_0 + I_{h_0,a_0})(x_1 + I_{h_1,a_1})} \\
& \subset & \sum_{i \in 2} \sum_{(b_0,b_1) \in \N^2} p^{-(B_i+b_{1-i})} I_{h_{1-i}(b_{1-i})} + p^{-(b_0+b_1)} I_{h_0(b_0)} I_{h_1(b_1)} \\
& & + p^{-B_i+a_{1-i}} W(R) + \sum_{i \in 2} p^{-b_i+a_{1-i}} I_{h_i(b_i)} + p^{a_1+a_2} W(R) \\
& \subset & \sum_{i \in 2} \sum_{(b_0,b_1) \in \N^2} p^{-(B_i+b_{1-i})} I_{h(B_i+b_{1-i})} + p^{-(b_0+b_1)} I_{h(2 \max \ens{b_0,b_1})} \\
& & + p^a W(R) + \sum_{i \in 2} p^{-b_i+a_{1-i}} I_{h(b_i)} + p^a W(R) \\
& \subset & I_{h,a},
\ee
and hence $(x_0 + I_{h_0,a_0})(x_1 + I_{h_1,a_1}) \subset x_0 x_1 + I_{h,a}$.
\end{proof}

By Proposition \ref{topology of W(R) 1/p}, the canonical embedding $W(R) \hookrightarrow W(R)[p^{-1}]$ is a homeomorphism onto the closed image, because the multiplication $W(R) \to W(R)$ by $p$ is a homeomorphism onto the closed image $p W(R)$. Since $W(R)_k[p^{-1}]$ is a free $W(R)[p^{-1}]$-module of finite rank, the topology of $W(R)_k[p^{-1}]$ coincides with the direct product topology with respect to a $W(R)[p^{-1}]$-linear basis, and hence coincides with the topology as a topological Abelian group generated by $\set{k_{\leq 1} I_{h,a}}{(h,a) \in \N^{\N} \times \N}$. By the definition of the inverse limit topology, the topology of $\BdR^{+}$ coincides with the topology as a topological Abelian group generated by $\set{I_{h,a,v}}{(h,a,v) \in \N^{\N} \times \N \times \N}$, where $I_{h,a,v} \subset \BdR^{+}$ denotes the $W(R)_k$-submodule $k_{\leq 1} I_{h,a} + \Fil^v \BdR^{+}$.

\begin{prp}
\label{topology of BdR}
The topology of $\BdR$ coincides with the topology as a topological Abelian group generated by $\set{I_{\xi,h,a,v}}{(h,a,v) \in (\N^{\N})^{\N} \times \N^{\N} \times \N}$ for any $\xi \in \Fil^1 \BdR \setminus \ens{0}$, where $I_{\xi,h,a,v}$ denotes the $W(R)$-submodule $(\sum_{c \in \N} \xi^{-c} k_{\leq 1} I_{h(c),a(c)}) + \Fil^v \BdR$, and with the direct limit topology of $\varinjlim_{x \mapsto \xi x} \BdR^{+}$.
\end{prp}

\begin{proof}
We denote by $\tau$ the topology of $\BdR$ as a topological Abelian group generated by $\set{I_{\xi,h,a,v}}{(h,a,v) \in (\N^{\N})^{\N} \times \N^{\N} \times \N}$. Then $\tau$ coincides with the direct limit topology of $\varinjlim_{x \mapsto \xi x} \BdR^{+}$. Since the multiplication of $\xi$ is a homeomorphism on any topology of $\BdR$ as a topological ring by $\xi \in \BdR^{\times}$, every topology of $\BdR$ as a topological $\BdR^{+}$-algebra is weaker than or equal to $\tau$. Therefore it suffices to show that $\tau$ is a topology as a topological $\BdR^{+}$-algebra. For this purpose, it suffices to show that for any $(x_i)_{i \in 2} \in \BdR^2$ and $(h,a,v) \in (\N^{\N})^{\N} \times \N^{\N} \times \N$, there exists an $(h_i,a_i,v_i)_{i \in 2} \in ((\N^{\N})^{\N} \times \N^{\N} \times \N)^2$ such that $(x_0 + I_{\xi,h_0,a_0,v_0})(x_1 + I_{\xi,h_1,a_1,v_1}) \subset x_0 x_1 + I_{\xi,h,a,v}$. For each $i \in 2$, denote by $V_i \in \N$ the least natural number satisfying $\xi^{V_i} x_i \in \BdR^{+}$, and put $v_i \coloneqq v+V_{1-i}$. For each $(i,c) \in 2 \times \N$, set $h_i(c) \coloneqq \sum_{c' = 0}^{V_{1-i}+2c} h(c') \in \N^{\N}$ and $a_i(c) \coloneqq \sum_{c' = 0}^{V_{1-i}+2c} a(c')$. Then we have
\be
& & \set{x - x_0 x_1}{x \in (x_0 + I_{\xi,h_0,a_0,v_0})(x_1 + I_{\xi,h_1,a_1,v_1})} \\
& \subset & \sum_{i \in 2} \sum_{(c_0,c_1) \in \N^2} \xi^{-(V_i+c_{1-i})} k_{\leq 1} I_{h_{1-i}(c_{1-i}),a_{1-i}(c_{1-i})} + \xi^{-(c_0+c_1)} k_{\leq 1} I_{h_0(c_0),a_0(c_0)} I_{h_1(c_1),a_1(c_1)} \\
& & + \Fil^{-V_i+v_{1-i}} \BdR + \sum_{i \in 2} \xi^{-c_i+v_{1-i}} k_{\leq 1} I_{h_i(c_i),a_i(c_i)} + \Fil^{v_1+v_2} \BdR \\
& \subset & \sum_{i \in 2} \sum_{(c_0,c_1) \in \N^2} \xi^{-(V_i+c_{1-i})} k_{\leq 1} I_{h(V_i+c_{1-i}),a(V_i+c_{1-i})} + \xi^{-(c_0+c_1)} k_{\leq 1} I_{h(2 \max \ens{c_0,c_1})} \\
& & + \Fil^v \BdR + \sum_{i \in 2} \xi^{-c_i+v_{1-i}} k_{\leq 1} I_{h(c_i),a(c_i)} + \Fil^v \BdR \\
& \subset & I_{\xi,h,a,v},
\ee
and hence $(x_0 + I_{\xi,h_0,a_0,v_0})(x_1 + I_{\xi,h_1,a_1,v_1}) \subset x_0 x_1 + I_{\xi,h,a,v}$.
\end{proof}

By Proposition \ref{topology of BdR}, the canonical embedding $\BdR^{+} \hookrightarrow \BdR$ is a homeomorphism onto the closed image, because the multiplication $\BdR^{+} \to \BdR^{+}$ by any $\xi \in \Fil^1 \BdR \setminus \ens{0}$ is a homeomorphism onto the closed image. This implies that $\log$ regarded as a $G_k$-equivariant group homomorphism $1 + \Fil^1 \BdR \to \Fil^1 \BdR$ is continuous with respect to the relative topologies as subsets of $\BdR$.

\vs
We fix a system $\pi_k \in R$ of $p$-power roots of a uniformiser of $k$. The following is an extension of the result that $W(R)_k \cap \ker(\theta^{\BdR})$ is a principal ideal of $W(R)_k$ generated by $[\pi_k] - \pi_k(0)$:

\begin{prp}
\label{kernel of theta}
Let $\xi \in W(R)_k \cap \ker(\theta^{\BdR})$. Let $\alpha \in R$ denote the image of $\xi$ by $\theta^k$, and put $\beta \coloneqq \pi_k^{-1} \alpha \in \Frac(R)$. Then $\beta$ belongs to $R$, and $\xi$ belongs to $[\beta] ([\pi_k] - \pi_k(0)) W(R)_k + \pi_k(0) ([\pi_k] - \pi_k(0)) W(R)_k$.
\end{prp}

\begin{proof}
By the definition of $\alpha$, there exists an $\eta \in W(R)_k$ such that $\xi = [\alpha] + \pi_k(0) \eta$. Therefore, we obtain
\be
0 = \theta^{\BdR}(\xi) = \theta^{\BdR}([\alpha] + \pi_k(0) \eta) = \alpha(0) + \pi_k(0) \theta^{\BdR}(\eta).
\ee
By $\eta \in W(R)_k$, we have $\theta^{\BdR}(\eta) \in \cC_{\leq 1}$. This implies $\alpha(0) = - \pi_k(0) \theta^{\BdR}(\eta) \in \pi_k(0) \cC_{\leq 1}$, and hence $\beta = \pi_k^{-1} \alpha \in R$. Since $W(R)_k$ is $\pi_k(0)$-adically complete, it suffices to show that $\xi$ belongs to $[\beta] ([\pi_k] - \pi_k(0)) W(R)_k + \pi_k(0) \ker(\theta^{\BdR})$. We obtain
\be
\xi = [\alpha] + \pi_k(0) \eta = [\beta] [\pi_k] + \pi_k(0) \eta = [\beta]([\pi_k] - \pi_k(0)) + \pi_k(0)(\eta + [\beta])
\ee
By $\xi \in \ker(\theta^{\BdR})$ and $[\pi_k] - \pi_k(0)) \in \ker(\theta^{\BdR})$, we have $\eta + [\beta] \in \ker(\theta^{\BdR})$. This implies $\xi \in [\beta]([\pi_k] - \pi_k(0)) + \pi_k(0) \ker(\theta^{\BdR})$.
\end{proof}

\begin{prp}
\label{inverse of Witt}
Let $\xi \in W(R)_k \setminus \ker(\theta^{\BdR})$ and $i \in \N$. Put $x \coloneqq \theta^{\BdR}(\xi) \in \cC_{\leq 1} \setminus \ens{0}$ and $\nu \coloneqq \lceil \log_{\v{\pi_k(0)}} \v{x} \rceil \in \N$. Then $\xi$ is invertible in $\BdR^{+}$, and $\xi^{-1}$ belongs to
\be
\sum_{j=0}^{i-1} \pi_k^{- \nu(j+1)}([\pi_k] - \pi_k(0))^j W(R)_k + \Fil^i \BdR.
\ee
\end{prp}

\begin{proof}
By $\xi \in W(R)_k \setminus \ker(\theta^{\BdR})$, we have $x \neq 0$ and $\xi \in \BdR^{+} \setminus \Fil^1 \BdR = (\BdR^{+})^{\times}$. Put $y \coloneqq \pi_k(0)^{\nu} x^{-1} \in \cC$. By $\nu \geq \log_{\v{\pi_k(0)}} \v{x}$, we have
\be
\v{y} = \v{\pi_k(0)}^{\nu} \v{x}^{-1} = \v{\pi_k(0)}^{\nu - \log_{\v{\pi_k(0)}} \v{x}} \leq 1.
\ee
Take a system $\eta \in R$ of $p$-power roots of $y$. By $xy = \pi_k(0)^{\nu}$, there exists a system $\epsilon \in R$ of (not necessarily primitive) $p$-power roots of $1$ such that $\xi \eta - [\epsilon] [\pi_k]^{\nu} \in \ker(\theta^{\BdR})$. Therefore, by Proposition \ref{kernel of theta}, there exists a $\zeta \in W(R)_k$ such that $\xi \eta = [\epsilon][\pi_k]^{\nu} + ([\pi_k] - \pi_k(0)) \zeta$. We have
\be
& & \xi \eta = [\epsilon][\pi_k]^{\nu} + ([\pi_k] - \pi_k(0)) \zeta = \pi_k(0)^{\nu} [\epsilon] \left( 1 + \frac{[\pi_k] - \pi_k(0)}{\pi_k(0)} \right)^{\nu} + ([\pi_k] - \pi_k(0)) \zeta \\
& = & \pi_k(0)^{\nu} [\epsilon] \left( 1 + (\pi_k(0)^{-1} + \pi_k(0)^{- \nu}[\epsilon^{-1}] \zeta)([\pi_k] - \pi_k(0)) + \sum_{v=2}^{\nu} \pi_k(0)^{-v} \binom{\nu}{v} ([\pi_k] - \pi_k(0))^v \right) \\
& = & \pi_k(0)^{\nu} [\epsilon] \left( 1 + \pi_k^{- \nu}([\pi_k] - \pi_k(0)) W(R)_k \right),
\ee
and hence
\be
\xi^{-1} = \eta (\xi \eta)^{-1} & \in & \pi_k(0)^{- \nu} \left( 1 + \sum_{j=1}^{i-1} \pi_k^{- \nu j}([\pi_k] - \pi_k(0))^j W(R)_k + \Fil^i \BdR \right) \\
& = & \sum_{j=0}^{i-1} \pi_k^{- \nu(j+1)}([\pi_k] - \pi_k(0))^j W(R)_k + \Fil^i \BdR
\ee
\end{proof}

We fix a system $\epsilon \in R$ of primitive $p$-power roots of $1 \in k$.

\begin{prp}
\label{valuation of log}
Let $e_k$ denote the maximum of a $q \in \N$ such that $q$ is a power of $p$ and $k$ admits a primitive $q$-th root of $1 \in k$. Let $\alpha \in \Frac(R)^{\times}$, $r \in \N$, and $i \in \N$ with $\alpha(0)^{p^r} \in k^{\times}$ and $\alpha^{\ell_k} \neq 1$. Put
\be
v_p & \coloneqq & \log_{\v{\pi_k(0)}} \v{p} \in \N \\
v_{\epsilon} & \coloneqq & \lceil \log_{\v{\pi_k(0)}} \v{\theta^{\BdR} \left( \frac{\log_R \epsilon}{[\pi_k] - \pi_k(0)} \right)} \rceil \in \N \\
\ell_i & \coloneqq & \lfloor \log_p \max \ens{i,1} \rfloor \in \N \\
C_{p,\epsilon,i} & \coloneqq & \max \ens{1, \ell_i v_p + v_{\epsilon}} \in \N \\
v_{\alpha,r} & \coloneqq & \log_{\pi_k(0)} \v{\alpha(0)^{p^r}} \in \R \\
u_{\alpha,r} & \coloneqq & \log_{\pi_k(0)} \v{\alpha(0)^{p^r \ell_k} - 1} \in \R \cup \ens{\infty} \\
w_{\alpha,r} & \coloneqq & \min \set{\nu \in \N}{\alpha(\nu)^{p^r \ell_k} \neq 1} \in \N \\
N_{\alpha,r} & \coloneqq &
\left\{
\begin{array}{ll}
\v{v_{\alpha,r}} & (\v{\alpha(0)} \neq 0) \\
u_{\alpha,r} & (\v{\alpha(0)} = 0 \land \alpha(0)^{p^r \ell_k} \neq 1) \\
e_k + w_{\alpha,r} & (\alpha(0)^{p^r \ell_k} = 1)
\end{array}
\right.
\in \R \cup \ens{\infty} \\
\rho_{\alpha,\epsilon} & \coloneqq & \frac{\log_R \alpha}{\log_R \epsilon} \in \BdR^{+}
\ee
\bi
\item[(1)] The value $N_{\alpha,r}$ is a positive integer, and the minimum $\nu_{\alpha,r}$ of a $\nu \in \N$ with $\nu > \log_p N_{\alpha,r}$
satisfies $\alpha(\nu)^{p^r} \notin k$.
\item[(2)] The inequality $\v{\theta^{\BdR}(\rho_{\alpha,\epsilon})} \geq \v{p}^{\nu_{\alpha,r} - r - 1}$ holds.
\item[(3)] The value $\rho_{\alpha,\epsilon}$ belongs to $\sum_{j=0}^{i-1} (j+1)^{-1} \pi_k(0)^{- C_{p,\epsilon,i}(j+1)} ([\pi_k] - \pi_k(0))^j W(R)_k + \Fil^i \BdR$.
\item[(4)] The value $\log_R \alpha$ is invertible in $\BdR$, and $(\log_R \alpha)^{-1}$ belongs to
\be
(\log_R \epsilon)^{-1} \left( \sum_{j=0}^{i-1} \pi_k^{-((\nu_{\alpha,r} - r - 1) v_p + C_{p,\epsilon,i}(i+1)) j + (\nu_{\alpha,r} - r - 1) v_p}([\pi_k] - \pi_k(0))^j W(R)_k + \Fil^i \BdR \right)
\ee
\ei
\end{prp}

\begin{proof}
(1) Since $\Frac(R)^{\times}$ is $p$-torsionfree, the condition $\alpha^{\ell_k} \neq 1$ implies $\alpha^{p^r \ell_k} \neq 1$. If $\v{\alpha(0)} \neq 0$, then $v_{\alpha,r}$ is a non-zero integer by $\alpha(0)^{p^r} \in k$, and hence $N_{\alpha,r} = \v{v_{\alpha,r}}$ is a positive integer. If $\v{\alpha(0)} = 0$ and $\alpha(0)^{p^r \ell_k} \neq 1$, then $N_{\alpha,r} = u_{\alpha,r}$ is a positive integer by $\alpha(0)^{p^r \ell_k} \in k_{\leq}^{\times} \setminus \ens{1}$. If $\alpha(0)^{p^r \ell_k} = 1$, then $N_{\alpha,r} = e_k + w_{\alpha,r}$ is a positive integer since $e_k$ is a positive integer and $w_{\alpha,r}$ is a non-negative integer. Therefore, in any case, $N_{\alpha,r}$ is a positive integer. In particular, there exists a $\nu \in \N$ with $p^{\nu} > N_{\alpha,r}$, and hence $\nu_{\alpha,r} \in \N$ makes sense.

\vs
Assume $\alpha(\nu_{\alpha,r})^{p^r} \in k$. Then $\alpha(0)^{p^r} = \alpha(\nu_{\alpha,r})^{p^{r + \nu_{\alpha,r}}}$ admits a $p^{\nu_{\alpha,r}}$-th root $\alpha(0)^{p^r}$ in $k$. First, suppose $\v{\alpha(0)} \neq 0$. We have
\be
p^{- \nu_{\alpha,r}} v_{\alpha,r} = p^{- \nu_{\alpha,r}} \log_{\v{\pi_k(0)}} \v{\alpha(0)^{p^r}} = p^{- \nu_{\alpha,r}} \log_{\v{\pi_k(0)}} \v{\alpha(\nu_{\alpha,r})^{p^{r + \nu_{\alpha,r}}}} = \log_{\v{\pi_k(0)}} \v{\alpha(\nu_{\alpha,r})^{p^r}}
\ee
and the right hand side is a non-zero integer by the assumptions $\alpha(\nu_{\alpha,r})^{p^r} \in k$ and $\v{\alpha(0)} \neq 0$. On the other hand, we have
\be
\v{p^{- \nu_{\alpha,r}} v_{\alpha,r}} < p^{- \log_p N_{\alpha,r}} \v{v_{\alpha,r}} = N_{\alpha,r}^{-1} \v{v_{\alpha,r}} = 1,
\ee
which contradicts that $p^{- \nu_{\alpha,r}} v_{\alpha,r}$ is a non-zero integer.

\vs
Next, suppose $\v{\alpha(0)} = 0$ and $\alpha(0)^{p^r \ell_k} \neq 1$. Then $\log_{\v{\pi_k(0)}} \v{\alpha(\nu_{\alpha,r})^{p^r \ell_k} - 1}$ is a non-zero integer by the assumptions $\alpha(\nu_{\alpha,r})^{p^r} \in k$ and $\alpha(\nu_{\alpha,r})^{p^r \ell_k} \neq 1$. We denote by $f$ the function
\be
(0,1) & \to (0,1) \\
x & \mapsto & \max \ens{\v{p} x,x^p}.
\ee
Then, for any $x \in (0,1)$, we have
\be
f(x) \leq \max \ens{x,\v{p}^{\frac{1}{p-1}}}^p.
\ee
Therefore, we obtain
\be
v & = & u_{\alpha,r} = \log_{\v{\pi_k(0)}} \v{\alpha(0)^{p^r \ell_k} - 1} = \log_{\v{\pi_k(0)}} \v{\alpha(\nu_{\alpha,r})^{p^{r+\nu_{\alpha,r}}\ell_k} - 1} \\
& \geq & \log_{\v{\pi_k(0)}} f^{\nu_{\alpha,r}} \left( \v{\alpha(\nu_{\alpha,r})^{p^r \ell_k} - 1} \right) \geq \log_{\v{\pi_k(0)}} f^{\nu_{\alpha,r}}(\v{\pi_k(0)}) \geq \log_{\v{\pi_k(0)}} \max \ens{\v{\pi_k(0)},\v{p}^{\frac{1}{p-1}}}^{p^{\nu_{\alpha,r}}} \\
& = & p^{\nu_{\alpha,r}} \max \ens{1,\frac{1}{p-1} \log_{\v{\pi_k(0)}} \v{p}} \geq p^{\nu_{\alpha,r}}.
\ee
This contradicts that $\nu_{\alpha,r} > \log_p v$.

\vs
Finally, suppose $\alpha(0)^{p^r \ell_k} = 1$. By $\alpha^{p^r \ell_k} \neq 1$, $\alpha^{p^r \ell_k}$ is a non-trivial system of (possibly non-primitive) $p$-power roots of $1$. By the definition of $w_{\alpha,r}$, $(\alpha(n+w_{\alpha,r})^{p^r \ell_k})_{n \in \N}$ is a system of primitive $p$-power roots of $1$. By the definition of $e_k$, its entry $\alpha(e_k+w_{\alpha,r})^{p^r \ell_k}$ at $e_k$ does not belong to $k$. By $\nu_{\alpha,r} > N_{\alpha,r} = e_k + w_{\alpha,r}$, $\alpha(\nu_{\alpha,r})^{p^r \ell_k}$ does not belong to $k$.

\vs
(2) Put $\xi \coloneqq \alpha(\nu_{\alpha,r})^{p^r}$. By (1), we have $\xi \notin k$, and hence there exists a $\sigma \in G_k$ such that $\sigma(\xi) \neq \xi$. Since $\sigma(\epsilon)$ is a system of primitive $p$-power roots of $1$, there exists a $u_{\sigma} \in \Zp^{\times}$ such that $\sigma(\epsilon) = \epsilon^{u_{\sigma}}$. Put $\epsilon_{\alpha,r} \coloneqq \alpha^{-p^r} \sigma(\alpha^{p^r}) \in \Frac(R)$. By $\alpha^{p^r}(0) \in k^{\times}$, $\epsilon_{\alpha,r}(\nu_{\alpha,r}) = \xi^{-1} \sigma(\xi)$ is a non-trivial (but possibly non-primitive) $p^{\nu_{\alpha,r}}$-th root of $\alpha^{-p^r}(0) \sigma(\alpha^{p^r}(0)) = 1$. Therefore, there exists an $m_{\alpha,r} \in \Zp \setminus p^{\nu_{\alpha,r}} \Zp$ such that $\epsilon_{\alpha,r} = \epsilon^{m_{\alpha,r}}$. We obtain
\be
\frac{\log_R \epsilon_{\alpha,r}}{\log_R \epsilon} & = & \frac{\log_R \epsilon^m}{\log_R \epsilon} = m_{\alpha,r},
\ee
and hence
\be
\v{\theta^{\BdR} \left( \frac{\log_R \epsilon_{\alpha,r}}{\log_R \epsilon} \right)} & = & \v{\theta^{\BdR}(m_{\alpha,r})} = \v{m_{\alpha,r}} > \v{p}^{\nu_{\alpha,r}-1}.
\ee
Since $\log_R$ is a $G_k$-equivariant group homomorphisms, we have
\be
\log_R \epsilon_{\alpha,r} = p^r(\sigma(\log_R \alpha) - \log_R \alpha)
\ee
and
\be
\frac{\log_R \epsilon_{\alpha,r}}{\log_R \epsilon} & = & p^r \left( \frac{\sigma(\log_R \alpha)}{\log_R \epsilon} - \frac{\log_R \alpha}{\log_R \epsilon} \right) = p^r \left( \frac{\sigma(\log_R \epsilon)}{\log_R \epsilon} \times \frac{\sigma(\log_R \alpha)}{\sigma(\log_R \epsilon)} - \frac{\log_R \alpha}{\log_R \epsilon} \right) \\
& = & p^r \left( \frac{\log_R \sigma(\epsilon)}{\log_R \epsilon} \sigma( \rho_{\alpha,\epsilon}) - \rho_{\alpha,\epsilon} \right) = p^r \left( \frac{\log_R \epsilon^{u_{\sigma}}}{\log_R \epsilon} \sigma(\rho_{\alpha,\epsilon}) - \rho_{\alpha,\epsilon} \right) \\
& = & p^r \left( u_{\sigma} \sigma(\rho_{\alpha,\epsilon}) - \rho_{\alpha,\epsilon} \right).
\ee
Since $\theta^{\BdR}$ is a $G_k$-equivariant ring homomorphism and the action of $G_k$ to $\cC$ preserves the valuation, if $\v{\theta^{\BdR}(\rho_{\alpha,\epsilon})} < \v{p}^{\nu_{\alpha,r} - r - 1}$, then we would have
\be
\v{\theta^{\BdR} \left( \frac{\log_R \epsilon_{\alpha,r}}{\log_R \epsilon} \right)} = \v{p}^r \v{\sigma(\theta^{\BdR}(\rho_{\alpha,\epsilon})) - \theta^{\BdR}(\rho_{\alpha,\epsilon})} \leq \v{p}^r \v{\theta^{\BdR}(\rho_{\alpha,\epsilon})} < \v{p}^{\nu_{\alpha,r} - 1},
\ee
which would contradict
\be
\v{\theta^{\BdR} \left( \frac{\log_R \epsilon_{\alpha,r}}{\log_R \epsilon} \right)} = \v{m_{\alpha,r}} \geq \v{p}^{\nu_{\alpha,r} - 1}.
\ee
Therefore, we obtain $\v{\theta^{\BdR}(\rho_{\alpha,\epsilon})} \geq \v{p}^{\nu_{\alpha,r} - r - 1}$.

\vs
(3), (4) Put $\beta \coloneqq \pi_k^{-v_{\alpha,r}} \alpha^{p^r} \in \Frac(R)$. By the definition of $v_{\alpha,r}$, $\beta$ is a system of $p$-power roots of elements of $\cC_{\leq 1}^{\times}$, and hence belongs to $R$. By Proposition \ref{kernel of theta}, we have
\be
& & \log_R \alpha = \log_R \pi_k^{v_{\alpha,r}} \beta = v_{\alpha,r} \log_R \pi_k + \log_R \beta \\
& \in & \sum_{j=1}^{i} \frac{(-1)^j}{j} (v_{\alpha,r} (\cN_R(\pi_k)-1)^j + (\cN_R(\beta)-1)^j) + \Fil^{i+1} \BdR \\
& = & \sum_{j=1}^{i} \frac{(-1)^j}{j} (\pi_k(0)^{-j}v_{\alpha,r}([\pi_k] - \pi_k(0))^j + \beta(0)^{-j}([\beta] - \beta(0))^j)) + \Fil^{i+1} \BdR \\
& \subset & ([\pi_k] - \pi_k(0)) \left( \sum_{j=0}^{i-1} (j+1)^{-1} \pi_k(0)^{-(j+1)} ([\pi_k] - \pi_k(0))^j W(R)_k + \Fil^i \BdR \right).
\ee
By Proposition \ref{kernel of theta}, we have
\be
p^{\ell_i} \frac{\log_R \epsilon}{[\pi_k] - \pi_k(0)} \in (W(R)_k \setminus \ker(\theta^{\BdR}))+ \Fil^i \BdR.
\ee
Therefore, by (2) and Proposition \ref{inverse of Witt} applied to the left hand side modulo $\Fil^i \BdR$, we have
\be
(\log_R \epsilon)^{-1} & \in & \frac{p^{\ell_i}}{[\pi_k] - \pi_k(0)} \left( \sum_{j=0}^{i-1} \pi_k^{-(\ell_i v_p + v_{\epsilon})(j+1)}([\pi_k] - \pi_k(0))^j W(R)_k + \Fil^i \BdR \right) \\
& = & \frac{1}{[\pi_k] - \pi_k(0)} \left( \sum_{j=0}^{i-1} \pi_k^{-C_{p,\epsilon,i}(j+1)}([\pi_k] - \pi_k(0))^j W(R)_k + \Fil^i \BdR \right).
\ee
We obtain
\be
\rho_{\alpha,\epsilon} & \in & ([\pi_k] - \pi_k(0)) \left( \sum_{j=0}^{i-1} (j+1)^{-1} \pi_k(0)^{-(j+1)} ([\pi_k] - \pi_k(0))^j W(R)_k + \Fil^i \BdR \right) \frac{1}{[\pi_k] - \pi_k(0)}\\
& &  \left( \sum_{j=0}^{i-1} \pi_k^{-C_{p,\epsilon,i}(j+1)}([\pi_k] - \pi_k(0))^j W(R)_k + \Fil^i \BdR \right) \\
& = & \left( \sum_{j=0}^{i-1} (j+1)^{-1} \pi_k(0)^{-(j+1)} ([\pi_k] - \pi_k(0))^j W(R)_k + \Fil^i \BdR \right) \\
& & \left( \sum_{j=0}^{i-1} \pi_k^{-C_{p,\epsilon,i}(j+1)}([\pi_k] - \pi_k(0))^j W(R)_k + \Fil^i \BdR \right) \\
& \subset & \sum_{j=0}^{i-1} (j+1)^{-1} \pi_k(0)^{-C_{p,\epsilon,i}(j+2)} ([\pi_k] - \pi_k(0))^j W(R)_k + \Fil^i \BdR.
\ee
This implies
\be
\pi_k(0)^{C_{p,\epsilon,i}(i+1)} \rho_{\alpha,\epsilon} \in W(R)_k + \Fil^i \BdR,
\ee
and we also have
\be
\v{\theta^{\BdR}(\pi_k(0)^{C_{p,\epsilon,i}(i+1)} \rho_{\alpha,\epsilon})} = \v{\pi_k(0)}^{C_{p,\epsilon,i}(i+1)} \v{\theta^{\BdR}(\rho_{\alpha,\epsilon}))} \geq \v{\pi_k(0)}^{C_{p,\epsilon,i}(i+1)} \v{p}^{\nu_{\alpha,r} - r - 1} = \v{\pi_k(0)}^{(\nu_{\alpha,r} - r - 1) v_p + C_{p,\epsilon,i}(i+1)}.
\ee
Especially, the inequality implies $\pi_k(0)^{C_{p,\epsilon,i}(i+1)} \rho_{\alpha,\epsilon} \notin \ker(\theta^{\BdR})$. By (2) and Proposition \ref{inverse of Witt} applied to $\pi_k(0)^{C_{p,\epsilon,i}(i+1)} \rho_{\alpha,\epsilon}$ modulo $\Fil^i \BdR$, $\rho_{\alpha,\epsilon}$ is invertible in $\BdR$ and we obtain
\be
(\rho_{\alpha,\epsilon})^{-1} & \in & \pi_k(0)^{C_{p,\epsilon,i}(i+1)} \left( \sum_{j=0}^{i-1} \pi_k^{-((\nu_{\alpha,r} - r - 1) v_p + C_{p,\epsilon,i}(i+1))(j+1)}([\pi_k] - \pi_k(0))^j W(R)_k + \Fil^i \BdR \right) \\
& \subset & \sum_{j=0}^{i-1} \pi_k^{-((\nu_{\alpha,r} - r - 1) v_p + C_{p,\epsilon,i}(i+1)) j + (\nu_{\alpha,r} - r - 1) v_p}([\pi_k] - \pi_k(0))^j W(R)_k + \Fil^i \BdR.
\ee
This implies that $\log_R \alpha$ is invertible in $\BdR$ and we obtain
\be
(\log_R \alpha)^{-1} \in (\log_R \epsilon)^{-1} \left( \sum_{j=0}^{i-1} \pi_k^{-((\nu_{\alpha,r} - r - 1) v_p + C_{p,\epsilon,i}(i+1)) j + (\nu_{\alpha,r} - r - 1) v_p}([\pi_k] - \pi_k(0))^j W(R)_k + \Fil^i \BdR \right).
\ee
\end{proof}

We recall that we skipped the proof of Proposition \ref{zero of log} because it is a very weak consequence of Proposition \ref{valuation of log}.

\begin{proof}[Proof of Proposition \ref{zero of log}.]
The condition (1) is equivalent to the condition (2) and the condition (3) implies the condition (1), because $\log_R$ is a group homomorphism and $\ell_k$ is invertible in $\BdR$. Suppose that $(\alpha,c)$ does not satisfy the condition (3). If $c = 0$, then the condition (1) holds by the definition of the power. If $c \neq 0$, then we have $\alpha^{\ell_k} \neq 1$, and hence $(\log_R \alpha)^c \neq 0$ by Proposition \ref{valuation of log} (4).
\end{proof}

\begin{prp}
\label{estimation of inverse of log}
Introduce the following conventions:
\bi
\item[(1)] (Constant depending only on $p$) Let $C$ denote the minimum of a $\nu \in \N$ with $x \log_p(1+x^{-1}) \leq \nu$ for any $x \in (0,1]$, which exists by the continuity of $\log_p$ and $\lim_{x \to +0} x \log_p(1+x^{-1}) = 0$.
\item[(2)] (Notions depending also on $k$) Let $v_{p,k}$ denote the $\pi_k(0)$-adic additive valuation $\log_{\v{\pi_k(0)}} \v{p} \in \Z$ of $p$. Let $\epsilon \in (\cC/k)^{\times}$ be a system of primitive $p$-power roots of $1$. Put
\be
b \coloneqq \theta^{\BdR} \left( \frac{[\pi_k] - \pi_k(0)}{\log_k \epsilon} \right) \in \cC,
\ee
and let $e$ denote the least integer greater than or equal to $- \log_{\v{p}} \v{b}$.
\item[(3)] (Variables and constants depending on them) Let $(\alpha,i) \in R \times \N$ with $\alpha(0) \in k_{\leq 1} \setminus \ens{0}$ and $\alpha^{\ell_k} \neq 1$. Let $M_{\alpha,k}$ denote the minimum of a $\nu \in \N$ with $\alpha(\nu) \notin k$, which exists by Proposition \ref{valuation of log} (1), $N_{\alpha,k}$ the minimum of a $\nu \in \N$ with $(e + M_{\alpha,k} + 1) v_{p,k} \leq \nu$, and $v_{\alpha,k}$ the $\pi_k(0)$-adic additive valuation $\log_{\v{\pi_k(0)}} \v{\alpha(0)} \in \Z$ of $\alpha(0)$.
\ei
Then the inverse of $[\alpha] - \alpha(0)$ in $\BdR$ belongs to
\be
\frac{1}{\pi_k(0)^{N_{\alpha,k}}([\pi_k] - \pi_k(0))} W(R)_k + \sum_{d=0}^{i-1} \frac{([\pi_k] - \pi_k(0))^d}{\pi_k(0)^{(d+1)N_{\alpha,k}}} W(R) + \Fil^i \BdR^{+},
\ee
and the inverse of $\log_R \alpha$ in $\BdR$ belongs to
\be
\frac{1}{\pi_k(0)^{\max \ens{Ci v_{p,k} + (i-1)v_{\alpha,k}, i N_{\alpha,k} - v_{\alpha,k}}}} \sum_{d=-1}^{i-1} ([\pi_k] - \pi_k(0))^d W(R)_k + \Fil^i \BdR^{+}.
\ee
\end{prp}

\begin{proof}
We denote by $\ol{\alpha(0)} \in \kappa$ the image of $\alpha(0)$ by the reduction $k_{\leq 1} \to \kappa$, and regard it as an element of $\cC_{\leq}/p \cC_{\leq}$ through the canonical embedding $W(\kappa) \hookrightarrow k_{\leq 1} \subset \cC_{\leq 1}$. We denote by $\alpha_0 \in R$ the system given as $((\ol{\alpha(0)})^{(p^{-1}q_k)^{\nu}})_{\nu \in \N}$. Then the image of $[\alpha] - \alpha(0)$ by $\theta^{\BdR}$ is $\alpha - \alpha_0$. We denote by $u \in k_{\leq 1}$ the image of $\ol{\alpha(0)}$ by the \Teichmuller embedding $\kappa \hookrightarrow W(\kappa) \hookrightarrow k_{\leq 1}$. Put
\be
a \coloneqq \theta^{\BdR} \left( \frac{[\alpha] - \alpha(0)}{\log_k \epsilon} \right) \in \cC
\ee
Since $[\pi_k] - \pi_k(0)$ generates $\Fil^1 \BdR$ by Proposition \ref{kernel of theta}, we have $b \neq 0$. By Proposition \ref{kernel of theta} applied to $[\alpha] - \alpha(0)$, $\pi_k^{-1}(\alpha - \alpha_0) \in \Frac(R)$ belongs $R$ and there exists an $\eta \in W(R)_k$ such that
\be
[\alpha] - \alpha(0) = [\pi_k^{-1}(\alpha - \alpha_0)]([\pi_k] - \pi_k(0)) + \pi_k(0)([\pi_k] - \pi_k(0)) \eta = ([\pi_k^{-1}(\alpha - \alpha_0)] + \pi_k(0) \eta)([\pi_k] - \pi_k(0)).
\ee
By the definition of $a$, we have
\be
a = \theta^{\BdR} \left( \frac{([\pi_k^{-1}(\alpha - \alpha_0)] + \pi_k(0) \eta)([\pi_k] - \pi_k(0))}{\log_k \epsilon} \right) = (\pi_k^{-1}(0)(\alpha(0) - u) + \pi_k(0) \theta^{\BdR}(\eta))b.
\ee
By Proposition \ref{valuation of log} (3) applied to $(\alpha,M_{\alpha,k})$, we have
\be
a & = & \theta^{\BdR} \left( \frac{[\alpha] - \alpha(0)}{\log_k \epsilon} \right) = - \theta^{\BdR} \left( \frac{1 - \frac{[\alpha]}{\alpha(0)}}{\log_k \epsilon} \right) = - \theta^{\BdR} \left( \frac{\log_R \alpha}{\log_k \epsilon} \right) \\
& \in & \cC \setminus p^{M_{\alpha,k} - 0 + 1} \cC_{\leq 1} = \cC \setminus p^{N_{\alpha,k}+1} \cC_{\leq 1}
\ee
Therefore, we obtain
\be
\frac{[\alpha] - \alpha(0)}{[\pi_k] - \pi_k(0)} \in (\theta^{\BdR})^{-1}(ab^{-1}) \cap W(R)_k \subset (\theta^{\BdR})^{-1}(\cC \setminus b^{-1} p^{M_{\alpha,k}+1} \cC_{\leq 1}) \cap W(R)_k.
\ee
By $[\alpha] - \alpha(0) \in ([\pi_k] - \pi_k(0)) W(R)_k \subset (\theta^{\BdR})^{-1}(\cC_{\leq 1})$ and $p^{e + M_{\alpha,k} + 1} \cC_{\leq 1} \subset b^{-1} p^{M_{\alpha,k}+1} \cC_{\leq 1}$,  we obtain
\be
ab^{-1} = \theta^{\BdR} \left( \frac{[\alpha] - \alpha(0)}{[\pi_k] - \pi_k(0)} \right) \in \cC_{\leq 1} \setminus p^{e + M_{\alpha,k} + 1} \cC_{\leq 1},
\ee
and hence
\be
\frac{[\alpha] - \alpha(0)}{[\pi_k] - \pi_k(0)} & \in & (\theta^{\BdR})^{-1}(\cC_{\leq 1} \setminus p^{e + M_{\alpha,k} + 1} \cC_{\leq 1}) \cap W(R)_k.
\ee
We fix a system $\beta \in R$ of primitive $p$-power roots of $ab^{-1} \in \cC_{\leq 1} \setminus p^{e + M_{\alpha,k} + 1} \cC_{\leq 1}$. By $\v{ab^{-1}} > \v{p}^{e + M_{\alpha,k} + 1} \geq \v{\pi_k(0)}^{N_{\alpha,k}}$, there exists a $\gamma \in R$ such that $[\pi_k]^{N_{\alpha,k}} = \beta \gamma$. We have
\be
\beta \gamma = [\pi_k]^{N_{\alpha,k}} = (\pi_k(0) + ([\pi_k] - \pi_k(0)))^{N_{\alpha,k}} \in \pi_k(0)^{N_{\alpha,k}} + ([\pi_k] - \pi_k(0)) W(R)_k.
\ee
This implies
\be
& & \frac{[\alpha] - \alpha(0)}{[\pi_k] - \pi_k(0)} \gamma = \beta \gamma + \left( \frac{[\alpha] - \alpha(0)}{[\pi_k] - \pi_k(0)} - \beta \right) \gamma \\
& \in & \pi_k(0)^{N_{\alpha,k}} + ([\pi_k] - \pi_k(0)) W(R)_k + (\ker(\theta^{\BdR}) \cap W(R)_k) = \pi_k(0)^{N_{\alpha,k}} + ([\pi_k] - \pi_k(0)) W(R)_k,
\ee
and hence
\be
\frac{1}{[\alpha] - \alpha(0)} & = & \frac{\gamma}{\pi_k(0)^{N_{\alpha,k}} [\pi_k] - \pi_k(0)} \left( \pi_k(0)^{-N_{\alpha,k}} \frac{[\alpha] - \alpha(0)}{[\pi_k] - \pi_k(0)} \gamma \right)^{-1} \\
& \in & \left( \frac{\gamma}{\pi_k(0)^{N}([\pi_k] - \pi_k(0))} \sum_{d=0}^{i} \left( \frac{[\pi_k] - \pi_k(0)}{\pi_k(0)^{N_{\alpha,k}}} \right)^d W(R) \right) + \Fil^i \BdR^{+} \\
& \in & \frac{1}{\pi_k(0)^{N_{\alpha,k}}([\pi_k] - \pi_k(0))} W(R)_k + \sum_{d=0}^{i-1} \frac{([\pi_k] - \pi_k(0))^d}{\pi_k(0)^{(d+1)N_{\alpha,k}}} W(R) + \Fil^i \BdR^{+}.
\ee
By the definition of the logarithm, we have
\be
\log_R \alpha = \sum_{d=1}^{\infty} \frac{(-1)^d}{d} \left( 1 - \frac{[\alpha]}{\alpha(0)} \right)^d = \frac{1}{\alpha(0)}([\alpha] - \alpha(0)) \sum_{d=0}^{\infty} \frac{1}{(d+1) \alpha(0)^d}([\alpha] - \alpha(0))^d.
\ee
For each $(m,l) \in \N^2$ with $m \leq l$, put $P_m(l) \coloneqq \set{d \in \N^m}{\sum_{j=0}^{m-1} d(j) = l}$. We obtain
\be
\frac{1}{\log_R \alpha} & = & \frac{\alpha(0)}{[\alpha] - \alpha(0)} \left( 1 + \sum_{d=1}^{\infty} \frac{1}{(d+1) \alpha(0)^d}([\alpha] - \alpha(0))^d \right)^{-1} \\
& = & \frac{\alpha(0)}{[\alpha] - \alpha(0)} \sum_{m=0}^{\infty} \left( - \sum_{d=1}^{\infty} \frac{1}{(d+1) \alpha(0)^d}([\alpha] - \alpha(0))^d \right)^m \\
& = & \frac{\alpha(0)}{[\alpha] - \alpha(0)} \sum_{l=0}^{\infty} \left( \sum_{m=0}^{l} \sum_{d \in P_m(l)} \frac{(-1)^m}{\prod_{j=0}^{m-1} (d(j) + 1)} \right) \frac{1}{\alpha(0)^l}([\alpha] - \alpha(0))^l.
\ee
For any $(m,l) \in \N^2$ with $1 \leq m \leq l$ and $d \in P_m(l)$, we have
\be
\prod_{j=0}^{m-1} (d(j) + 1) \leq \left( \frac{1}{m} \sum_{j=0}^{m-1} (d(j) + 1) \right)^m = \left( 1 + \frac{l}{m} \right)^m = \left( \left( 1 + \frac{l}{m} \right)^{\frac{m}{l}} \right)^l \leq p^{M l}.
\ee
We obtain
\be
& & \frac{1}{\log_R \alpha} = \frac{\alpha(0)}{[\alpha] - \alpha(0)} \sum_{l=0}^{\infty} \left( \sum_{m=0}^{l} \sum_{d \in P_m(l)} \frac{(-1)^m}{\prod_{j=0}^{m-1} (d(j) + 1)} \right) \frac{1}{\alpha(0)^l}([\alpha] - \alpha(0))^l \\
& \in & \frac{\alpha(0)}{[\alpha] - \alpha(0)} + \sum_{d=1}^{i} \frac{1}{p^{M d} \alpha(0)^{d-1}}([\alpha] - \alpha(0))^{d-1} W(R)_k + \Fil^i \BdR^{+} \\
& \subset & \frac{\alpha(0)}{\pi_k(0)^{N_{\alpha,k}}([\pi_k] - \pi_k(0))} W(R)_k \\
& & + \sum_{d=0}^{i-1} \left( \frac{1}{p^{M(d+1)} \alpha(0)^d} ([\pi_k] - \pi_k(0))^d W(R)_k + \frac{\alpha(0)}{\pi_k(0)^{(d+1)N_{\alpha,k}}} ([\pi_k] - \pi_k(0))^d W(R)_k \right) + \Fil^i \BdR^{+}. \\
& \subset & \frac{1}{\pi_k(0)^{\max \ens{M(d+1) v_{p,k} +  dv_{\alpha,k}, i N_{\alpha,k} - v_{\alpha,k}}}} \sum_{d=-1}^{i-1} ([\pi_k] - \pi_k(0))^d W(R)_k + \Fil^i \BdR^{+}.
\ee
\end{proof}

\begin{prp}
\label{estimation of inverse of log 2}
Following the conventions (1) and (2) in Proposition \ref{estimation of inverse of log}, introduce the following convention:
\bi
\item[(3)'] (Variables and constants depending on them) Let $(\alpha,r,i) \in \Frac(R) \times \N \times \N$ with $\alpha(0)^{p^r} \in k^{\times}$ and $\alpha^{\ell_k} \neq 1$. Let $M_{\alpha,k,r}$ denote the minimum of a $\nu \in \N$ with $\alpha(\nu)^{p^r} \notin k$, which exists by Proposition \ref{valuation of log} (1), $N_{\alpha,k,r}$ the minimum of a $\nu \in \N$ with $(e + M_{\alpha,k,r} + 1) v_{p,k} \leq \nu$, and $v_{\alpha,k,r}$ the rational number $\log_{\v{\pi_k(0)}} \v{\alpha(0)} \in p^{-r} \Z$.
\ei
Then the inverse of $\log_R \alpha$ in $\BdR$ belongs to
\be
\frac{p^r}{\pi_k(0)^{\max \ens{Ci v_{p,k} + (i-1) p^r \v{v_{\alpha,k,r}}, i N_{\alpha,k,r} - p^r \v{v_{\alpha,k,r}}}}} \sum_{d=-1}^{i-1} ([\pi_k] - \pi_k(0))^d W(R)_k + \Fil^i \BdR^{+}.
\ee
\end{prp}

\begin{proof}
The assertion immediately follows from Proposition \ref{estimation of inverse of log}, $\log_R \alpha = p^{-r} \log_R \alpha^{p^r}$, and $\log_R \alpha = - \log_R \alpha^{-1}$, since $\alpha^{\ell_k} \neq 1$ implies $\alpha^{p^r \ell_k} \neq 1$ as we have already shown in the beginning of the proof of Proposition \ref{valuation of log} (1).
\end{proof}

%% file: References.tex
\vspace{0.3in}
\noindent {\Large \bf Acknowledgements}
\vs

\noindent
This work is partially based on my master thesis written twelve years ago (2012), which was not related to perfectoid theory because it was written before the publication of \cite{Sch12}. Therefore, I first acknowledge people who helped me when I was writing the master thesis.

\vs
I am thankful to my supervisor Takeshi Tsuji for his great help in seminars, lectures, and so on. I improved the master thesis thanks to his various feedback. I would like to thank my friends. In particular, I asked many times for help to Kazuki Tokimoto in the field of algebra, to Yusuke Isono in the field of analysis, to Fumihiko Sanda in the field of geometry, and to Yuichiro Tanaka in the total point of view. Several works related to the master thesis on Tate's acyclicity and $p$-adic representation, through which I came across the idea to apply perfectoid theory in the formulation of analytic singular homology, was supported by JSPS KAKENHI Grant Number 13J08878 (2013/04--2015/03).

\vs
Also, I thank Peter Scholze for the investigation of perfectoid theory. When I wrote the master thesis, I have no idea on how to prove Tate's acyclicity for candidates of cosimplicial objects. Therefore, I needed to try a formulation an analytic singular homology in a way abandoning the sheaf condition. Thanks to perfectoid theory, I obtain a better formulation of an analytic singular homology, as is written in this paper.

\vs
This work was supported by JSPS KAKENHI Grant Number 21K13763 (2021/04--2024/03). I am also thankful to my family for their deep affection.